\documentclass[review]{elsarticle}

\usepackage{lineno,hyperref}
\usepackage{graphicx}
\usepackage[mathscr]{eucal}
\usepackage{multirow}
\usepackage{amsmath,amssymb,amsthm,caption}
\usepackage{enumitem}
\usepackage{color}
\usepackage[labelformat=simple]{subcaption}

\def\be{\begin{equation}\label}
\def\ee{\end{equation}}

\def\barr{\begin{array}}
\def\earr{\end{array}}

\def\ber{\begin{eqnarray}}
\def\eer{\end{eqnarray}}

\def\bers{\begin{eqnarray*}}
\def\eers{\end{eqnarray*}}

\def\bpf{\begin{pf}}
\def\epf{\end{pf}}

\modulolinenumbers[5]

\journal{Journal of Economic Behavior \& Organization}

\begin{document}
\rmfamily

\begin{frontmatter}
	
	\title{A continuum of path-dependent equilibrium solutions induced by sticky 	expectations}
	
	
	\author[IMCAS]{Pavel Krej\v{c}\'\i}
	\author[RMU]{Eyram Kwame}
	\author[GMU]{ Harbir Lamba}
	\author[UTD]{ Dmitrii  Rachinskii}
	
	
	\address[IMCAS]{Institute of Mathematics of the Czech Academy of Sciences}
	\address[RMU]{Regional Maritime University, Accra, Ghana}
	\address[GMU]{Department of Mathematical Sciences,	George Mason University, VA, USA}
	\address[UTD]{Department of Mathematical Sciences, University of Texas at Dallas, TX, USA}
				
	\begin{abstract}
		We analyze a simple macroeconomic model where
		rational inflation expectations is replaced by a boundedly
		rational, and genuinely sticky, response
		to changes in the actual inflation rate. The stickiness is
		introduced in a novel way using a 
		mathematical operator that is
		amenable to rigorous analysis. We prove that, when exogenous
		noise is absent from the system, the unique equilibrium of the
		rational expectations model is replaced by an entire line segment
		of possible equilibria with the one chosen depending, in a
		deterministic way, upon the
		previous states of the system. The agents are sufficiently
		far-removed from the rational expectations paradigm that problems of
		indeterminacy do not arise.
		
		The response to exogenous noise is far more subtle than in a
		unique equilibrium model. After sufficiently small shocks the system
		will indeed revert to the same equilibrium but larger ones will move
		the system to a different one (at the same model
		parameters).  The path to this new equilibrium may be very long
		with a highly unpredictable endpoint. At certain model parameters
		exogenously-triggered runaway
		inflation can occur.
		
		Finally, we analyze a variant model in which the same form of sticky
		response is introduced into the interest rate rule instead.
		
	\end{abstract}
	
	
\end{frontmatter}



\section{Introduction}\label{intro}

Modern macroeconomics has been dominated by a
modeling framework in which the economy is assumed always to be
at (or rapidly moving back towards) a unique and stable
equilibrium. This has had profound implications both for the way in which
the modelers perceive real-world events and their policy prescriptions
for dealing with them.

The critiquing of equilibrium models has a long 
history which we shall not attempt to detail here. But many
antagonists, see for example
\cite{robinson1974,n72,s97a,colander2008beyond}, have  eloquently
pointed out profound issues concerning the assumed equilibrating
processes and the ways in which the `aggregation problem' was being
solved. 
In this paper we will focus upon one specific pillar of the
equilibrium approach which is the assumption of Rational Expectations
introduced by Muth in 1961 \cite{m61}. This posits that not only are
individuals perfectly rational, optimizing, far-sighted and
independent of each other but that their expectations about future
uncertainties are in agreement with the model itself. 

Our mathematical analysis, and the supporting numerics, rigorously
show that, when rational expectations about future inflation are
replaced by an aggregated `sticky' expectation, a simple macroeconomic
model changes from a unique equilibrium system to one with an entire
continuum of path-dependent equilibria. The form of stickiness that we
use is, to our knowledge, new in a macroeconomic setting and differs
from, for example, the stickiness of the Calvo pricing model
\cite{calvo1983staggered} where hypothetical agents are only allowed
to adjust (to the correct price) at a fixed rate.

The way in which we incorporate stickiness into the model will be
justified and described more fully below but, briefly,
our sticky variables can only be in one of two modes. They are either
currently `stuck' at some value or they are being `dragged' along by
some other (related) variable because the maximum allowable difference between
them has been reached. Each of these modes (which we shall also refer
to as the `inner' and 'outer' modes respectively) can be  
analyzed separately as linear systems using standard stability
techniques. However the full `hybrid' system is nonlinear and displays
far richer dynamics in the presence of exogenous noise and shocks.

It must be emphasized right away that our modelling approach and
analytical tools are not restricted to inflation expectations or even
to macroeconomics. The form of stickiness described above belongs to a
class of operators that have well-understood and very desirable
properties. These have already been used to develop
non-equilibrium asset-pricing models \cite{nsf5} that have (almost-)
analytic solutions.


Here we are able to prove the existence of an entire line interval of
feasible equilibrium points, examine their stability, and identify
some important consequences of path dependence regarding the effects of
exogenous shocks and policy changes upon the state of the
system. Furthermore, these changes are realistic in that they both
correspond closely to observed, but potentially puzzling, economic
situations and are robust enough to be observed numerically in more
sophisticated variants of the model.

The level of mathematical knowledge required to follow most of the
arguments is not much more than is needed to examine the existence and
stability of equilibria in more traditional, fully linear,
models. Another useful aspect of this simple model is that the
stickiness can be smoothly `dialed back' to zero and the unique
equilibrium case is recovered.  Or, to put it another way, we can
rigorously show that a plausible, boundedly rational {\em yet fully
analyzable}, change to a fully rational model significantly
alters the qualitative behaviour of the system in recognizable ways.

Before introducing the model and starting the mathematical analysis,
it is worth stepping back to consider the effects of stickiness and
friction in physical rather than economic systems. This helps develop
our intuition about the nature of equilibria in such systems but the
comparison also offers a high-level explanation of the failure of
mainstream economics to foresee the recent economic crisis and the difficulties
in coping with the aftermath.


\subsection{Economics, Earthquakes and Friction} \label{quake}

In early 2009, Alan Greenspan, former Chairman of the Federal reserve, wrote
the following:

``We can model the euphoria and the fear stage of the business cycle. Their
parameters are quite different. ... we have never successfully modelled
the transition from euphoria to fear.'' \\ --- Alan Greenspan,
{\em Financial 
  Times}, March 27th 2009.

The implication is that Central Bank models work well
`most of the time' with suitably calibrated parameters. Occasionally
the parameters suddenly change but once these are measured
the model again works well in the neighbourhood of a new
equilibrium.

The above response to models that suddenly fail is only justified when
the transitions between euphoria and fear and the changes in
parameters are truly exogenously triggered. If they are due to endogenous
causes then the model was never really working before the transition
and it probably won't after the transition either!

A useful analogy here is with earthquakes and seismology. Earthquake
zones {\em appear} to be stable (i.e. in an equilibrium) for very long
periods of time with only very brief, but violent, `transitions'. A
tectonic-plate-denying `equilibrium seismologist'
might argue that the earthquake-free equilibrium model was
essentially correct except for some occasional unpredictable exogenous
events (unobserved meteorite strikes!?) that didn't in any way cast
doubt on the modelling assumptions.

Of course, earthquakes are almost always endogenously generated and
the analogy can be pushed further. An earthquake is  a
very fast shift from one (meta-)stable\footnote{Metastability in
  physics is when a system can stay in a particular state for an
  indefinite amount of time even though it is not the state of lowest
  energy. It occurs when there is some kind of barrier to true
  equilibration.} internal configuration to another and this leads us
consider the concept of `balance-of-forces' in both physics and
economics.

Ever since the time of Walras and Jevons  the idea that there should be a
complete and unique set of equilibrium prices that exactly balances
all of the competing needs and desires of economic agents has offered a
compelling view of a perfectly balanced economy with t\^{a}tonnement
processes somehow achieving this outcome. But this view is based
upon a comparison with physical systems that is misleading. A spring
or piece
of elastic subject to competing forces will achieve a unique
equilibrium but this is because there is no complex internal structure
capable of absorbing any of the stresses without yielding.

A more complicated physical system such as a tectonic fault line
has myriad internal configurations capable of balancing
the forces applied to it --- up to a point. Which
particular configuration exists at any given moment will depend upon the
previous states of the system. And when one small part of the fault
line suddenly shifts this can transfer excess stress to neighbouring
parts resulting in a large cascading failure/earthquake. There is a
balance of forces before the earthquake and after the earthquake but
not {\em during} the earthquake!

A modern economy is arguably the most
complicated man-made construct on the planet with an immensely
intricate internal description which cannot simply be averaged away. The
analogy is also useful in that the fundamental source of earthquakes
is friction. Without it, continental plates would gracefully and
safely glide rather than stick and then briefly grind. Frictions and
stickiness are present in many forms in an economy or financial system
and it should not be a surprise if they cause similar qualitative
effects.

This brings us to the notion of timescales.  In an equilibrium
system there is no notion of {\em any} timescale except for ones imposed
exogenously\footnote{There is no notion of history either. If a system
  is at its unique equilibrium there is no way of telling where it has
  been}.  If one examines an earthquake fault line on a
long-enough timescale, maybe tens of thousands of years, then it
doesn't look like an equilibrium at all. The mere presence of
frictional effects can introduce surprisingly long timescales into a
system via the existence of metastable states. 

If economies feel like they are close to a unique equilibrium maybe
that's because most of the time tomorrow does indeed turn out to
be a lot like yesterday! Over short timescales unique equilibrium
models will frequently appear to be working --- especially when their parameters
are  being updated  to match incoming real-world data!

Finally, it must be pointed out that the analogy between earthquakes
and the models that we will analyze below is not perfect. Fault lines are being
consistently forced in a single direction while the changes experienced
by economies are more random. Also, our main model has a very
small number of variables and only one sticky component and so
`slippage cascades' aren't possible. However even a single sticky
component allows for the existence of an entire interval of equilibria
and complicated transitions between them.

\subsection{Permanence and Path-Dependence}

If the presence of stickiness/frictions in economics does indeed
induce a myriad of coexisting equilibria then phenomena that are not
possible (or require a posteriori model adjustments) in unique
equilibrium models become not just feasible but inevitable. Perhaps
the most obvious of these is {\em permanence}, also known as
remanence, where a system does not revert to its previous state after
an exogenous shock is removed. It is of course a central concern
of macreconomics whether or not economies affected by, say, significant
negative shocks can be expected to have permanently reduced
productivity levels.

For the models studied in this paper, sufficiently small shocks
(whether exogenous or applied by policy makers) will not change the
equilibrium point and a standard linear stability analysis determines
the rate at which the system returns to it.  Larger shocks will move
the equilibrium point along a line of potential equilibria in the
expected direction. But even larger shocks may move the system far
enough away from the equilibrium interval that the return path and
ending point on the interval are very hard to predict. Furthermore, in
neither of the last two cases will the system exhbit a tendency to
return to its pre-shocked state --- the model displays true
permanence. And the model parameters alone cannot determine which
equilibrium a system is currently in without knowing important
information about the prior states of the system --- true path
dependence. This does not, however, prevent the system from being
iterated once the intial conditions are {\em fully} specified.

\subsection{Sticky Models and Indeterminacy}

The most widely-used sticky models are the sticky-prices of Calvo
\cite{calvo1983staggered} and the sticky-information of Mankiw and
Reis \cite{mankiw2002sticky}. These models are conceptually very
similar to each other in that agents do not instantaneously move to
the `correct' price or opinion but rather do so at a fixed rate and
can be represented mathematically by introducing a delay term into the
relevant equations. In the absence of noise the same optimal
equilibrium solution will be reached as if the stickiness were absent.

Continua of possible equilibria can also occur in such models (see for
example \cite{benhabib1999indeterminacy,evans2015observability}) and
is considered an extreme form of {\em indeterminacy}. This is
especially problematic within a Rational Expectations framework since
it makes it (even) harder to justify how the agents' expectations can be
consistent with the model.

Our hypothetical agents are less rational than those above. They are
truly stuck (not just delayed) until forced to adjust by the
discrepancy with the actual inflation rate.\footnote{This is now very
  close to the situation where a frictional force has to be overcome
  before an object will move.} 
If an equilibrium is reached it is chosen by the prior states of the
system and not by modeling assumptions about the future and, as we shall see, a continuum of equilibria
is an intrinsic feature and not an inconvenience that
occurs only in certain special cases (such as a passive interest-rate
policy \cite{calvo1983staggered,antinolfi2007monetary}).

The research into how expectations are formed is extensive but far
from conclusive, see for example
\cite{curtin2010inflation,rudd2006can,branch2007sticky,carroll2003macroeconomic,mankiw2003disagreement}. However
the idea of threshold effects and a `harmless interval' of inflation
is not new in economics
\cite{Threshold2010,kremer2013inflation,2013inflation,F2010inflation,2001threshold}.
In the absence of any exogenous forcing it would be very easy to
distinguish between Calvo-type stickiness and the stuck-then-dragged
behavior we investigate here --- indeed Calvo stickiness would most
likely be observed since agents could tell far more easily over time
that, for example, their wage demands were too low and they were
losing purchasing power. However, given the uncertainty of reality and
the very limited cognitive skills or interest in forecasting
of most economic agents, that may no longer hold.

Our model of expectation formation is thus both mathematically
tractable and has some basis in both observed data (see also \cite{g, laura})
and models of bounded
rationality. As such it provides a potentially useful,
analytically tractable, alternative
to staggered/delayed models --- and one with additional
complexity and explanatory power.

\subsection{Bounded Rationality and Aggregation}

As mentioned above, the standard approach to the problem of
aggregating expectations is to introduce a `Representative Agent'
whose expectations are fully-informed and rational and consistent with
the model itself. Here, an aggregation of {\em boundedly}
rational agents into a similar Representative  is required.

Our approach is similar in spirit to that of De Grauwe
\cite{DeGrauwe2012}.
 In \cite{DeGrauwe2012} both
the expectations terms in inflation and output gap are linear
combinations of the expectations of two kinds of agent --- rational
`fundamentalists' and boundedly rational `extrapolators' --- with the
probability of an agent using each being dictated by discrete choice
theory \cite{anderson1993discrete,brock1997rational}. He then showed numerically that cycles of booms-and-busts
occurred with changes in the `animal spirits' and corresponding
non-Gaussian `fat-tailed' disributions for the model variables.
Discrete choice theory is the aggregating mechanism that De Grauwe
uses to avoid ending up with an agent-based model where each agent has
to be individually simulated.

We start from the empirical evidence cited above that individual agents'
expectations are often sticky and may lag behind the currently
observable values before they start to move. We also posit that
this gap between future expectations and current reality cannot grow
too large. We then imbue our now boundedly rational Representative
Agent with these same properties.
This leads us in a very natural way to the play operator that is described
fully in Section \ref{sec_play}. And while it is certainly not a fully-justified
aggregation procedure neither are the others mentioned above!


\bigskip
\subsection{Outline of the paper}

We start from a dynamic stochastic general equilibrium (DSGE)
macroeconomics model, which includes aggregate demand and aggregate
supply equations
\begin{equation}\label{eqn:M1}
\begin{array}{c}
\hspace{-16mm} y_t=y_{t-1}-a(r_t-p_t)+\epsilon_t,\\
x_t=b_1p_t+(1-b_1)x_{t-1}+b_2y_t+\eta_t\\
\end{array}
\end{equation}
augmented with the rate-setting rule 
\begin{equation}\label{eqn:M1'}
r_t=c_1x_t+c_2y_t,
\end{equation}
where $y_t$ is output gap (or unemployment rate, or another measure of
economic activity such as gross domestic product), $x_t$ is inflation
rate, $r_t$ is interest rate, $p_t$ is the economic agents' aggregate
expectation of future inflation rate and $\epsilon_t$, $\eta_t$ are
exogenous noise terms. All the parameters are non-negative and in
addition, $b_1<1$.  This model is close to the starting
model used in \cite{DeGrauwe2012} but simpler in that we do not
include the aggregate expectation of the output gap and the
correlation between the subsequent values of the interest rate. We
also choose to remove the noise term from the interest rate update
rule.  The inclusion of such factors does not affect our most
significant qualitative observations, but would complicate some
aspects of the rigorous analysis that we present.

The novelty of our modeling strategy is in how we define the
relationship between the aggregate expectation of inflation $p_t$ and
the inflation rate $x_t$.  This relationship is defined precisely in the next
section where we introduce the play operator to model the economic
agents' aggregate expectation of future inflation.

In Sections \ref{staloc}-\ref{staglo} we present the main stability analysis for various
parameter regimes, with some
details relegated to Appendices. The stability properties of the
system are not as clear cut as in a truly linear system. 
In fact, our equations define a {\em piecewise linear} (PWL) system, and
certain nonlinear effects come into play. 
In particular, in nonlinear
systems an
equilibrium may only be {\em locally} stable. This means that the
equilibrium is only stable to perturbations of a certain size --- ones
that don't move  the system outside of a `basin of attraction' --- and
this phenomenon is responsible for much of the interesting dynamics
in the presence of shocks.

In Sections \ref{numerics1}-\ref{numerics6} we present various numerical simulations. We are
particularly interested in the transitions between equilibrium states
caused by exogenous shocks, and the effects of increasing or
decreasing stickiness. Where possible we compare results against the
non-sticky model. Permanence is the rule not the exception and there
are even parameter regimes where a large enough shock will completely
destabilise an apparently stable system via a runaway inflation
mechanism.

We also compare the statistical output of the model against that of De
Grauwe \cite{DeGrauwe2012} at similar parameters and see the same boom-and-bust
cyclicality and heavy-tailed distributions.

Then, in Section \ref{multi} we briefly consider a more complicated
version of the model with three representative agents all with
different levels of stickiness. This is primarily to demonstrate that
multiple play operators can indeed be used together to simulate
different representative agents within a model and that the most
important qualitative features are unchanged.

Finally, in Section \ref{interest} we emphasize that play
operators are not just a potential tool for modeling expectations --- we
remove the stickiness from the inflation expectations and add
it into the response of the Central Bank instead. We perform a second
stability analysis and obtain some interesting new effects ---  there
is the possibility of (quasi)-periodic behavior in the absence of noise and the
stickiness does appear to destabilize  equilibria. We
conclude with a summary of the main results and some suggestions for
future work.

\section{The model} \label{sec:model}



\subsection{Play and Stop Operators}\label{sec_play}

We assume the following rules that define the variations of the
expectation of future inflation rate $p_t$ with the actual inflation
rate $x_t$ at integer times $t$:
\begin{itemize}
\item[(i)] The value of the difference $|p_t-x_t|$ never exceeds a certain bound $\rho$;

\item[(ii)] As long as the above restriction is satisfied, the expectation does not change, i.e. $|x_t-p_{t-1}|\le \rho$ implies $p_t=p_{t-1}$;

\item[(iii)] If the expectation has to change, it makes the minimal
  increment consistent with constraint (i).
\end{itemize}

Rule (ii) introduces stickiness in the dependence of $p_t$ on $x_t$,
while (i) states that the expected inflation rate cannot
deviate from the actual rate more than prescribed by a threshold value
$\rho$.  Hence $p_t$ follows $x_t$ reasonably closely but on the other
hand is conservative because it remains indifferent to variations of
$x_t$ limited to a (moving) window $p-\rho\le x\le p+\rho$.  The last
rule (iii) enforces continuity of the relationship between $p_t$ and
$x_t$ and, in this sense, can be considered as a technical modeling
assumption that is mathematically convenient.

\begin{figure}[h!]
\centering
\begin{subfigure}{0.45\textwidth}
\includegraphics[width=\textwidth]{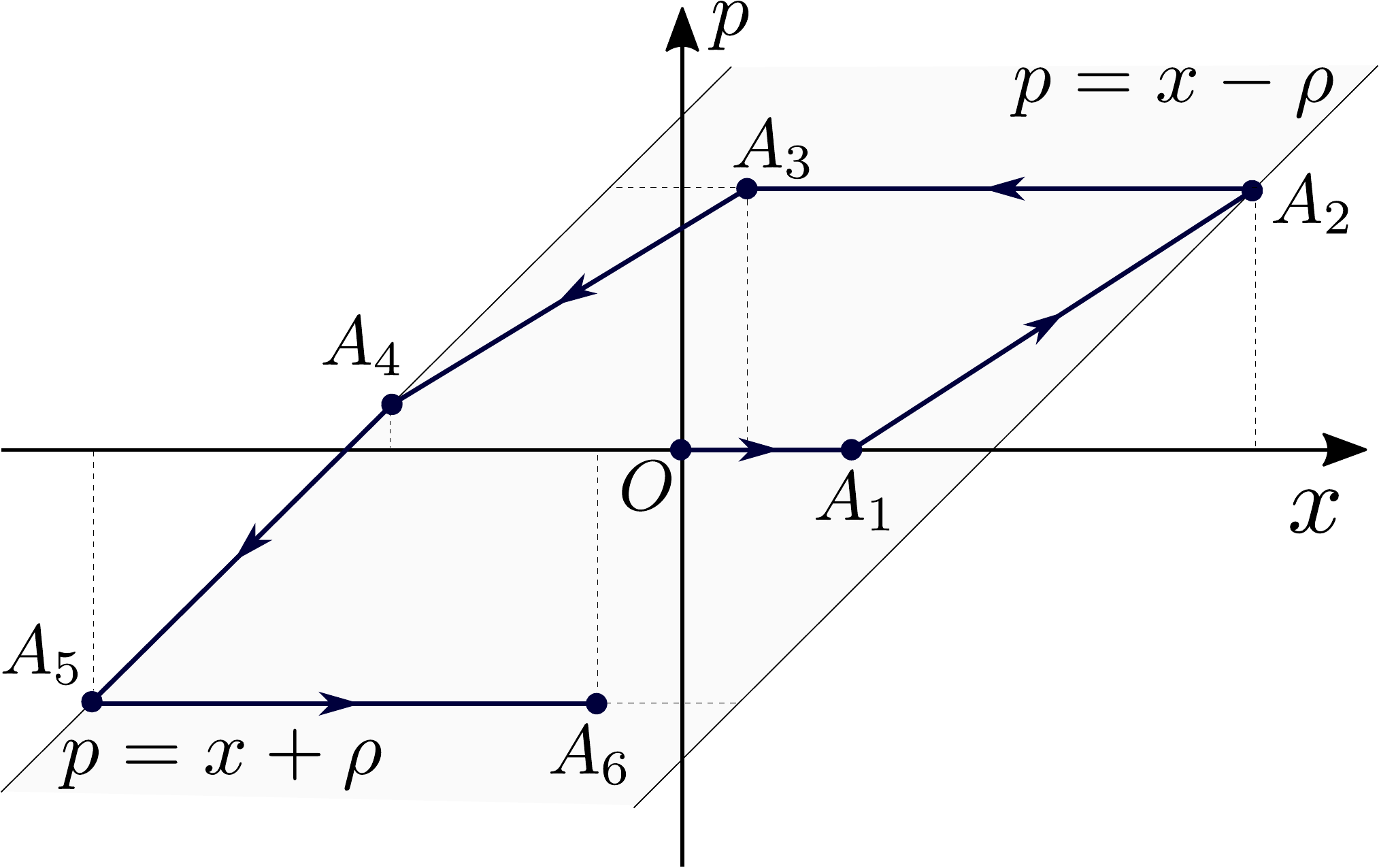}
\caption{} 
\label{fig:PlayO}
\end{subfigure} 
\quad
\begin{subfigure}{0.45\textwidth}
\includegraphics[width=\textwidth]{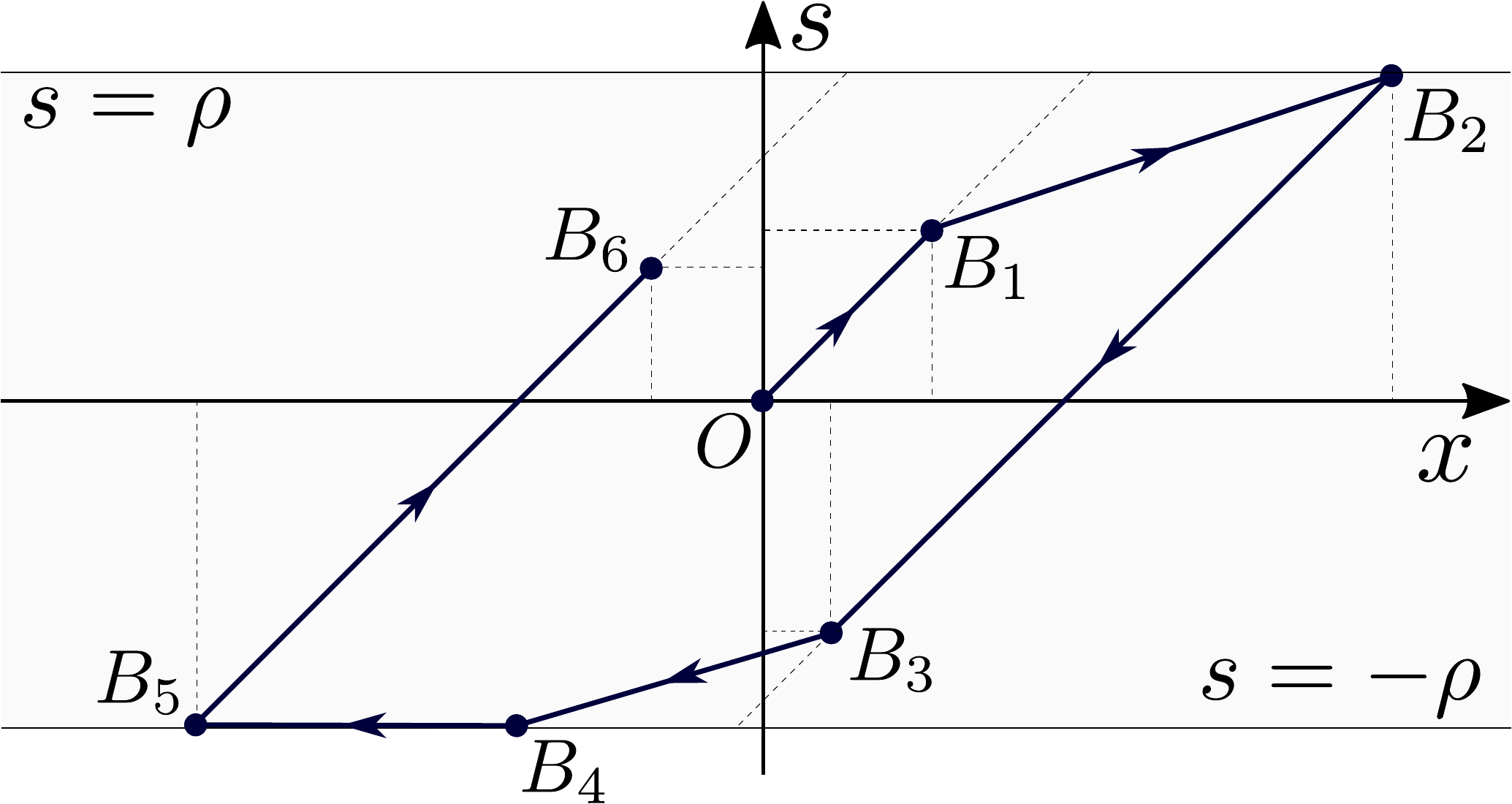}
\caption{} 
\label{fig:StopO}
\end{subfigure}
\caption{(a) An illustration of the input-output sequence of the (a) play operator and (b) stop operator. 
(a)	The polyline $OA_1A_2A_3A_4A_5A_6$ represents a sample input-output trajectory for the play operator.
	 	The input-output pair $(x,p)$ is bounded to the gray strip between the two parallel lines $p=x\pm\rho$. 
		In \cite{laura}, this strip is called {\em band of inactivity},  the line $x=x-\rho$ is called {\em upward spurt} line  while the line $p=x+\rho$ is called {\em downward spurt} line. 
		The ouptut $p$ remains unchanged for a transition from $(x_{t-1},p_{t-1})$ to the next point $(x_t,p_t)$ as long as the pair $(x_{t},p_{t-1})$ fits to the band of inactivity (for example,		the transitions from $A_2=(x_2,p_2)$ to $A_3=(x_3,p_3)$ with $p_2=p_3$ or from $A_5=(x_5,p_5)$ to $A_6=(x_6,p_6)$ with $p_5=p_6$). 
		If $x_t>x_{t-1}$ and the point $(x_{t},p_{t-1})$ lies to the right of the inactivity band, then the output increases resulting in 		the point $(x_t,p_t)$ to lie on the upward spurt curve (for example, the transition from $A_1=(x_1,p_1)$ to $A_2=(x_2,p_2)$). 		
		Similarly, if $x_t<x_{t-1}$ and the point $(x_{t},p_{t-1})$ lies to the left of the inactivity band, then 		the output decreases and the point $(x_t,p_t)$ lies on the downward spurt line (for example,  the transition from $A_3=(x_3,p_3)$ to $A_4=(x_4,p_4)$).
		(b) The input-output trajectory of the dual stop operator corresponding to the trajectory of the play operator shown in panel (a). Here $s_t=x_t-p_t$; the trajectory is limited to the horizontal strip $-\rho\le s\le \rho$ at all times.
		}
		\label{fig:HOperators}
\end{figure}

Rules (i)--(iii) are expressed by the formula
\begin{equation}\label{formula}
p_{t}=x_t+\Phi_\rho(p_{t-1}-x_t)
\end{equation}
with the piecewise linear saturation function 
\begin{equation}\label{formula'}
\Phi_\rho(x)=\left\{\begin{array}{rlc} 
\rho & {\rm if} & x\ge \rho,\\
x & {\rm if} & -\rho<x<\rho,\\
-\rho & {\rm if} & x\le -\rho.
\end{array}
\right.
\end{equation}
Relationship \eqref{formula} is known as the {\em play} operator with {\em threshold} $\rho$, see Fig.~\ref{fig:HOperators}(a).
A dual relationship 
\begin{equation}\label{formula''}
s_{t}=\Phi_\rho(x_t-x_{t-1}+s_{t-1})
\end{equation}
between $x_t$ and the variable 
\[
s_t=x_t-p_t
\]
is referred to as the {\em stop} operator, see
Fig.~\ref{fig:HOperators}(b). In the context of our model, $s_t$
measures the difference between the inflation rate and the expectation
of the future inflation rate, hence $s_t$ remains within the bound
$|s_t|\le \rho$ at all times.  Interestingly the explicit relationship
\eqref{formula} has been observed in actual economic data \cite{g, laura}.

One can think of the play operator as having two modes. A `stuck mode'
where it will not respond to small changes in the input and a
`dragged mode' where the absolute difference between the input and output are
at the maximum allowable and changes to the input, in the correct
direction, will drag the output along with it.


Equations \eqref{formula} and \eqref{formula''} will now be denoted by 
\begin{equation}\label{playstop}
p_t=\mathscr{P}_\rho[x_t], \qquad s_t=x_t-p_t=\mathscr{S}_\rho[x_t],
\end{equation}
where $\mathscr{P}_\rho$ and $\mathscr{S}_\rho$ are  the \textit{play} and \textit{stop} operators with threshold $\rho$, respectively. 


\subsection{A model with sticky inflation expectations}\label{inve}
Equations \eqref{eqn:M1} and \eqref{eqn:M1'}, completed with formulas
\eqref{formula} and \eqref{formula'}, form a closed model for the
evolution of the aggregated variables $x_t, y_t, r_t, p_t$.  However,
the dependence of these quantities at time $t$ upon their values at
 time $t-1$ is implicit. In order to implement the model, we
proceed by solving equations \eqref{eqn:M1}--\eqref{formula'} with
respect to the variables $x_t, y_t$.  As shown in Appendix A, the
model can be written in the following equivalent form:
\begin{equation}\label{eqn:SA1}
z_t={A}z_{t-1}+s_t{d}+N\xi_t
\end{equation}
where $z_t=(y_t,x_t)^\top$, $\xi_t=(\epsilon_t,\eta_t)^\top$, the superscript ${}^\top$ denotes transposition, the matrices $A,N$ and the column vector $d$ are defined by
\begin{equation}\label{AA}
{A}=\frac1\Delta\left(\begin{array}{cc}
{1-b_1}&{a(1-b_1)(1-c_1)}\\
{b_2}&{(1-b_1)(1+ac_2)}
\end{array}\right), \ \ \ 
N=\frac1\Delta\left(
\begin{array}{cc}
{1-b_1}&{a(1-c_1)}\\
{ b_2}&{ 1+ac_2}
\end{array}\right),
\end{equation}
\[
{d}=\frac1\Delta\left(\begin{array}{c}
{a(b_1c_1-1)}\\
-({ab_2+b_1(1+ac_2)})
\end{array}
\right)
\]
with
\begin{equation}\label{Delta}
\Delta={(1-b_1)(1+ac_2)+ab_2(c_1-1)}
\end{equation}
and 
$s_t=x_t-p_t$ is defined by the equation
\begin{equation}\label{eqn:SA1'}
s_t=\frac1{1+\alpha}\,\Phi_{(1+\alpha)\rho}(f_t-f_{t-1}+s_{t-1})
\end{equation}
with
\begin{equation}\label{alpha}
\alpha=\frac{\Delta}{b_1(1+ac_2)+ab_2},
\end{equation}
\begin{equation}\label{f}
f_t=\frac{\alpha}{\Delta}\bigl({b_2}y_{t-1} +{(1-b_1)(1+ac_2)}x_{t-1}+{b_2\epsilon_t+(1+ac_2)\eta_t}\bigr).
\end{equation}
Equations \eqref{eqn:SA1}, \eqref{eqn:SA1'} express $y_t$, $x_t$ and $s_t=x_t-p_t$
explicitly in terms of the previous values of the same variables 
and the exogenous noise  $\epsilon_t$, $\eta_t$. We use these equations in
all the simulations that follow.

We shall refer to the variable $s_t = x_t-p_t$ as the {\em perception gap}.
Note that \eqref{eqn:SA1'} defines a stop operator with input $f_t$
and threshold $(1+\alpha)\rho$, which is different from $\rho$
(cf.~\eqref{formula'}) and so \eqref{eqn:SA1'} can be written as
\[
s_t=\frac{1}{1+\alpha}\mathscr{S}_{(1+\alpha)\rho}[f_t]
\]
using the notation \eqref{playstop}. It is important to note that 
the transition to equations  \eqref{eqn:SA1}, \eqref{eqn:SA1'} is justified under the condition that $\alpha$ is positive,
and we assume this constraint to hold in the rest of the paper. In
particular, $\alpha>0$ whenever $c_1>1$ (see Section \ref{staglo}).

\subsection{An entire line segment of equilibrium points}

\begin{figure}[ht]
\begin{center}
	\begin{subfigure}{0.45\textwidth}
		\includegraphics[scale=0.3]{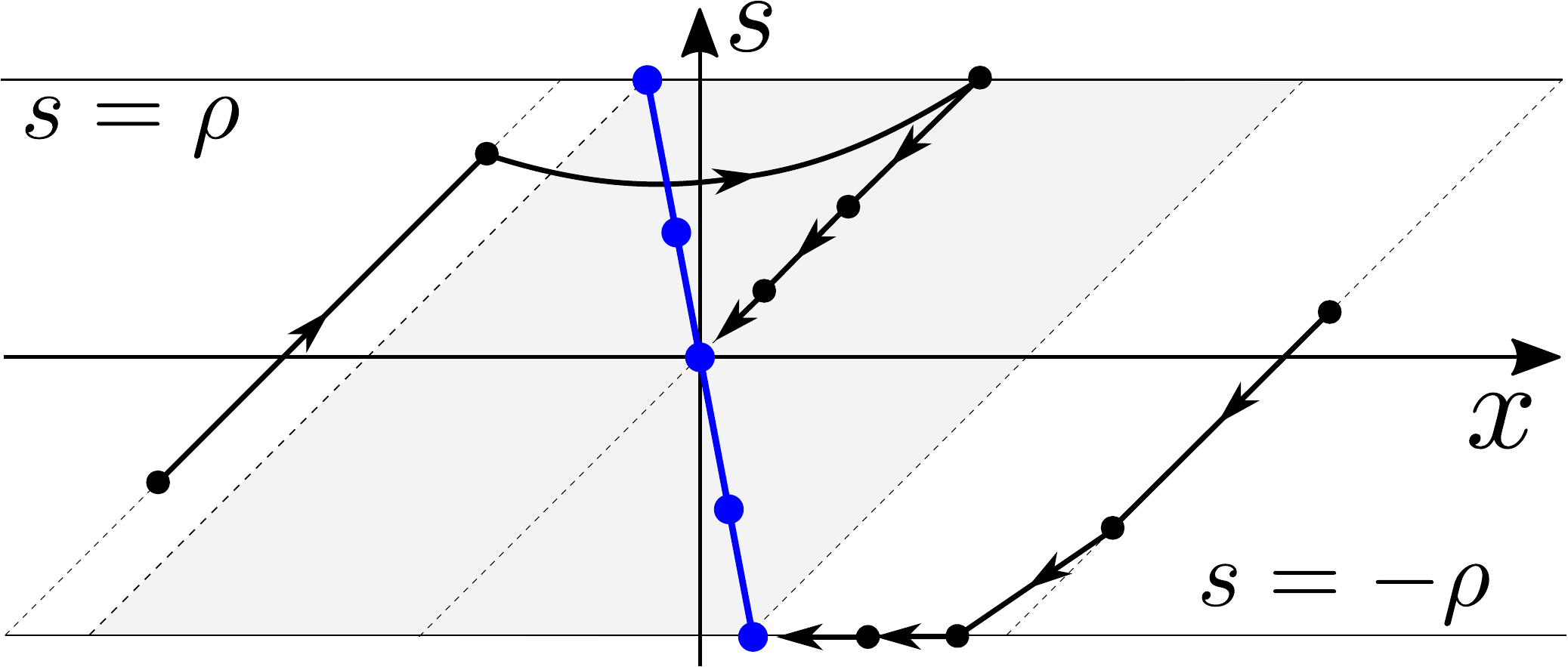}
		\caption{}
	\end{subfigure} 
 \quad
\begin{subfigure}{0.45\textwidth}
		\includegraphics[scale=0.3]{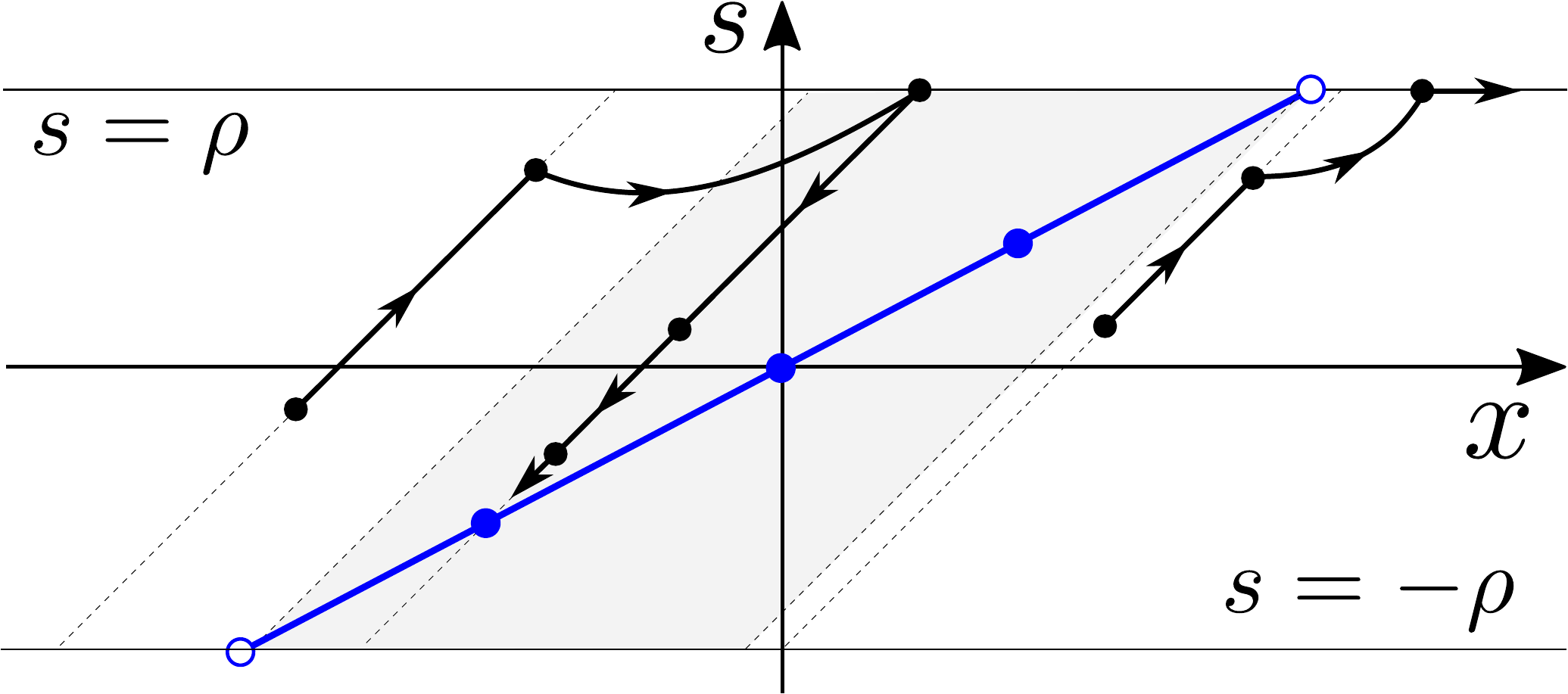}
		\caption{}
	\end{subfigure} 
\vspace{-1mm}
\caption{The projection of the line segment of equilibrium points (blue) onto the $(x,s)$ plane for (a) $c_1>1$ and (b) $c_1<1$. 
The segment has a negative slope in (a) and a positive slope in (b).
Sample trajectories of system \eqref{eqn:SA1''} are shown in black. 
}
\label{fig:C1Less1}
\end{center}
\end{figure}

We begin the analysis of the model \eqref{eqn:SA1}, \eqref{eqn:SA1'}
by looking at the case of no exogenous noise, i.e. we set $\xi_t=0$
and consider the equation
\begin{equation}\label{eqn:SA1''}
z_t={A}z_{t-1}+s_t{d},\qquad z_t=(y_t,x_t)^\top
\end{equation} 
instead of \eqref{eqn:SA1} with $s_t$ defined by 
\eqref{eqn:SA1'}, \eqref{alpha} and 
\begin{equation}\label{f}
f_t=\frac{\alpha}\Delta\bigl({b_2} y_{t-1} +{(1-b_1)(1+ac_2)}x_{t-1}\bigr).
\end{equation}
This model has an entire line segment of equilibrium points which
corresponds to a continuum of feasible equilibrium states of the
economy as a function of the inflation expectations of economic
agents.  Indeed, equation \eqref{eqn:SA1''} implies
\begin{equation}\label{eqn:SA1E}
z_*=s_*(\mathbb{I}-{A})^{-1}{d}=s_*\begin{pmatrix}
\frac{b_1}{b_2},&\frac{b_2+b_1c_2}{b_2(1-c_1)}\end{pmatrix}^\top 
\end{equation}
for an equilibrium point $z_*=(x_*,y_*)^\top$,
where $\mathbb{I}$ is the $2\times2$ identity matrix. 
Hence one obtains a different equilibrium for each admissible value 
of the perception gap variable $s_*$, i.e.
$-\rho\le s_*\le \rho$. Thus, the set of all equilibrium points, which can be denoted as $z_*(s_*)$
for different $s_*$,
can be naturally thought of as a line segment in the phase space of the system, see Fig.~\ref{fig:C1Less1}.
In particular, the value of the output gap at an equilibrium, $y_*(s_*)$
ranges over the interval $[-\rho b_1/b_2,\rho b_1/b_2]$ 
and the equilibrium value of the actual inflation  
belongs to the range 
\[
x_*(s_*)=s_*\,\frac{b_2+b_1c_2}{b_2(1-c_1)} \quad {\rm with} \quad -\rho\le s_*\le \rho.
\]
Interestingly, at least in this simple model, the range of equilibrium
values of the output gap is unaffected by the controls $c_1$, $c_2$
applied by the regulator through Taylor's rule \eqref{eqn:M1'}.
However, these controls do affect the range of possible values of the
equilibrium inflation rate.

Equation \eqref{eqn:SA1E} indicates the difference between the cases  $c_1>1$ and $c_1<1$.
	When $c_1>1$, the equilibrium $z_*(\rho)$ corresponding to the lowest expectation of inflation
	has the highest value of the output gap and the lowest inflation of all the equilibrium points.
	Similarly, 
	the equilibrium $z_*(-\rho)$ with the highest expectation of inflation
	has the lowest value of the output gap and the highest inflation.
	On the other hand, in case $c_1<1$, the equilibrium $z_*(\rho)$ with the highest output gap value 
	has simultaneously the highest inflation rate.
	
	The difference between the cases $c_1>1$ and $c_1<1$ will be
        further highlighted in Section \ref{staglo}.

\subsection{Local stability analysis}\label{staloc}
System \eqref{eqn:SA1}, \eqref{eqn:SA1'} is locally linear in some neighborhood of any equilibrium
point from the linear segment \eqref{eqn:SA1E} with the exception of the two end points 
$z_*(\pm \rho)$ corresponding to equilibria where the play is right
at one end of its inactive band. In other words, for sufficiently small deviations of the vector
$z_t=(y_t,x_t)^\top$ from an interior equilibrium $z_*(s_*)$, system \eqref{eqn:SA1''} is equivalent to
\begin{equation}\label{eqn:SA2}
z_t-z_*(s_*)=B(z_t-z_*(s_*))
\end{equation}
where 
\[{B}=\begin{pmatrix}
\frac{1}{1+a(b_2c_1+c_2)}&\frac{a(b_1-1)c_1}{1+a(b_2c_1+c_2)}\\
\frac{b_2}{1+a(b_2c_1+c_2)}&\frac{(1-b_1)(1+ac_2)}{1+a(b_2c_1+c_2)}
\end{pmatrix}
\]
As shown in Appendix B, the matrix $B$ is stable for any admissible set of
parameter values, hence every equilibrium with $|s_*|<\rho$ is locally
stable. This local stability ensures that if a {\em sufficiently small}
perturbation is applied to the system residing at an equilibrium
$z_*(s_*)$, removing the perturbation returns the system to
the same equilibrium. Further, the eigenvalues of the matrix $B$ determine
how quickly (or slowly) the system returns to the equilibrium
state. This situation is of course very similar to the expected response in a
fully linear equilibrium model. The dependence of the eigenvalues of the parameters of the
system is discussed in Appendix C.

However, the situation for these interior equilibria changes markedly for larger perturbations. This
is related to the stability properties of the two extreme equilibria
$z_*(\pm \rho)$ and  is far more subtle as discussed in the next section.
In particular, the basin of attraction of the equilibrium 
decreases and finally vanishes as one approaches either of the extreme
equilibrium points along the line segment \eqref{eqn:SA1E} (the
extreme equilibria themselves are stable but not asymptotically stable).

\subsection{Global stability analysis}\label{staglo}
System \eqref{eqn:SA1''} without stickiness ($\rho=0$) simply has the form
\begin{equation}\label{eqn:SA1'''}
z_t={A}z_{t-1}.
\end{equation} 
As shown in Appendix B, its unique zero equilibrium is globally stable if $c_1>1$ and is unstable if $c_1<1$.

For system \eqref{eqn:SA1''} with stickiness ($\rho>0$), equation \eqref{eqn:SA1'''} approximates the
dynamics far from equilibrium points because the term $s_t$ in
\eqref{eqn:SA1''} is bounded in absolute value by $\rho$. In
particular, since \eqref{eqn:SA1'''} is unstable for $c_1<1$, so is
system \eqref{eqn:SA1''}. This creates the possibility of run-away
inflation at these  values of $c_1$ (see Section \ref{infln}).

Interestingly, the same condition $c_1>1$ that ensures the global stability of system \eqref{eqn:SA1'''}, also guarantees the global stability of the set of equilibrium states for the sticky nonlinear system 
\eqref{eqn:SA1''}. 
In order to show this, one can use a family of {\em Lyapunov functions}
\[
\begin{array}{rcl}
V(x_t,s_t,\nabla_t x, \nabla_t s)&=&\frac12\bigl(C (\nabla_t x)^2+ G (\nabla_t s)^2 + (C x_t+Gs_t)^2\bigr)\\~\\
&+& \gamma \bigl( (C x_t+Gs_t)\nabla_t x  +\frac{H}{2C} (C x_t+Gs_t)^2\bigr),
\end{array}
\]
where $\nabla_t u=u_{t}-u_{t-1}$, $u=x,s$. 
A proper choice of the parameters $C, G,H,\gamma$ ensures that such a function is non-negative, achieves its minimum zero value on the linear interval of equilibrium states, and decreases to zero along every trajectory of system \eqref{eqn:SA1''}. This allows us to prove that 
 every trajectory of system \eqref{eqn:SA1''} converges to one of the equilibrium states \eqref{eqn:SA1E}. 
 In the interest of space, details of the proof are omitted here and will be presented elsewhere.
 
 For system \eqref{eqn:SA1} with noise, this global stability result implies that trajectories tend to return towards the segment of equilibrium points after large fluctuations and hover in a vicinity of equilibrium states for extended periods of time. The rate with which the system returns towards the line segment of equilibrium states after a large perturbation is removed is determined by the eigenvalues of the matrix $A$, see Appendix C.

 \section{Numerical results}\label{numerics}
 
 \subsection{Parameter values}\label{numerics1}
 
 The default parameter set that we use for numerical simulation is the
 same as in \cite{DeGrauwe2012}, see Table \ref{tab:DS1}, and we shall
 explore in detail the surrounding parameter space.
 \begin{table}[h]
 	\begin{center}
 		\begin{tabular}{|c||c|c|c|c|c|}
 			\hline
 			Parameters &  $a$&$b_1$ & $b_2$ & $c_1 $& $c_2$ \\
 			\hline
 			Values & $0.2$ & $0.5$ & $0.05$ & $1.5$ & $0.5$  \\
 			\hline
 		\end{tabular}
 		\caption{The set of parameter values. }\label{tab:DS1}
 	\end{center} 
 \end{table}
  Note that, as an example, if with the above parameters we choose $\rho=\frac12$ then the components of the  equilibrium points  $z_*(s_*)=(y_*(s_*),x_*(s_*)^{\top}$ range over the intervals 
 \[
y_*(s_*)\in[-5,5],\qquad x_*(s_*)\in[-6,6].
 \]

The choice of $\rho$ is somewhat arbitrary as there is of course no
corresponding reference parameter in \cite{DeGrauwe2012} and so in many
of the simulations it will be varied. Also it should be emphasized
that these reference parameters are motivated by \cite{DeGrauwe2012} but
very similar numerical results were obtained for other choices. 

 \subsection{Lower inflation volatility due to stickiness}
 The range of the equilibrium points of the system is directly
 proportional to the threshold value $\rho$ of the play operator
 because the perception gap $s_*$ in \eqref{eqn:SA1E} can take
 any value in the interval $-\rho\le s_*\le \rho$. In particular,
 $\rho=0$ corresponds to the system without stickiness in which the
 expectation of inflation coincides with the current inflation rate,
 $p=x$.  This system is simply described by the equation
  \begin{equation}\label{eq0}
z_t=Az_{t-1}+N\xi_t
  \end{equation}
    (cf.~\eqref{eqn:SA1}). In the absence of noise, it
  has a unique equlibrium at $x=y=0$.
 
\begin{figure}[h!]
 		\centering
 		\begin{subfigure}{0.45\textwidth}
 			\centering
 			\includegraphics[width=\textwidth]{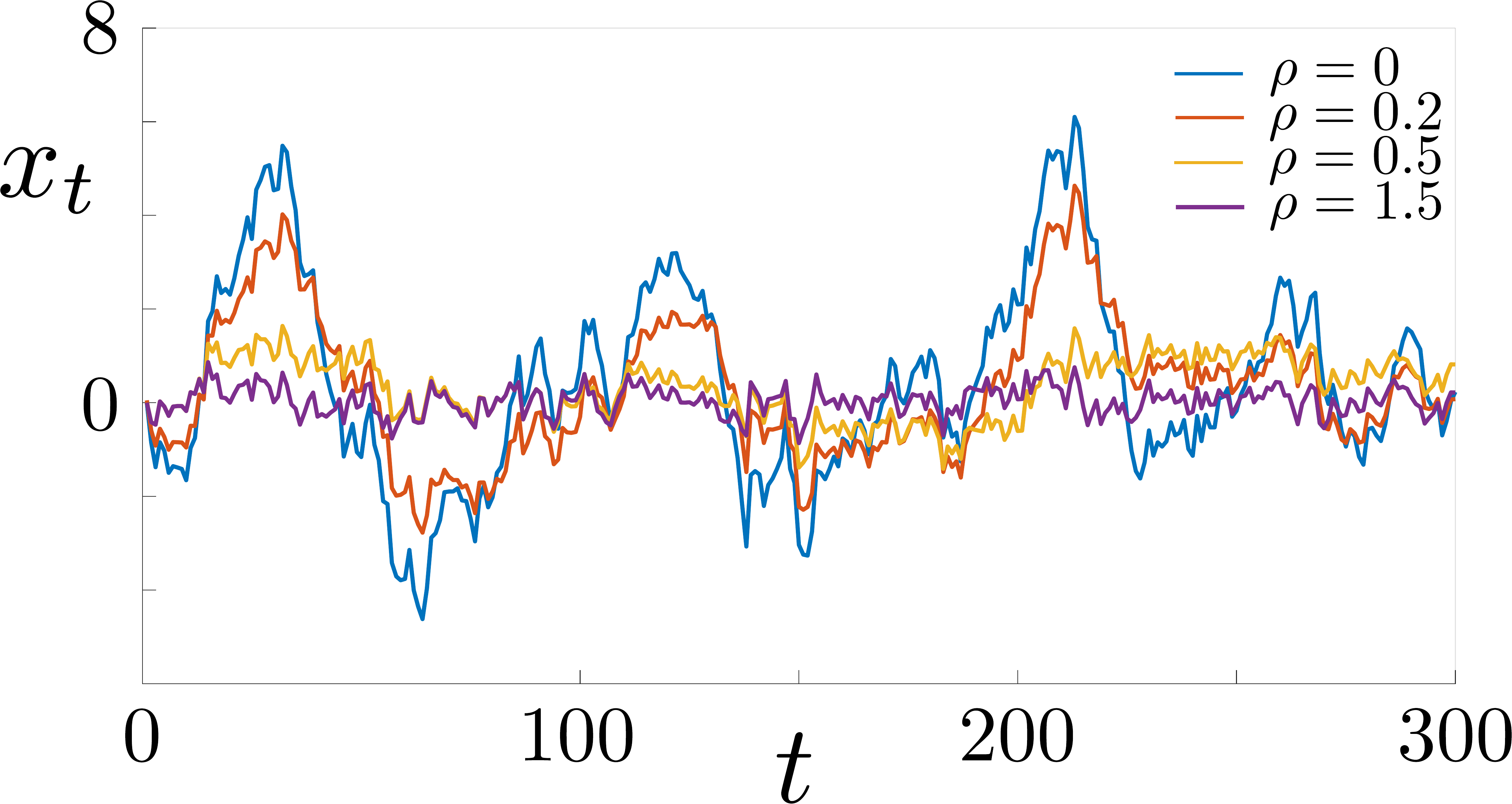}
 			\caption{}
			\label{fig:2DWSOx}
 		\end{subfigure}
 		\qquad	
 		\begin{subfigure}{0.45\textwidth}
 			\centering
 			\includegraphics[width=\textwidth]{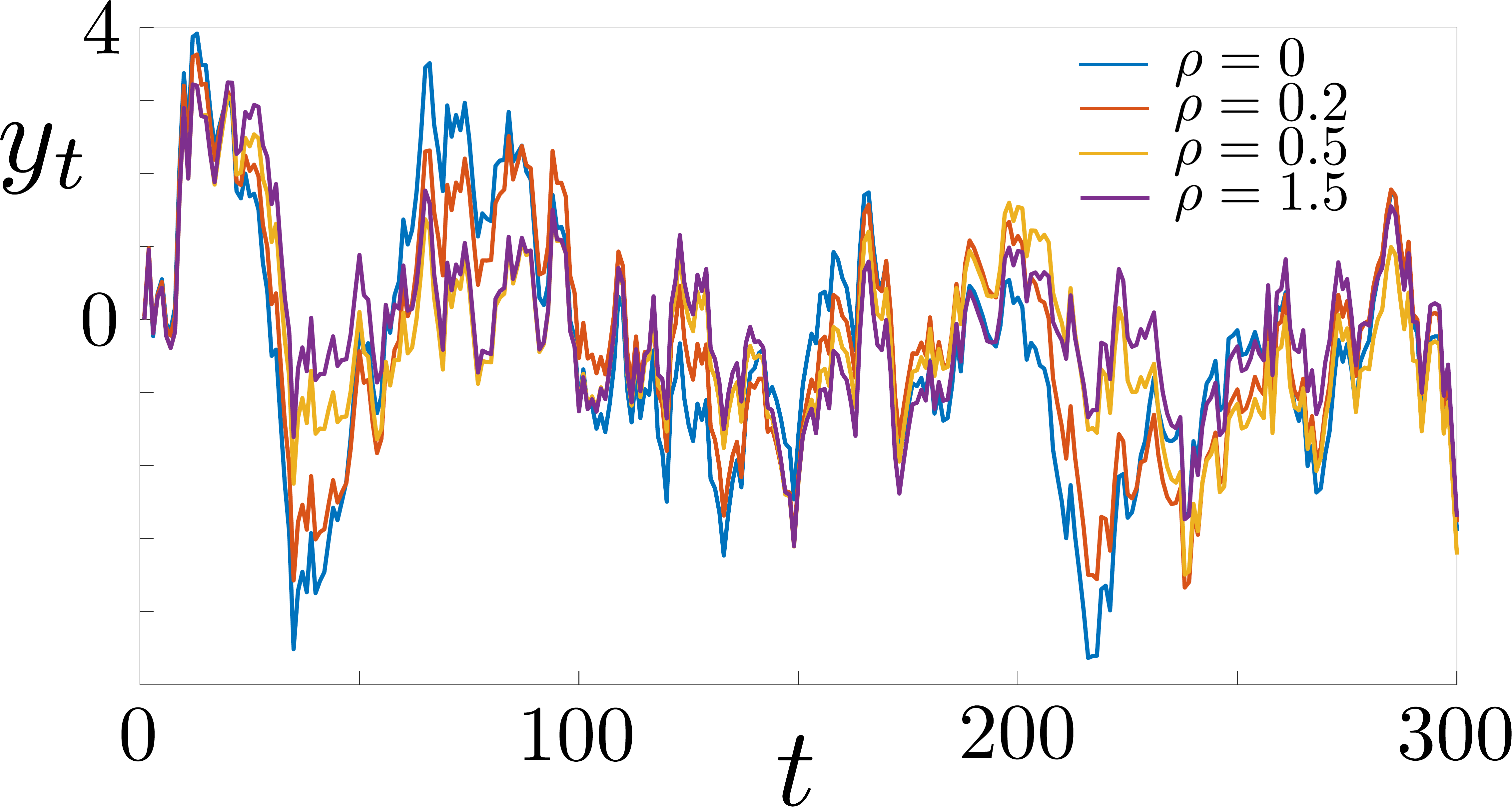}
 			\caption{} 
			\label{fig:2DWSOy}
 		\end{subfigure}
		
		\begin{subfigure}{0.45\textwidth}
 			\centering
 			\includegraphics[width=\textwidth]{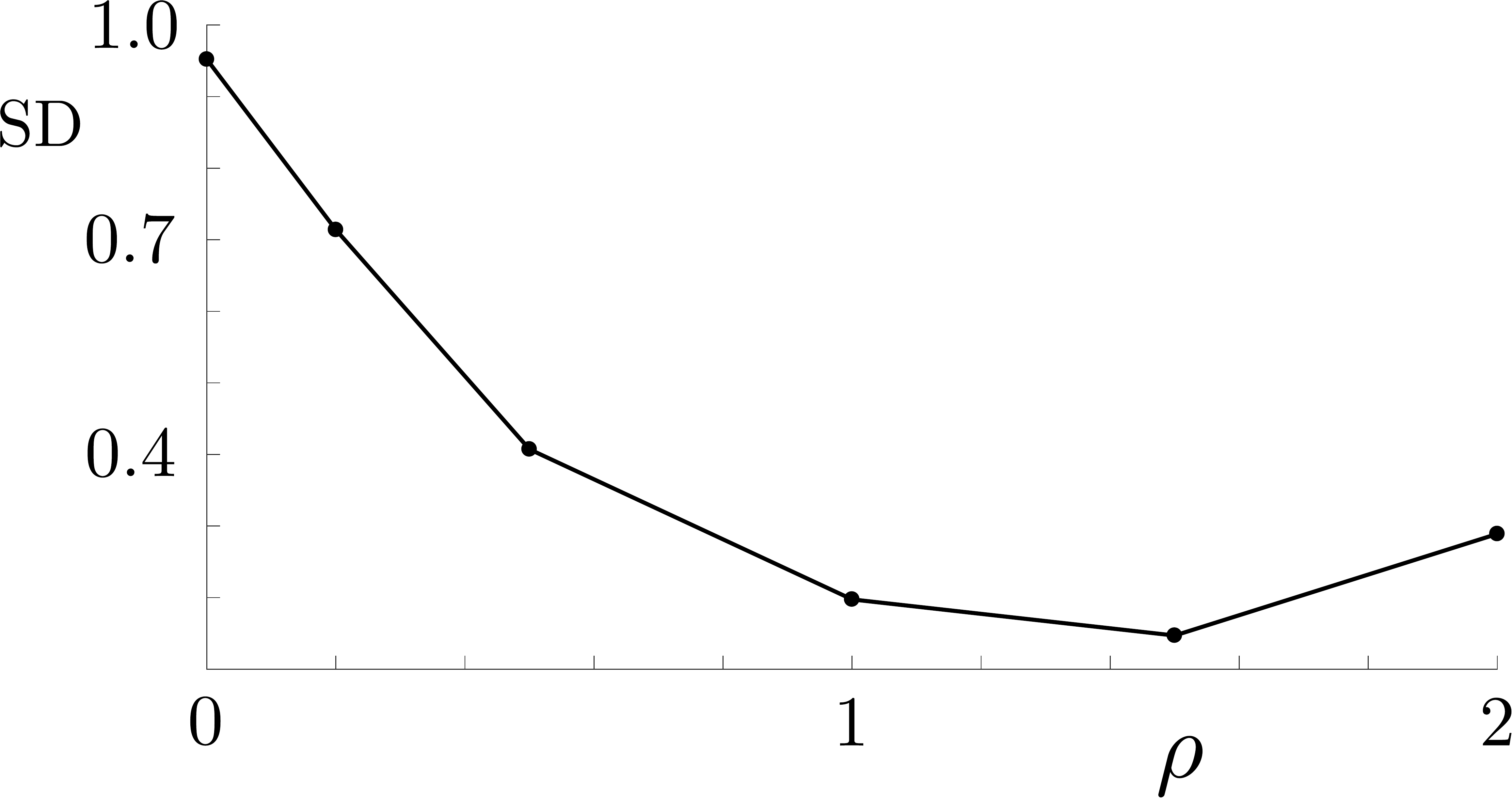}
 			\caption{}
			\label{fig:2DWSOx1}
 		\end{subfigure}	
 		\qquad			
 		\begin{subfigure}{0.45\textwidth}
 			\includegraphics[width=\textwidth]{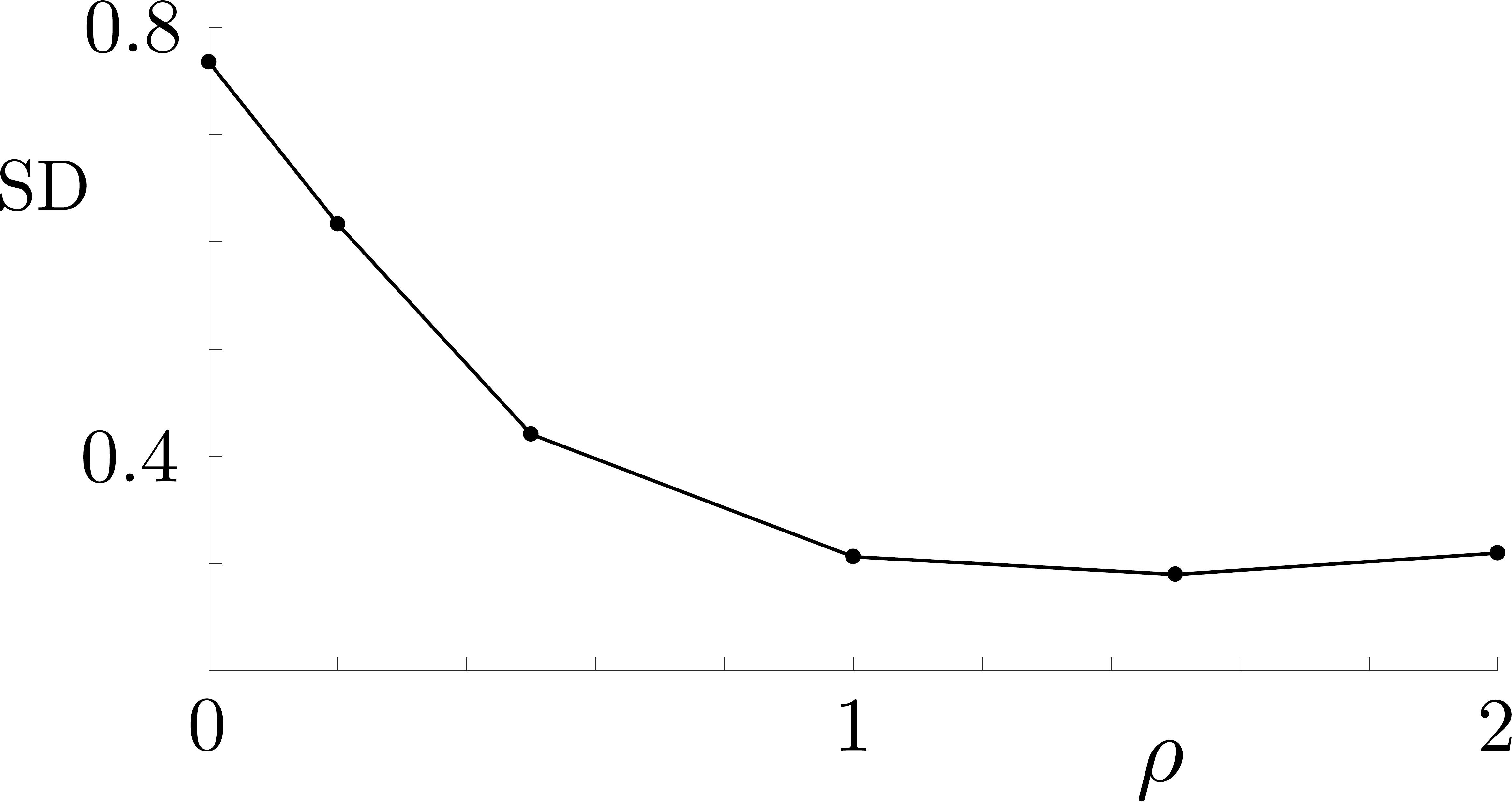}
 			\caption{} 
			\label{fig:2DWSOy1}
 		\end{subfigure}
 		\caption{ Trajectories of (a)  inflation rate $x_t$ and (b) output gap $y_t$. Measure of 
		volatility of (c) inflation rate  and (d) output gap  for different values of $\rho$ with standard deviation (SD).}
	\label{fig:2DWSO}
 	\end{figure}

  The sticky system exhibits lower volatility in the  inflation rate than the
  system without stickiness, see Fig.~\ref{fig:2DWSO}. This can be
  explained by the stability properties of matrices $A$ and $B$ where
  $B$ is the linearization matrix of \eqref{eqn:SA2} for the sticky
  system at an equilibrium. For the parameter values of Table
  \ref{tab:DS1}, the spectral radius of the matrix $B$ is smaller than
  the spectral radius of $A$ (see Appendix C), hence the sticky system
  tries to revert to equilibrium more strongly within the basin of
  attraction of individual equilibria, i.e.~as long as the
  perception gap does not become extreme. Fig.~\ref{fig:2DWSO} shows that
  the volatility decreases with $\rho$. For large (compared to $\rho$)
  deviations of $z_t$ from the set of equilibrium points, system
  \eqref{eqn:SA1} behaves as \eqref{eq0}.
 
 \subsection{Transitions between equilibrium states}

\begin{figure}[h!]
 	\centering
 	\begin{subfigure}{0.305\textwidth}
 		\includegraphics[width=\textwidth]{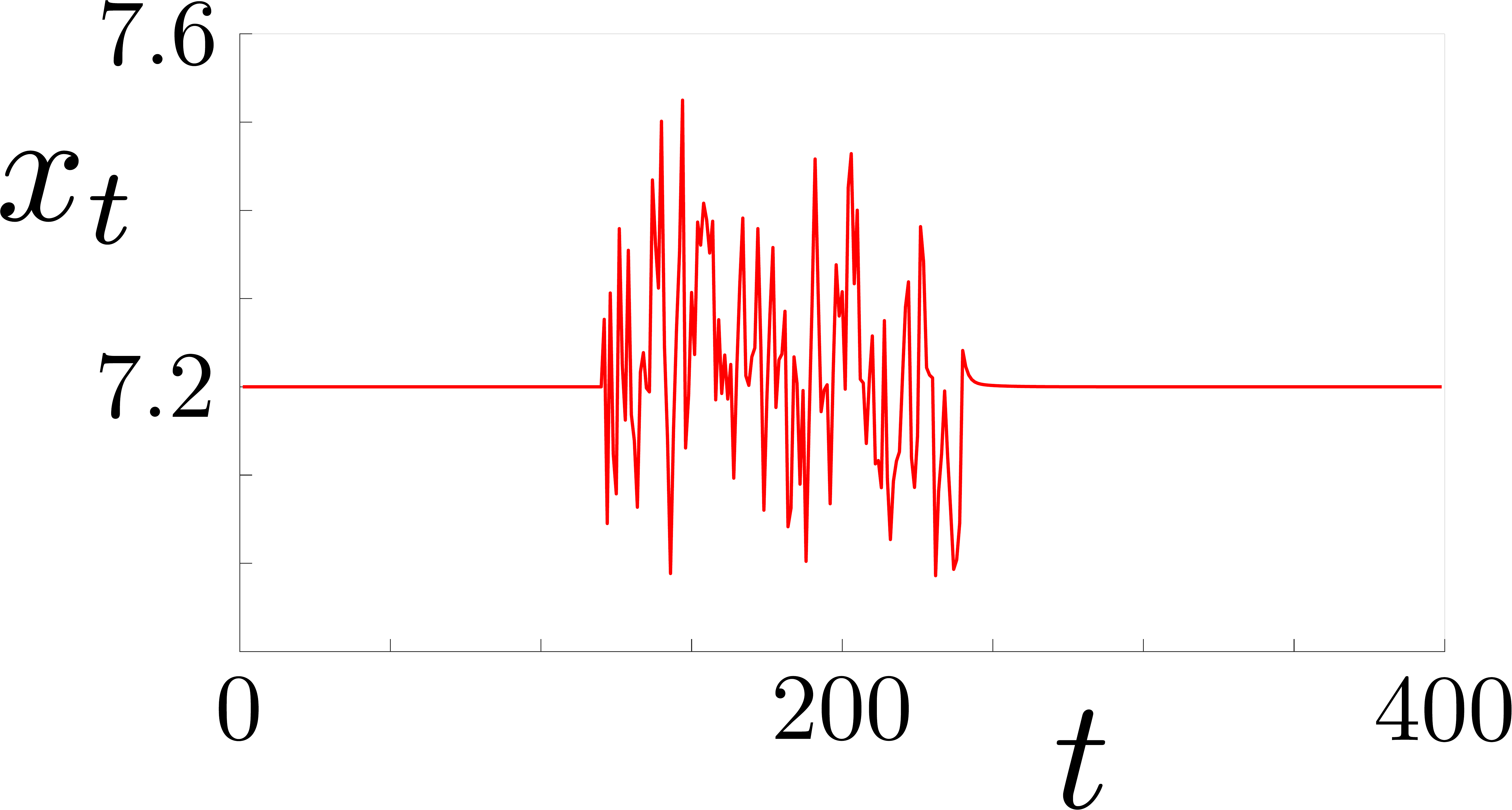}
 		\caption{ }
 		\label{fig:C1G1B1}
 	\end{subfigure}
 	\quad
 	\begin{subfigure}{0.305\textwidth}
 		\includegraphics[width=\textwidth]{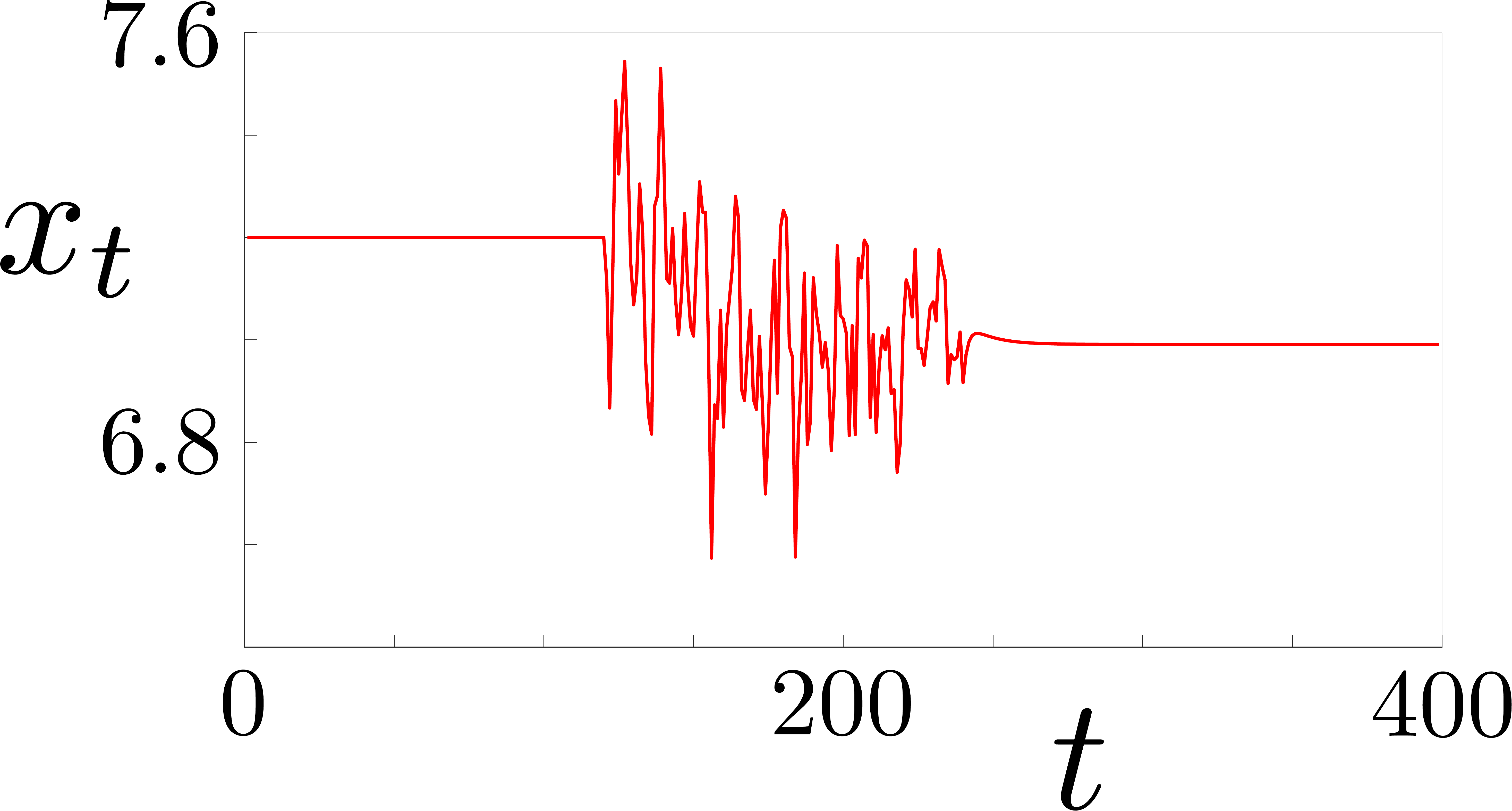}
 		\caption{}
 		\label{fig:C1G1B3}
 	\end{subfigure} 
 	\quad
 	\begin{subfigure}{0.305\textwidth}
 		\includegraphics[width=\textwidth]{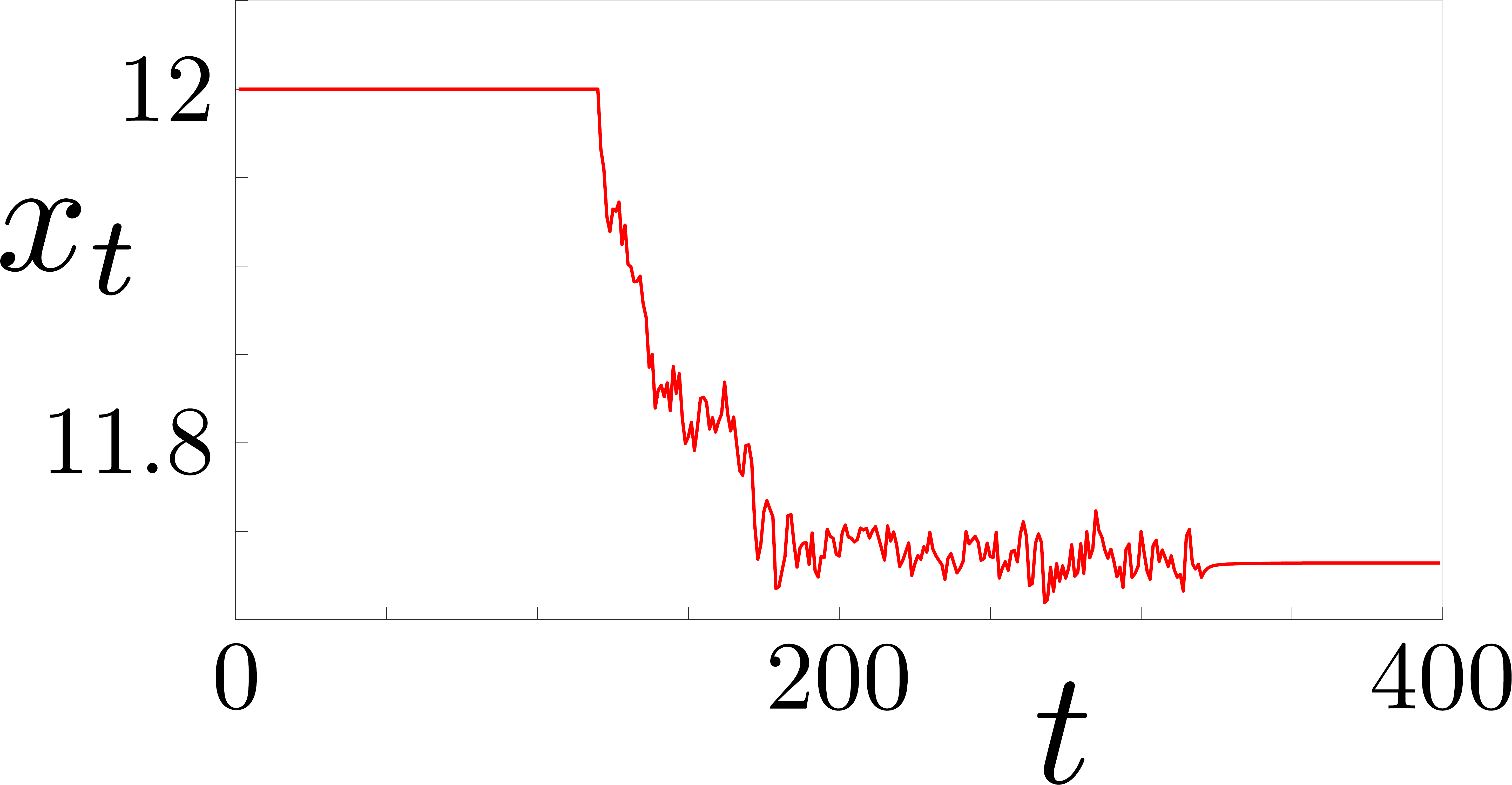}
 		\caption{}
 		\label{fig:C1G1c}
 	\end{subfigure}

 	 	\begin{subfigure}{0.305\textwidth}
 		\includegraphics[width=\textwidth]{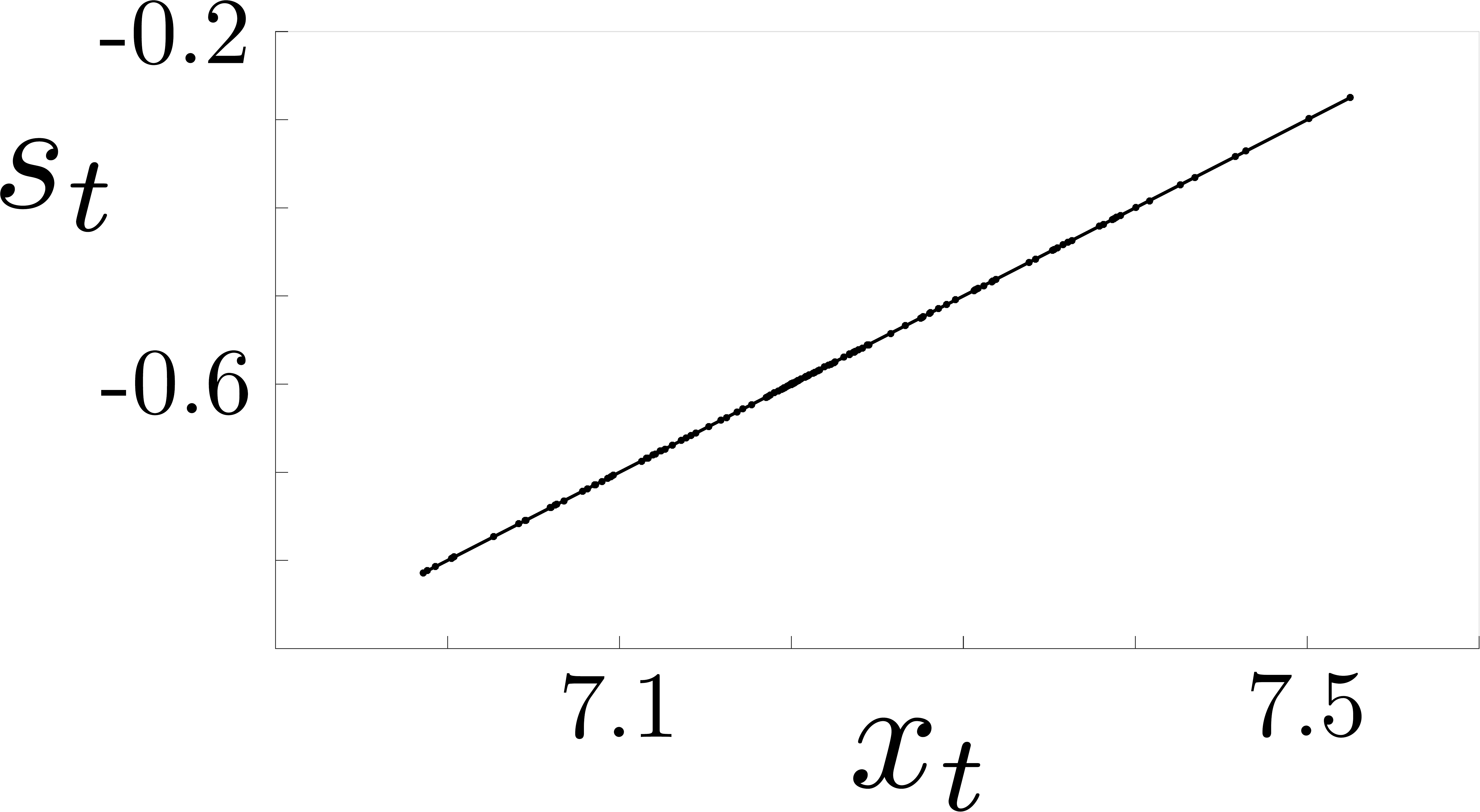}
 		\caption{}
 		\label{fig:C1G1B2}
 	\end{subfigure}   
 	\quad
 	\begin{subfigure}{0.305\textwidth}
 		\includegraphics[width=\textwidth]{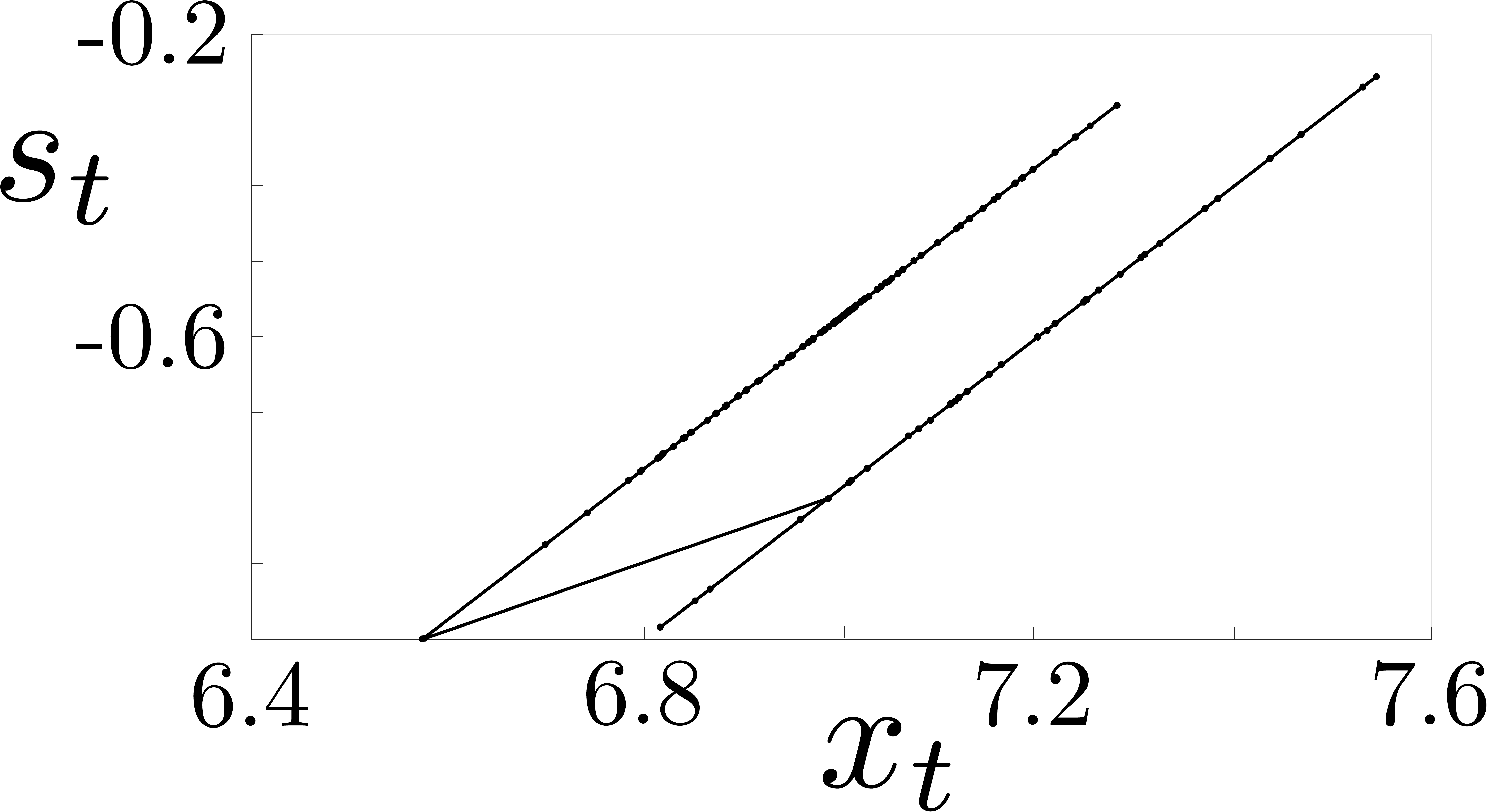}
 		\caption{}
 		\label{fig:C1G1B4}
 	\end{subfigure}
 	\quad
 	\begin{subfigure}{0.305\textwidth}
 		\includegraphics[width=\textwidth]{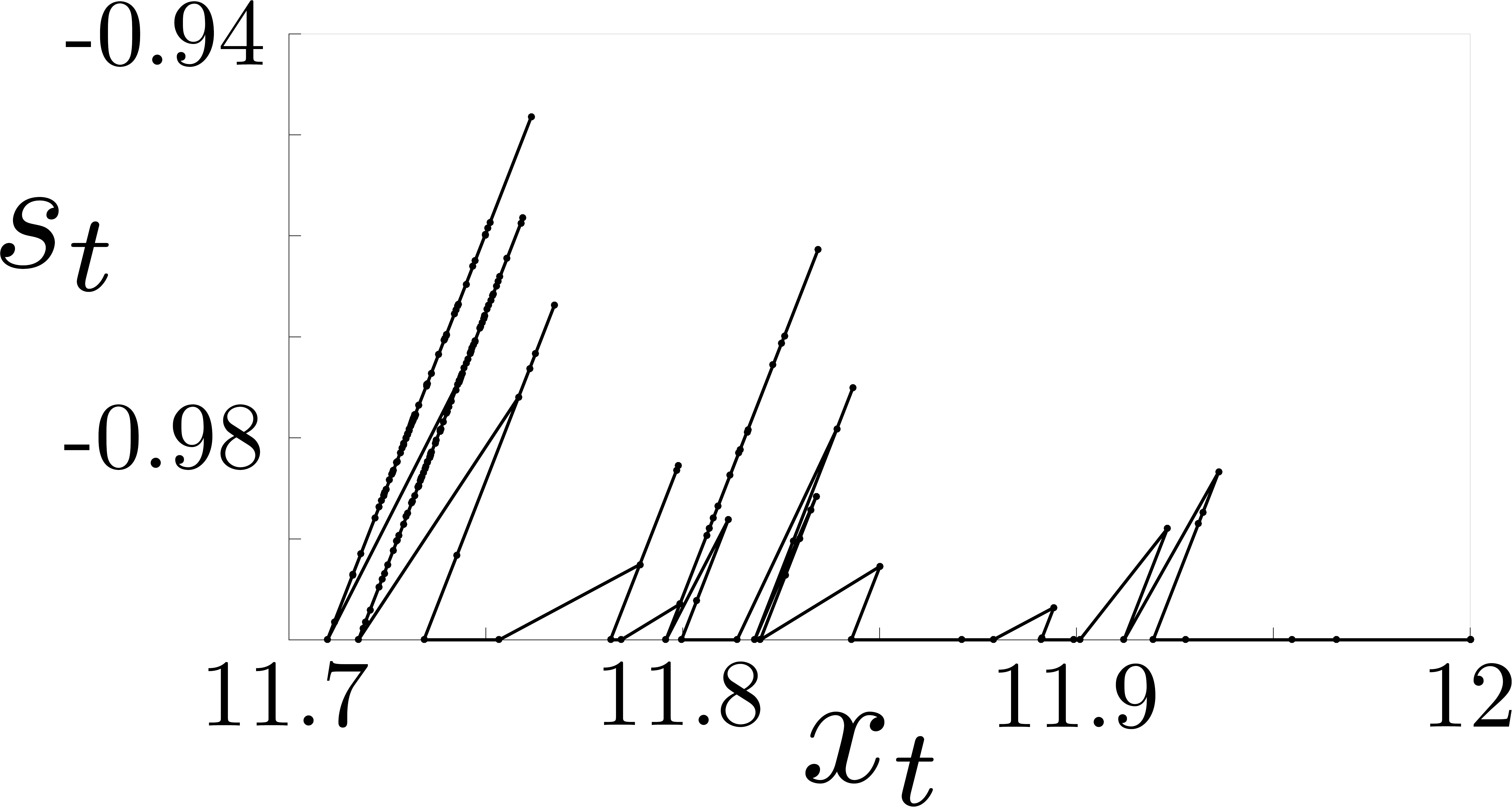}
 		\caption{} 
 		\label{fig:C1G1d}
 	\end{subfigure}
 	\caption{Transitions between the equilibrium states.   (a -- c) Time traces of inflation rate; (d -- f) the corresponding plots in the $(x,s)$-space exhibiting different 	transition scenarios.
	The noise is turned off before and after the interval of time of interest in order to show the equilibrium state at the ends of this interval.
	(a, d) The perception gap remains within the bounds $|s_t|<\rho$, and the system stays in the basin of attraction of one equilibrium point.
	The inflation rate $x_*(s_*)$ is the same before and after the noisy interlude. (b, e) The perception gap reaches the extreme value $-\rho$ (the highest expectation of inflation), and the trajectory transits
	from the basin of attraction of an equilibrium state with higher inflation rate and lower output gap (the right slanted segment in (e)) to the basin of attraction of an equilibrium state with a lower inflation rate and higher output gap (the left slanted
 	segment in (e)). (c, f). A transition from the equilibrium with the highest inflation rate (the rightmost point in (f))
	to an equilibrium state with a more moderate inflation rate through the basins of attraction of several other equilibrium states.
	}\label{fig:C1G1B}
 \end{figure}

  
      The system remains within the basin of attraction of a
      particular equilibrium state $z_*(s_*)$ as long as the
      perception gap $s_t$ does not reach either of the extreme values
      $\pm\rho$ and remains confined to the interval $|s_t|<\rho$, see
      Fig.~\ref{fig:C1G1B}(a,d).  But as soon as the perception gap
      hits the end of its range and starts being `dragged' by the
      actual inflation rate (Fig.~\ref{fig:C1G1B}(b,e)) the system
      transitions to the basin of attraction of a different
      equilibrium state where $s_t$ becomes `stuck' again.  For
      this reason, the system stays near equilibrium states which
      correspond to non-extreme perception gaps for longer periods of time
      than near extreme ones.  Figures~\ref{fig:C1G1B}(c,f) illustrate a
      transition from the equilibrium state with an extreme
      perception gap, $z_*(\rho)$, to one with a more moderate
      perception gap.
 
\subsection{Response to shocks}
We shall stress the system by applying supply shocks through the 
term $\eta_t$. The response of the system to demand shocks applied through the term $\epsilon_t$ is similar.  However, the parameter regime being considered  diminishes the effect of relatively small demand shocks due to the small value of $b_2=0.05$.
\begin{figure}[h!]
\centering
\begin{subfigure}{0.45\textwidth}
\includegraphics[width=\textwidth]{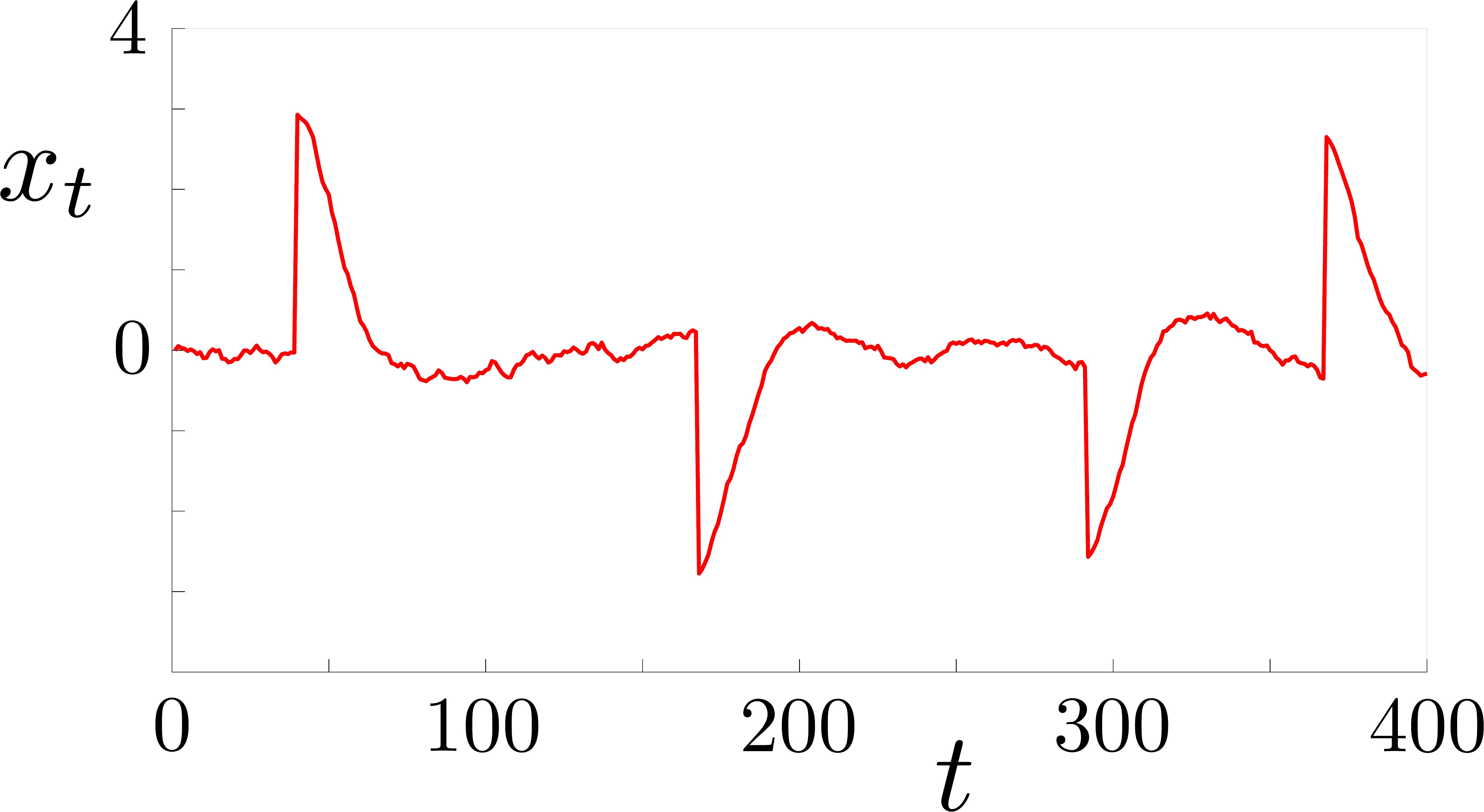}
\caption{ }
\label{fig:ShWoS}
\end{subfigure} 
\quad
\begin{subfigure}{0.45\textwidth}
\includegraphics[width=\textwidth]{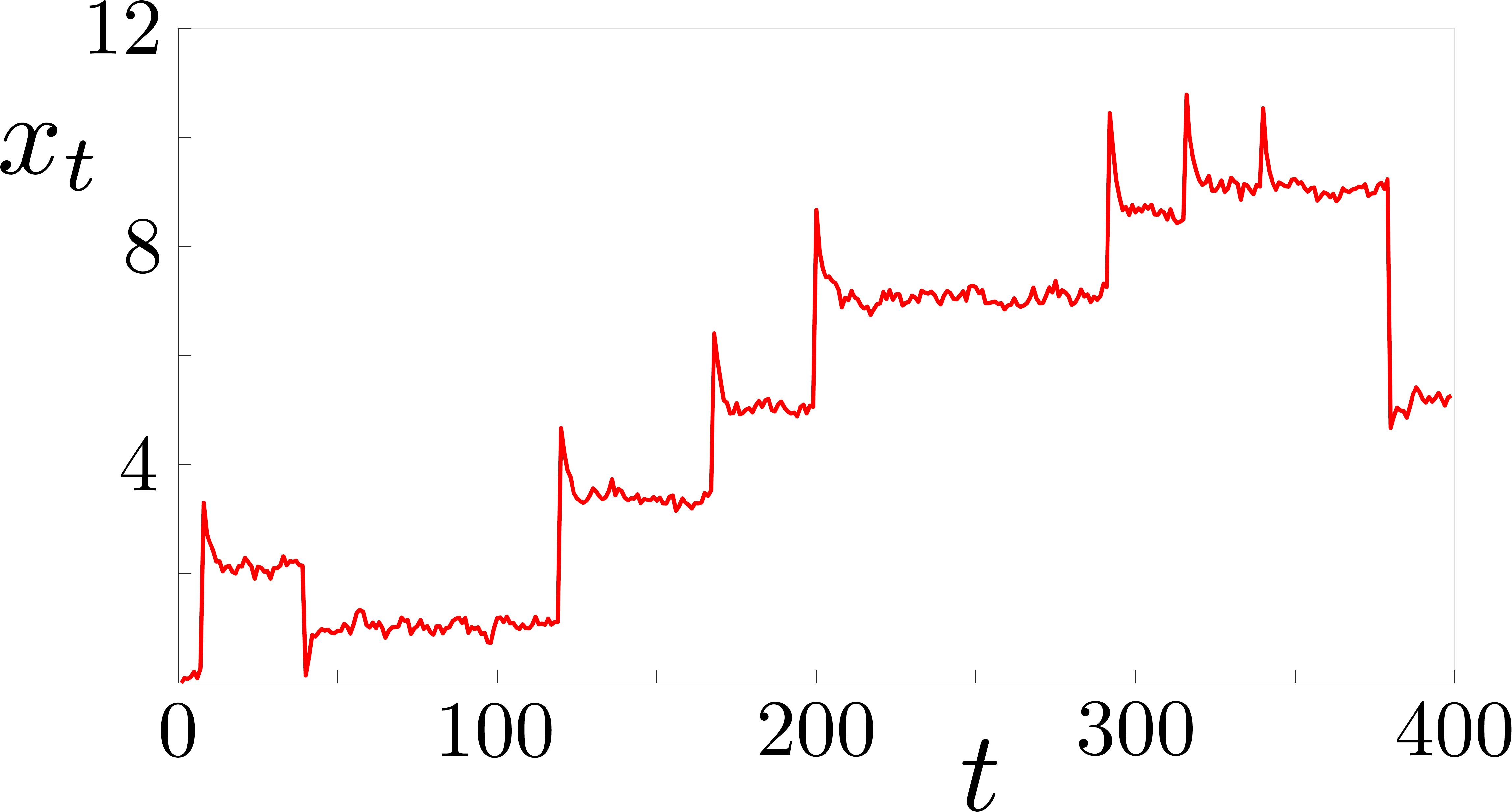}
\caption{}
\label{fig:ShWS}
\end{subfigure}
\caption{Response to shocks. (a) The system without stickiness ($\rho=0$) settles to the same unique equilibrium after each shock. (b) The system with stickiness ($\rho=1$)
settles to a new equilibrium after a shock is applied.}\label{fig:CompareShocks}
\end{figure}

\begin{figure}[h!]
\centering
\begin{subfigure}{0.45\textwidth}
\includegraphics[width=\textwidth]{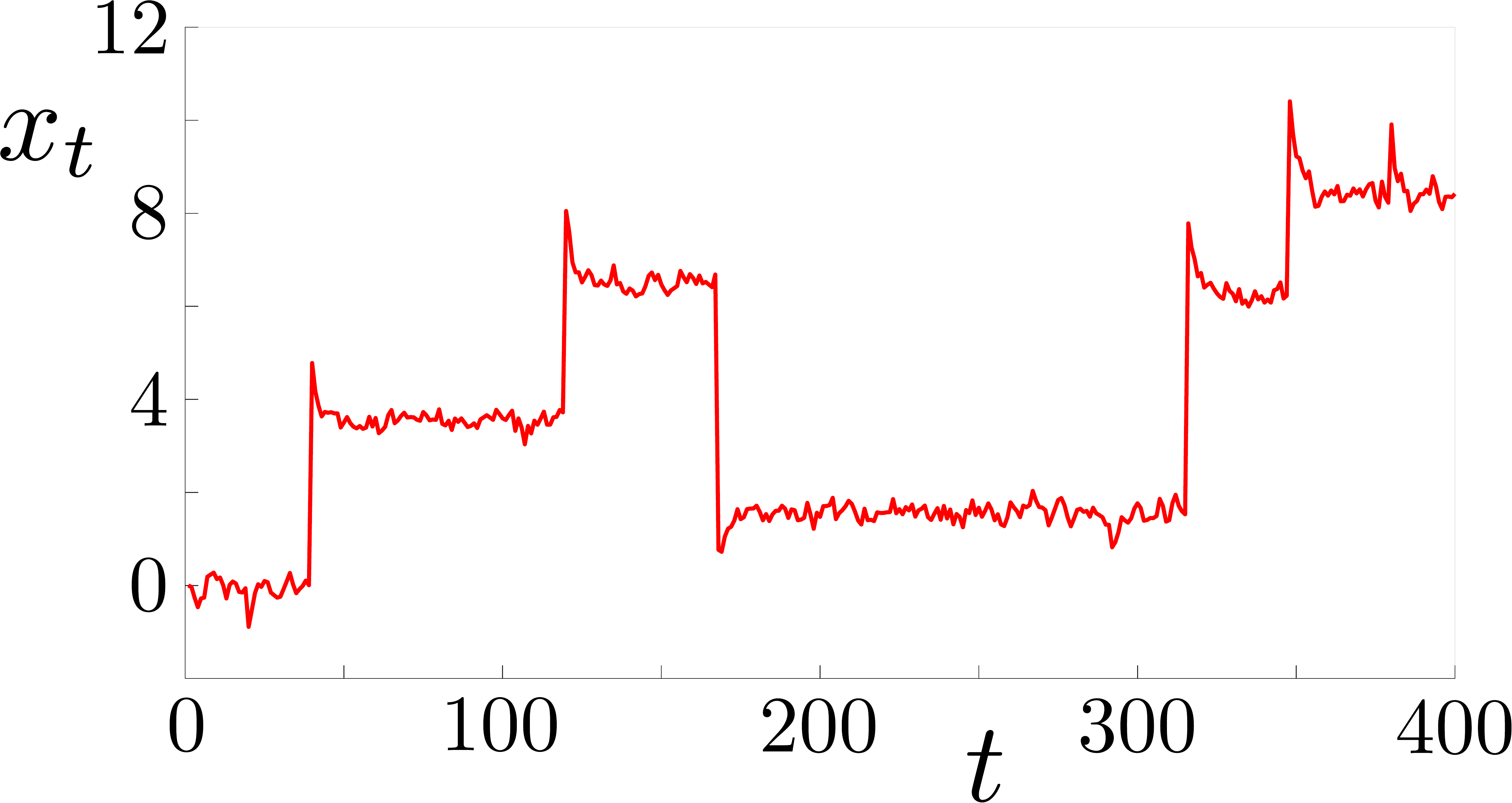}
\caption{ }
\label{fig:2DRS1sx}
\end{subfigure} 
\quad
\begin{subfigure}{0.45\textwidth}
\includegraphics[width=\textwidth]{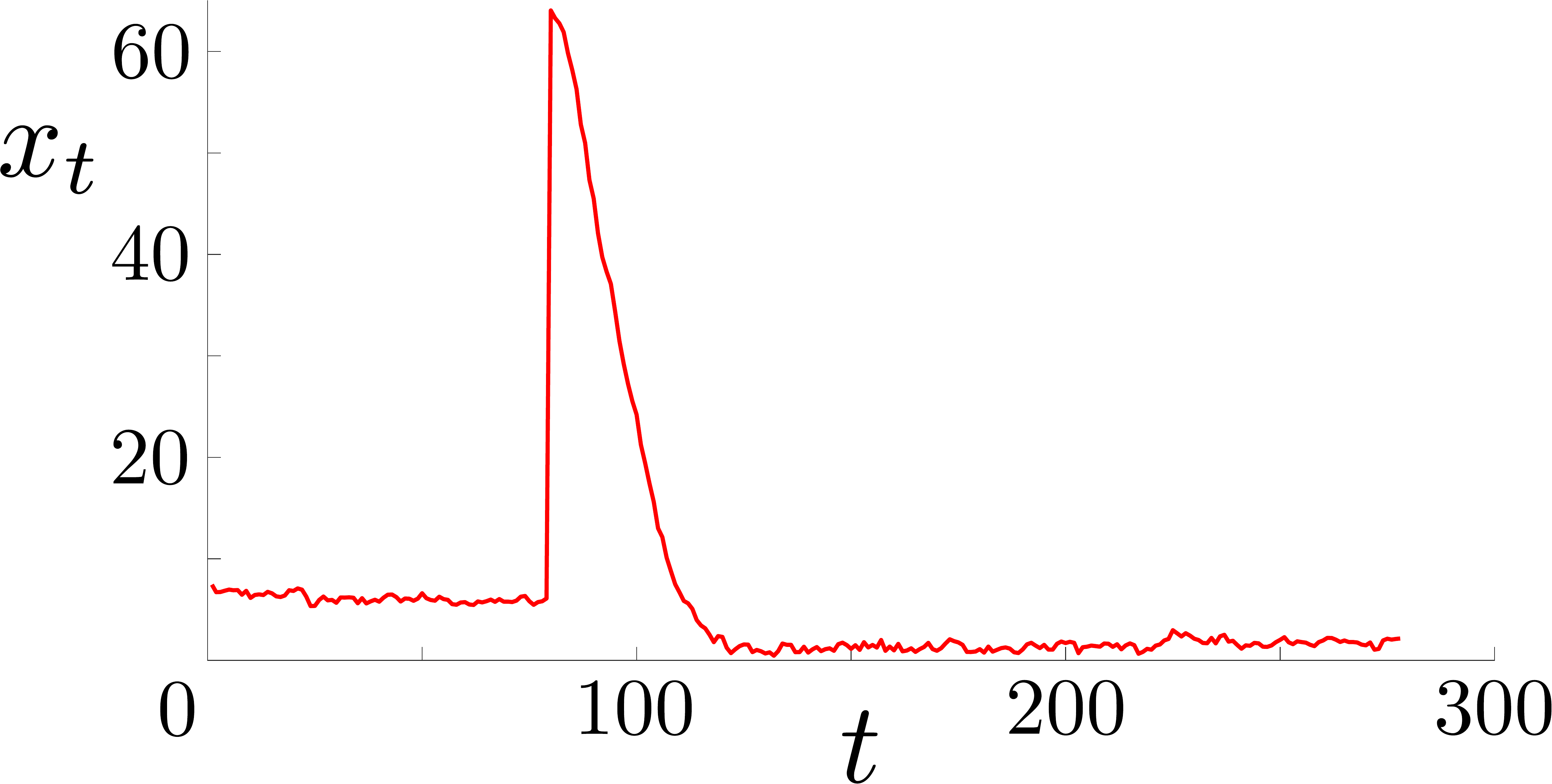}
\caption{}
\label{fig:2DRS1lx}
\end{subfigure}
\caption{Response to shocks of (a) small and (b) large magnitude. 
}\label{fig:2DRS1}
\end{figure}

System \eqref{eq0} without stickiness, which has a unique globally
stable equilibrium state $x_*=y_*= 0$, as expected returns to the
equilibrium (and hovers near it due to noise) after each shock,
see  Fig.~\ref{fig:CompareShocks}(a). 
Shocks applied to the sticky system \eqref{eqn:SA1}, \eqref{eqn:SA1'}
result in transitions between equilibrium states, see Figure~\ref{fig:CompareShocks}(b).
Numerical simulation show that shocks of small magnitude typically move the system in the direction of the shock (see Fig.~\ref{fig:2DRS1sx}). For example,
after a shock that pushes up the inflation rate the system settles to a new equilibrium state, which has higher inflation rate (and lower output gap) than the equilibrium occupied prior to the shock.
On the other hand, shocks of larger magnitude cause a transition to an
equilibrium state that can be hard to predict because such shocks
cause a longer and more complex excursion into the phase space far from equilibrium set.
In Fig.~\ref{fig:2DRS1lx}, the system resides near an equilibrium with high inflation rate before a shock is applied. Although the shock pushes the inflation even higher,
the system eventually settles to an equilibrium with nearly zero inflation rate after the shock is removed.

\subsection{The possibility of runaway inflation}\label{infln}

\begin{figure}[h!]
\centering
\begin{subfigure}{0.45\textwidth}
\includegraphics[width=\textwidth]{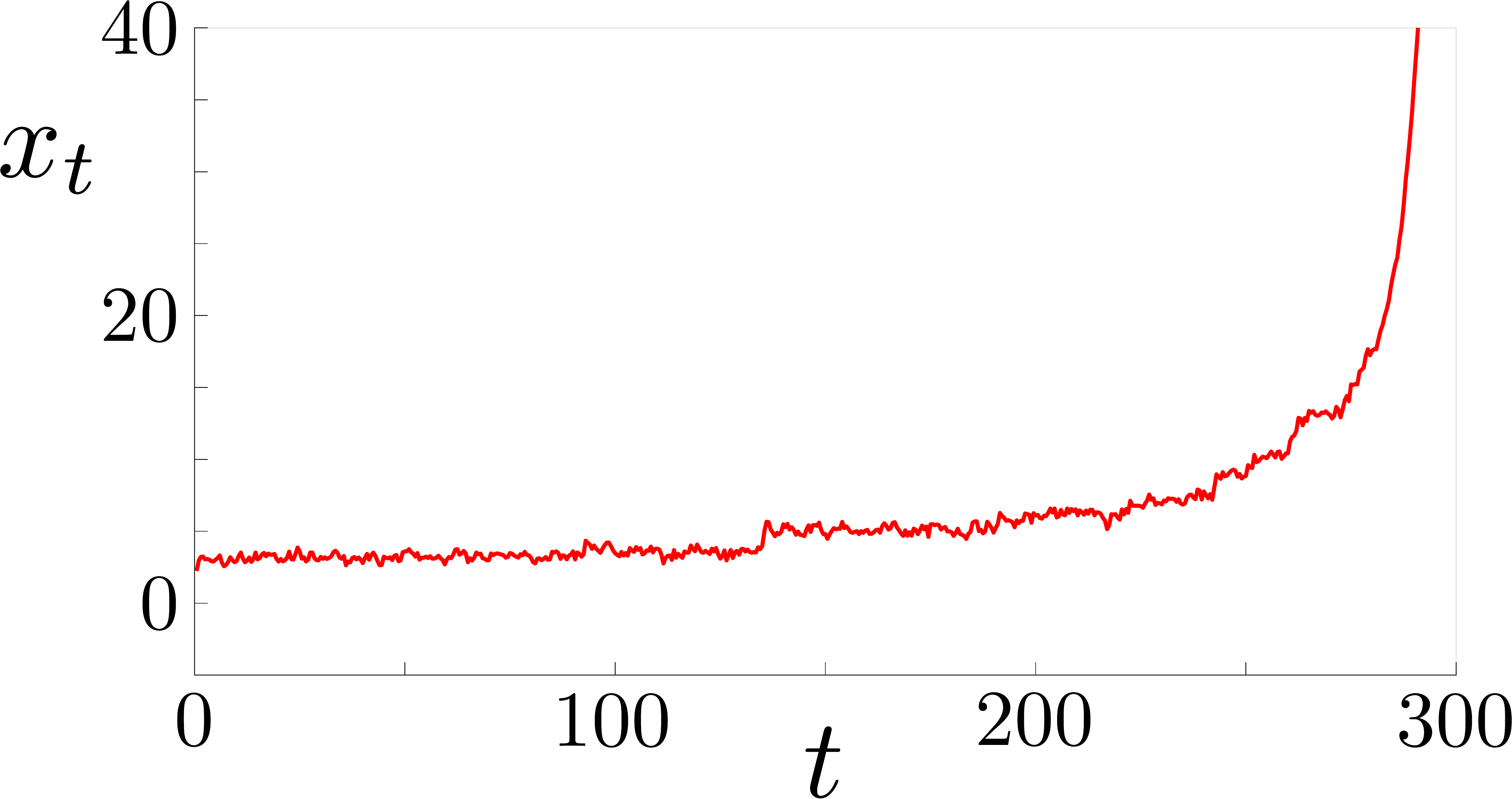}
\caption{}
\label{fig:C1L1c}
\end{subfigure}
\quad
\begin{subfigure}{0.45\textwidth}
\includegraphics[width=\textwidth]{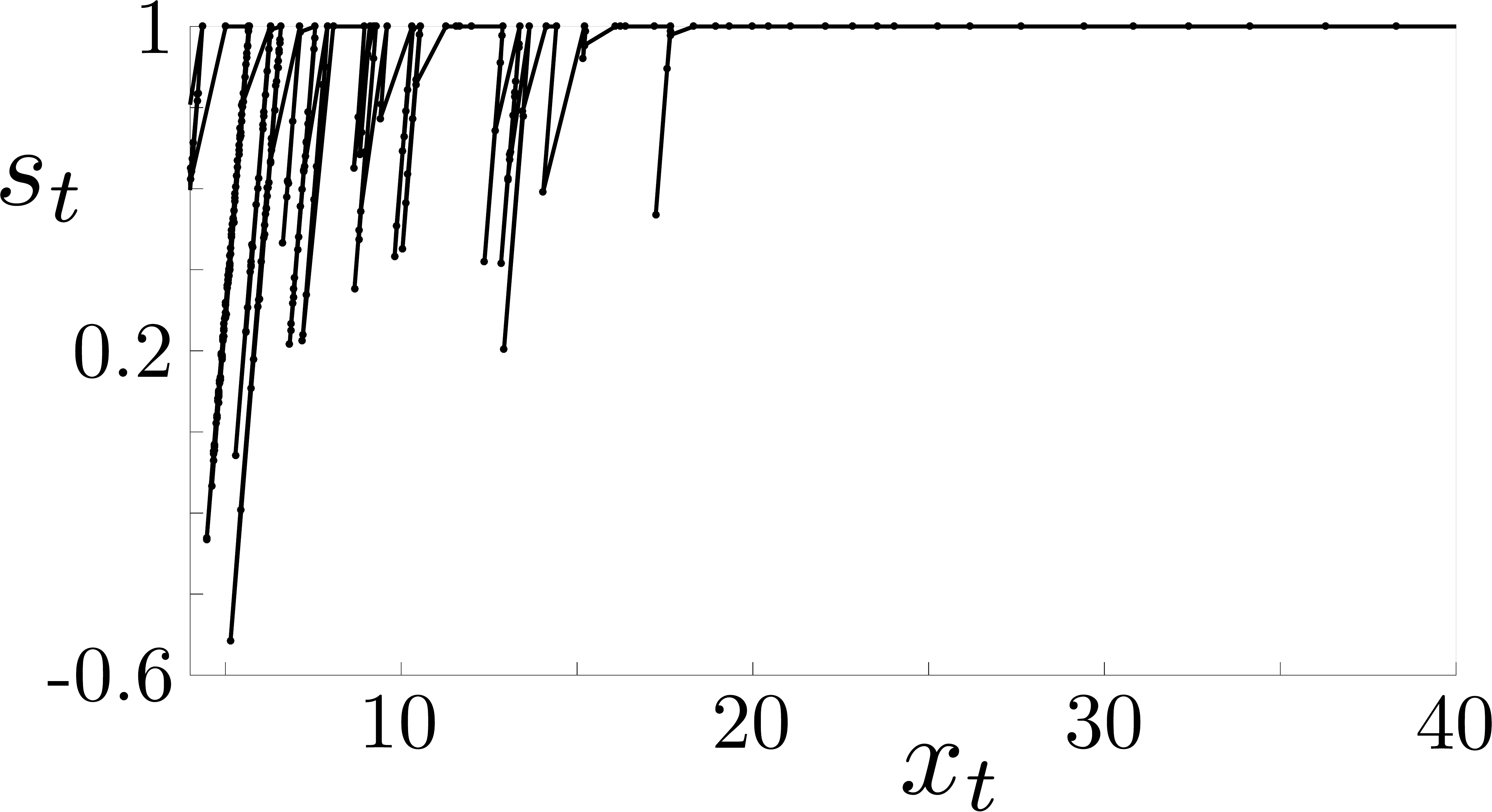}
\caption{}
\label{fig:C1L1d}
\end{subfigure}
\caption{Run-away inflation scenario. 
Parameter are $\rho=1$, $a=0.3$, $b_1=0.5$, $b_2=0.05$, $c_1=0.9$, $c_2=0.01$. The ranges of inflation rate and output gap values at equilibrium states  for this set parameter  are $x_*\in[-11,11]$ and $y_*\in[-10,10]$, respectively. 
(a) Time series of inflation rate $x_t$.   
(b) Trajectory in the  $(x,s)$ space.}
\label{fig:C1L1B}
\end{figure}

According to Section \ref{staglo} the system is globally stable for $c_1>1$, but becomes unstable for $c_1<1$.
The latter case creates a possibility of the  run-away inflation scenario. 
It is interesting that as shown in Section \ref{staloc} all the equilibrium points are {\em locally} stable even if $c_1<1$.
As a result, dynamics appear to be stable as long as the trajectory is confined to the basin of attraction of an equilibrium state.
However, when noise or a shock or another fluctuation drives the trajectory outside this bounded stability domain,
the run-away scenario may and is likely to start, see Fig.~\ref{fig:C1L1B}.
Just to be clear, the behavior is stable while the perception gap is
not extreme, but if a shock causes that to change then the runaway instability
can suddenly occur with no change in the system parameters.

\subsection{A trade-off between inflation and output gap volatility}\label{numerics6}

\begin{figure}[h!]
\centering
\begin{subfigure}{0.48\textwidth}
\includegraphics[width=\textwidth]{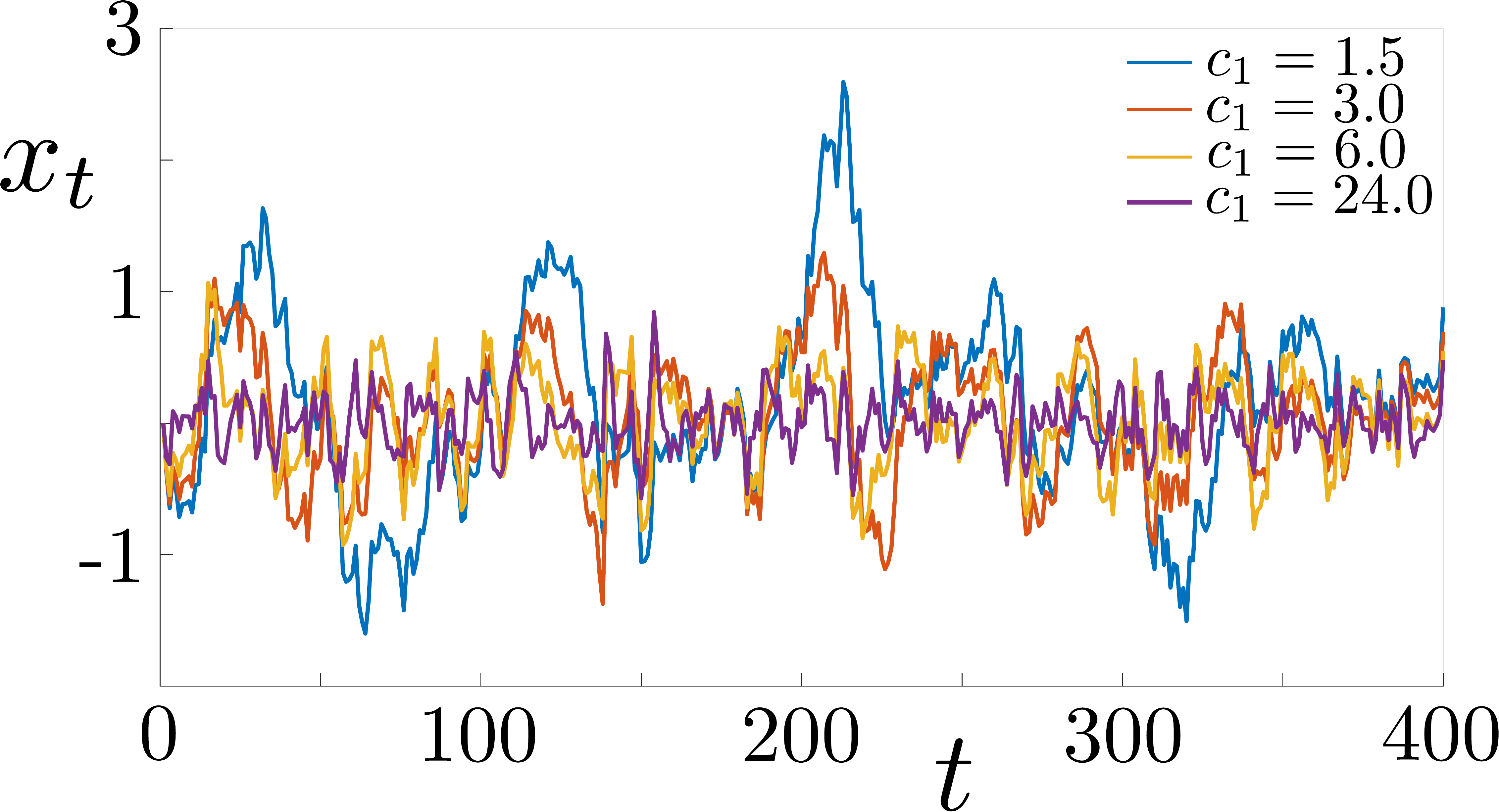}
\caption{} 
\label{fig:C1Inc3a}
\end{subfigure}
\quad
\begin{subfigure}{0.48\textwidth}
\includegraphics[width=\textwidth]{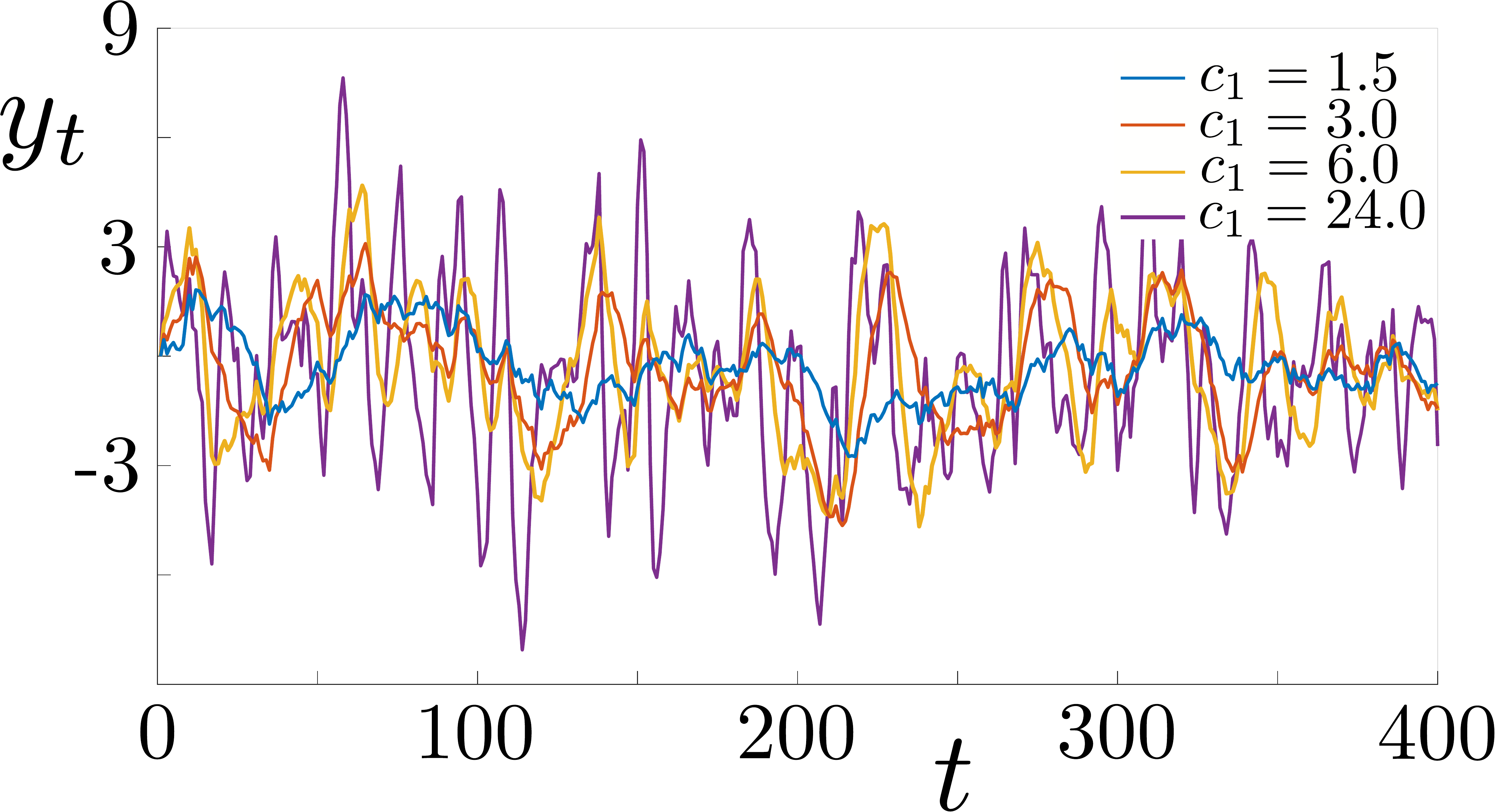}
\caption{}
\label{fig:C1Inc3b}
\end{subfigure}        
\caption{Numerical simulations of (a) inflation rate, $x_t$ and (b) output gap, $y_t$ for $\rho=1$ and various values of $c_1$. The remaining parameters values are from Table \ref{tab:DS1}. }\label{fig:C1Inc3}
\end{figure}

Parameters $c_1$ and $c_2$  of Taylor's rule \eqref{eqn:M1'} control the volatility level of inflation and output gap near an equilibrium state. 
Numerical simulations of the model with sticky inflation expectation show that
when $c_1$ increases (which corresponds to stronger inflation targeting by the Central Bank), 
the volatility of the inflation rate decreases, see Fig.~\ref{fig:C1Inc3a}. However, at the same time, the output gap becomes highly volatile with increasing $c_1$, see Fig.~\ref{fig:C1Inc3b}. 

\begin{figure}[h!]
\centering
\begin{subfigure}{0.48\textwidth}
\includegraphics[width=\textwidth]{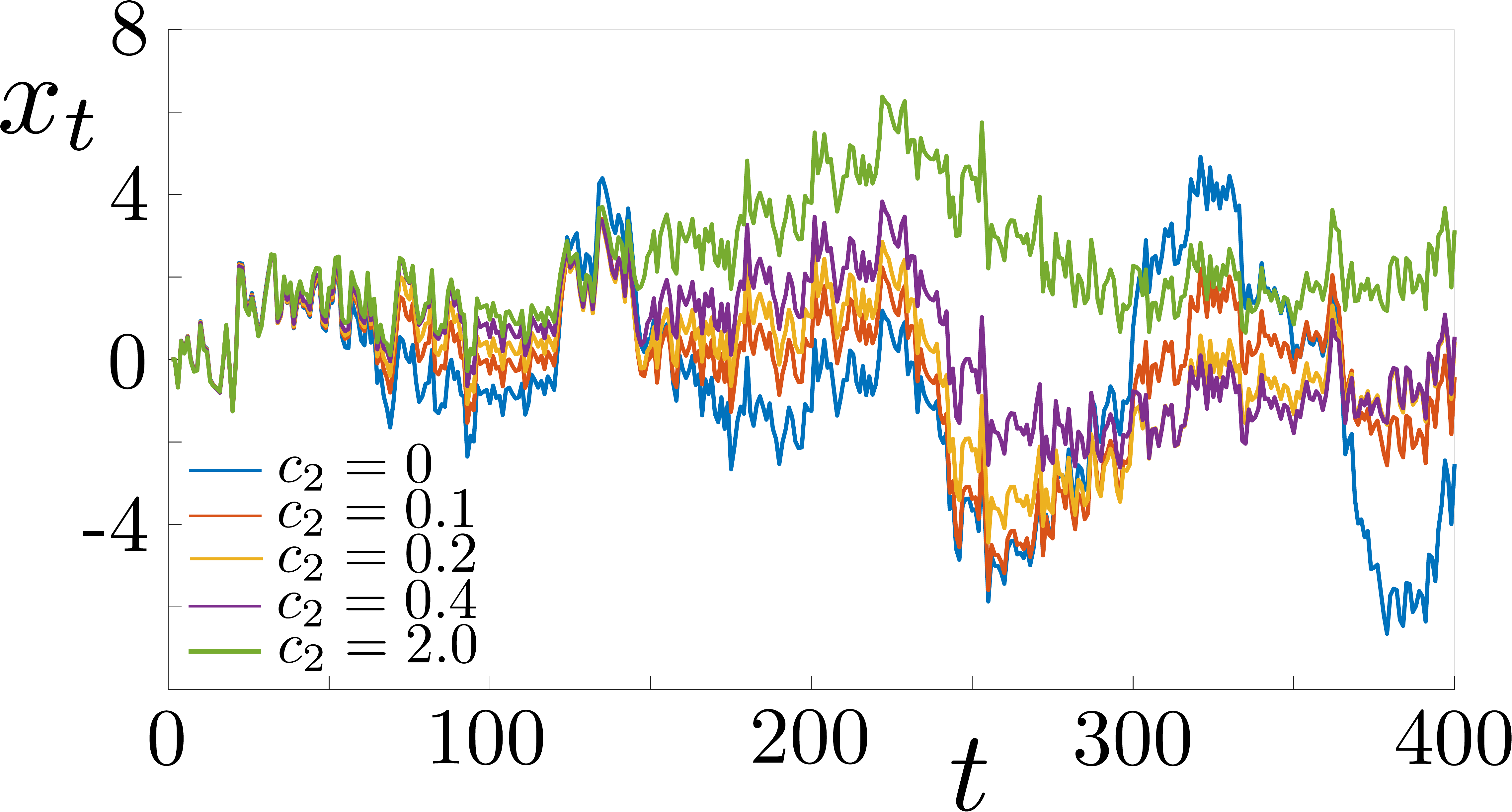}
\caption{ } 
\label{fig:C2Inc3a}
\end{subfigure}
\quad
\begin{subfigure}{0.48\textwidth}
\includegraphics[width=\textwidth]{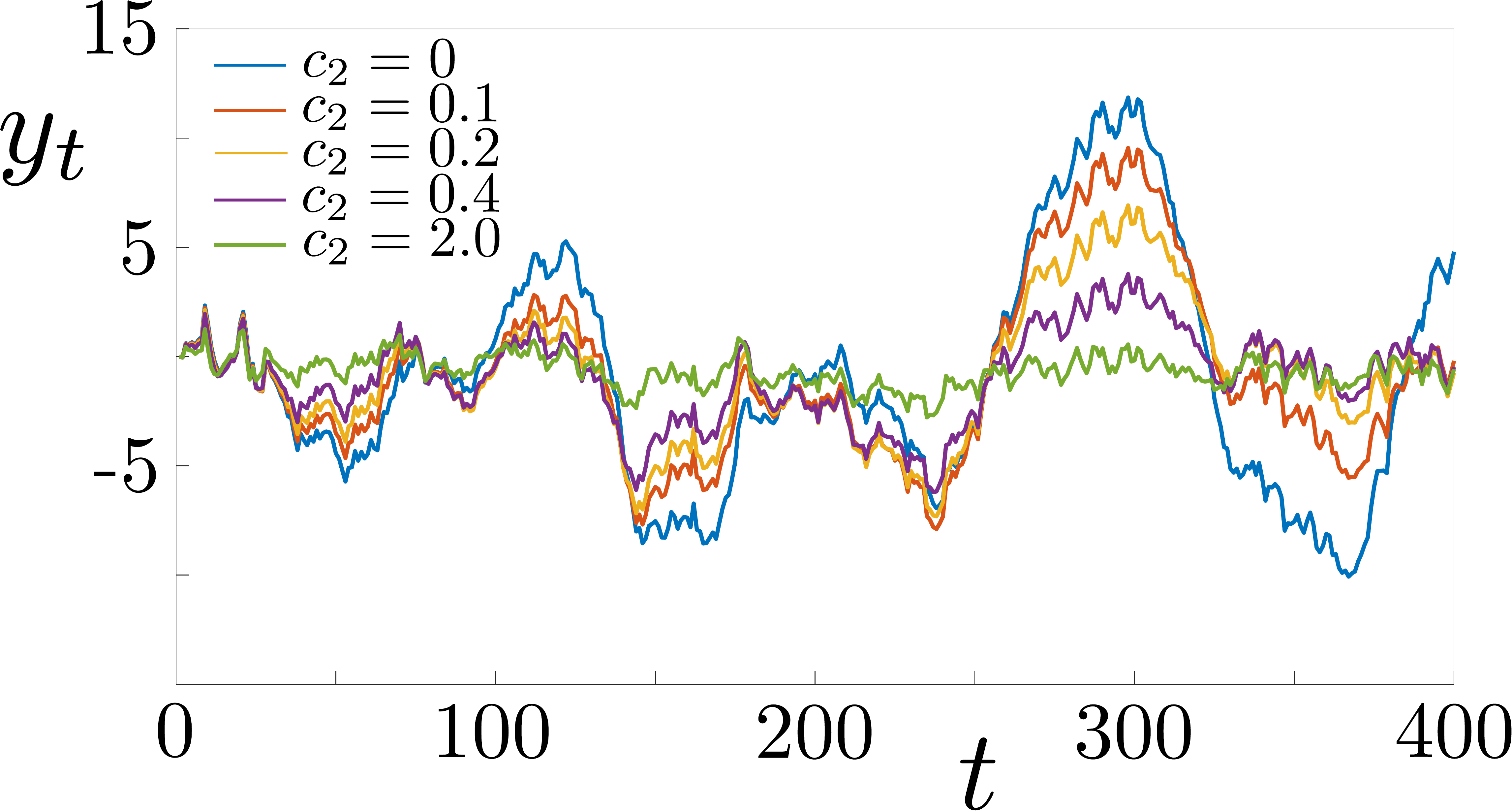}
\caption{}
\label{fig:C2Inc3b}
\end{subfigure}           
\caption{Numerical simulations of  (a) inflation rate, $x_t$ and (b) output gap, $y_t$ for $\rho=1$ and various values of $c_2$. The remaining parameter values are from Table \ref{tab:DS1}. }\label{fig:C2Inc3A}
\end{figure}

\begin{figure}[h!]
\centering
\begin{subfigure}{0.48\textwidth}
\includegraphics[width=\textwidth]{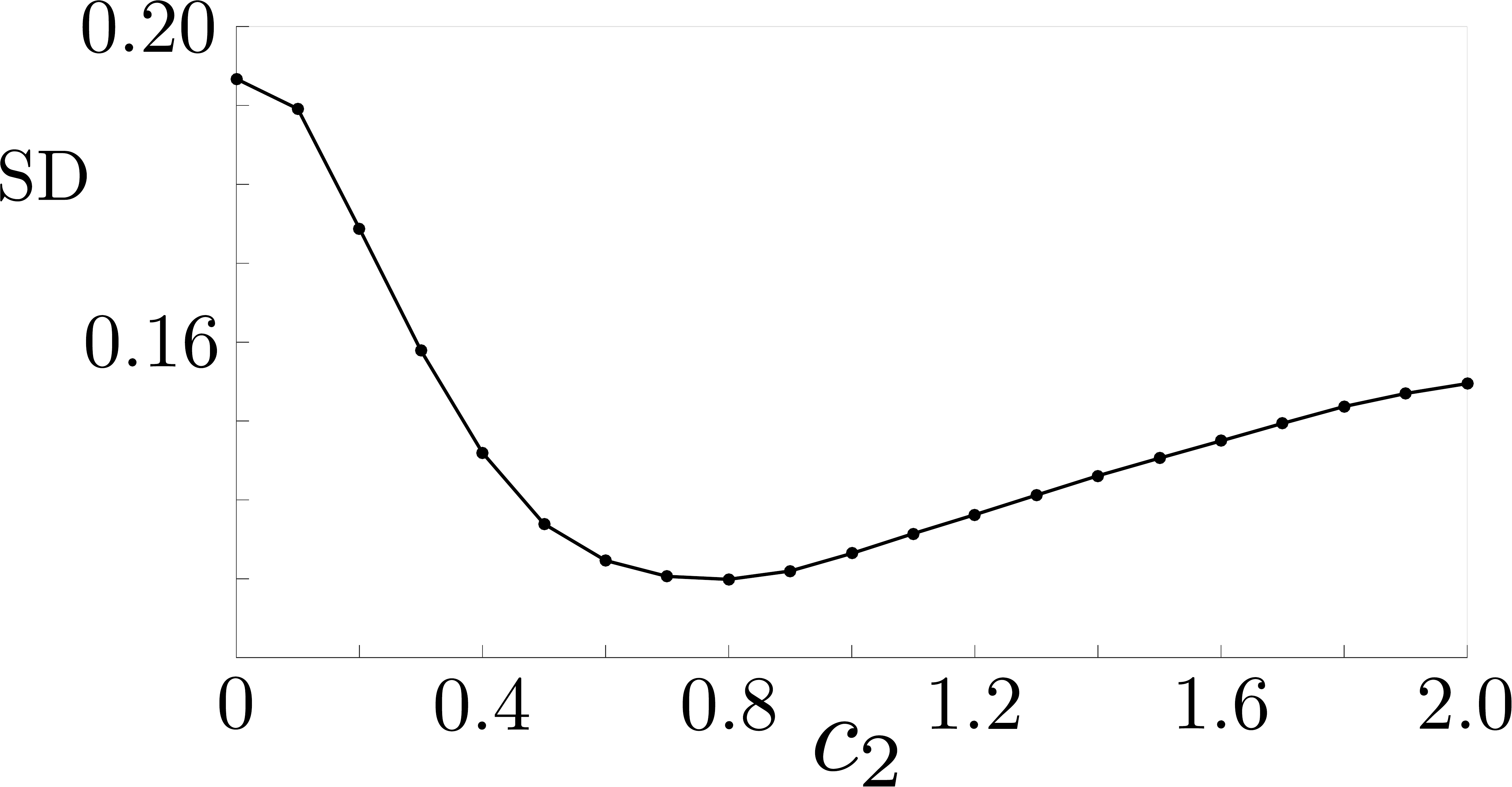}
\caption{ }
\label{fig:C2STDInc3a}
\end{subfigure}
\quad
\begin{subfigure}{0.48\textwidth}
\includegraphics[width=\textwidth]{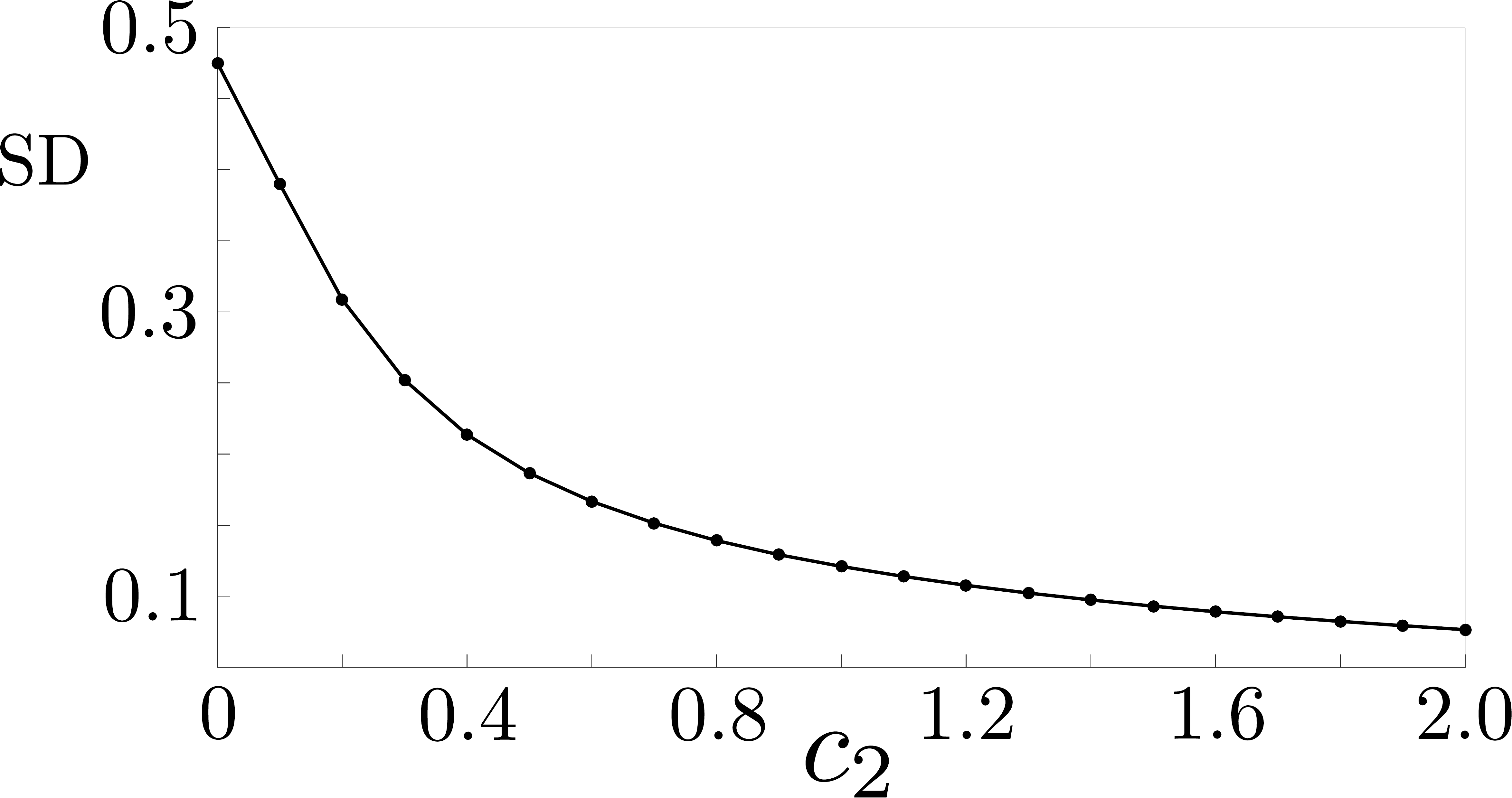}
\caption{ }
\label{fig:C2STDInc3b}
\end{subfigure}           
\caption{Measure of the effect of $c_2$ on volatility of (a) $x_t$ and (b) $y_t$ with standard deviation (SD).
	}\label{fig:C2STDInc3A}
\end{figure}

When $c_2$ increases (stronger output gap targeting), the output gap volatility decreases, see Fig.~\ref{fig:C2Inc3b}. 
In particular, the case $c_2=0$ corresponding to pure inflation targeting in Taylor's rule is characterized by the highest volatility of the output gap.
However, from Fig.~\ref{fig:C2Inc3a}, it appears that the inflation rate volatility exhibits a non-monotone behavior with $c_2$.
This is confirmed by Fig.~\ref{fig:C2STDInc3A}, which shows the dependence of the standard deviation of $x_t$ and $y_t$ on $c_2$ for the
trajectories presented in  Fig.~\ref{fig:C2Inc3A}.  The inflation rate
volatility reaches its minimum for $c_2\approx 0.8$ for the parameter values $a,b_1,b_2,c_1$ from Table \ref{tab:DS1} and $\rho=1$. 

All the above results are in agreement with \cite{DeGrauwe2012}. In addition, $c_1$ and $c_2$ affect the range of the inflation rate value
at the equilibrium states
for the model \eqref{eqn:SA1}. According to \eqref{eqn:SA1E}, this range increases with $c_2$ and decreases with $c_1-1$ (for $c_1>1$).
At the same time, the range of output gap equilibrium values is unaffected by the parameters of Taylor's rule.

\subsection{A multi-agent model}\label{multi}
Model \eqref{eqn:SA1} can be easily extended to account for 
differing types of agent with different inflation rate expectation rules/thresholds.
To this end, we replace the simple relationship \eqref{playstop}
between $p_t$ and $x_t$ with the equation
\begin{equation}\label{eqn:MStop}
p_t=\sum_{i=1}^n\mu_i\mathscr{P}_{\rho_i}[x_t]=x_t-\sum_{i=1}^n\mu_i\mathscr{S}_{\rho_i}[x_t]
\end{equation}
with 
\begin{equation}\label{eqn:MStop'}
\sum_{i=1}^n \mu_i=1.
\end{equation}
Here the play operator $\mathscr{P}_{\rho_i}$ models the expectation
of inflation by the $i$-th agent; $p_t$ is the aggregate expectation of inflation; 
$\mu_i>0$ is a weight measuring the contribution of agent's expectation 
of inflation to the aggregate quantity; and, $\rho_i$ is an individual threshold characterizing
the behavior of the $i$-th agent. 
Relation \eqref{eqn:MStop} is equivalent to the formula
\begin{equation}\label{eqn:MStop''}
s_t={\mathcal I}[x_t]:=\sum_{i=1}^n\mu_i\mathscr{S}_{\rho_i}[x_t],
\end{equation}
which is a (discrete) Prandtl-Ishlinskii (PI) operator with thresholds $\rho_i$ and weights $\mu_i$ \cite{AYIshlinskiiStatMethods,prandtl1928gedankenmodell,KP},
where $s_t=x_t-p_t$.

The implicit system \eqref{eqn:M1}, \eqref{eqn:M1'}, \eqref{eqn:MStop} with multiple agents can be converted into an explicit form
using the same technique as we used for the system with one play
operator. Again this involves the inversion of the PI operator.
The explicit system
\begin{equation}\label{eqn:SA111}
z_t={A}z_{t-1}+\mathcal{\hat I} [c\cdot z_{t-1}+\hat\xi_t]\,d+N\xi_t,
\end{equation}
which is similar to its counterpart \eqref{eqn:SA1}, includes 
a PI operator with rescaled thresholds $\hat \rho_i$ and weights $\hat \mu_i$,
see Appendix D for details; $\xi_t$, $\hat \xi_t$ denote the noise terms.

\begin{figure}[h!]
\centering
\begin{subfigure}{0.45\textwidth}
\includegraphics[width=\textwidth]{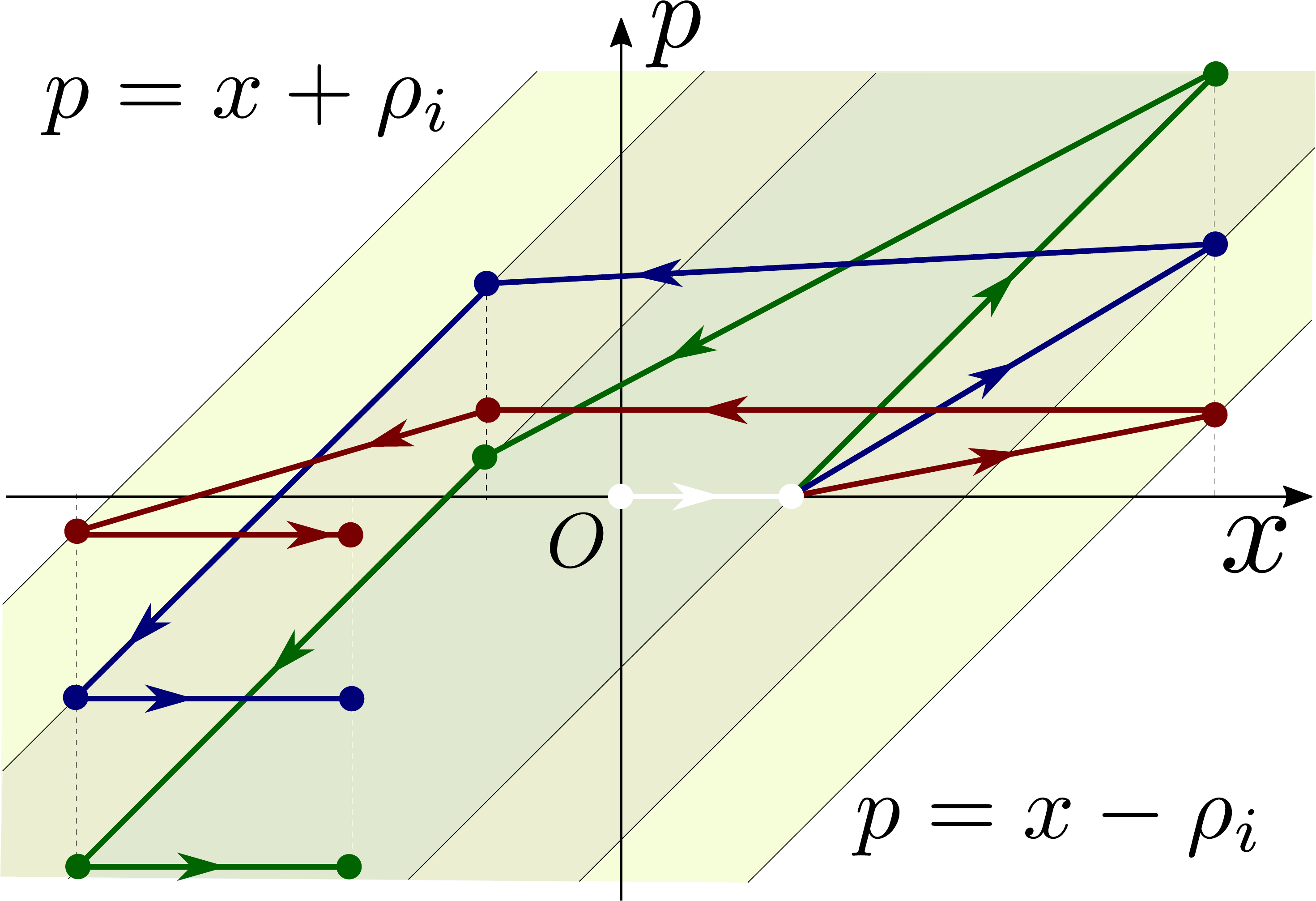}
\caption{}
\label{fig:PlayO3A}
\end{subfigure} 
\quad
\begin{subfigure}{0.44\textwidth}
\includegraphics[width=\textwidth]{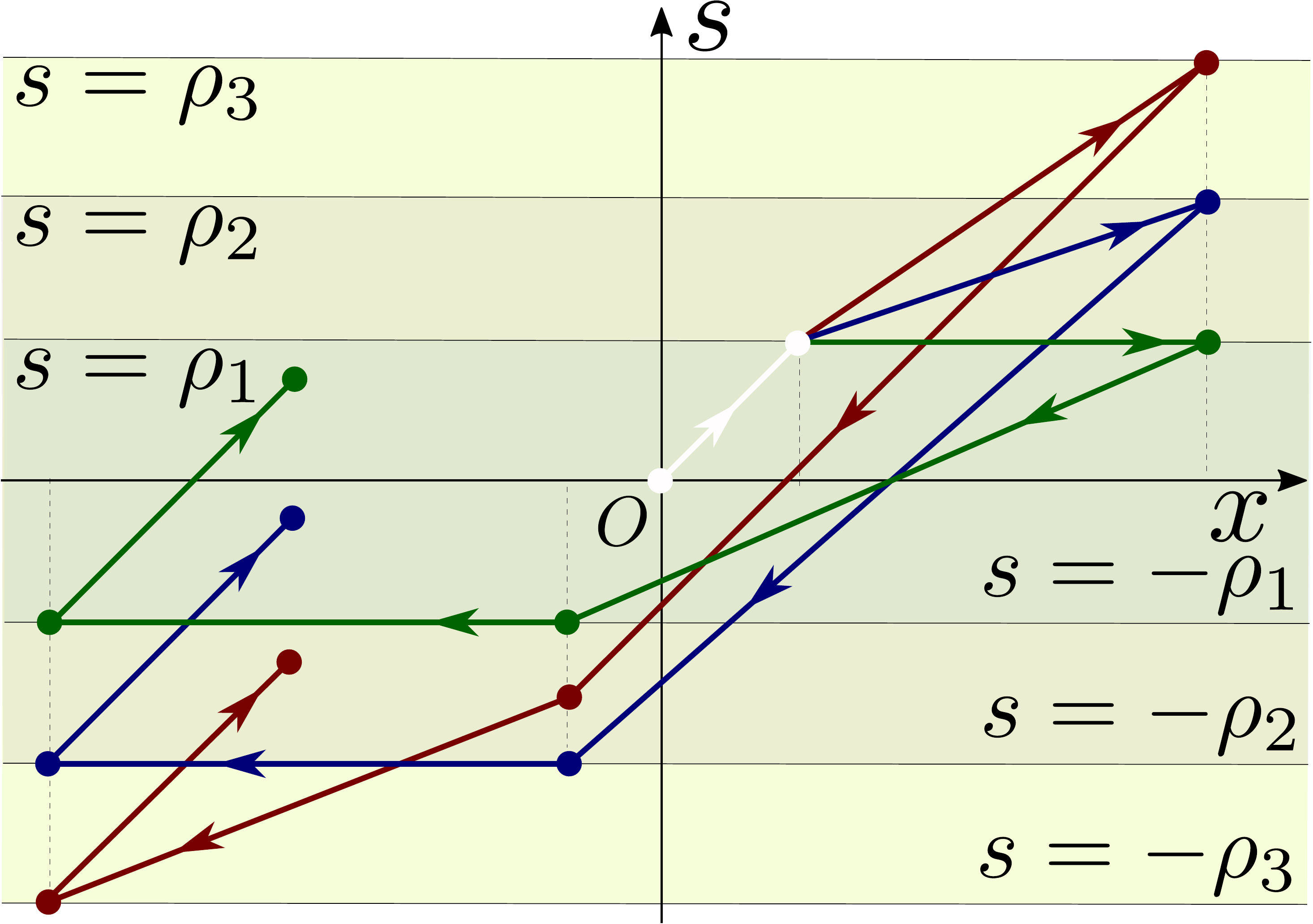}
\caption{}
\label{fig:StopO3A}
\end{subfigure}
\caption{Different expectations of agents based on three thresholds
 $\rho_1<\rho_2<\rho_3$ of (a) play and (b) stop operators   with a single input $x_t$.
 }\label{fig:HOperators3A}
\end{figure}

The stability properties of the equilibrium states of system \eqref{eqn:SA111} with multiple agents are similar to the stability properties considered above in Section 
\ref{staglo}. In particular, if we consider the system without external noise for $c_1>1$, then the set of equilibrium states is globally stable, and every trajectory converges to an equilibrium state.

In the simulations of this section, we 
classify economic agents into three categories, strongly, moderately, and weakly sensitive to inflation rate variations
(hence $n=3$), by assigning thresholds $\rho_1<\rho_2<\rho_3$, respectively, to these groups, see Fig.~\ref{fig:HOperators3A}.  Further, the contribution of each group
to the aggregate expectation of inflation carries equal weight, $\mu_i=1/3$.

\begin{figure}[h!]
	\centering
	\begin{subfigure}{0.48\textwidth}
		\includegraphics[width=\textwidth]{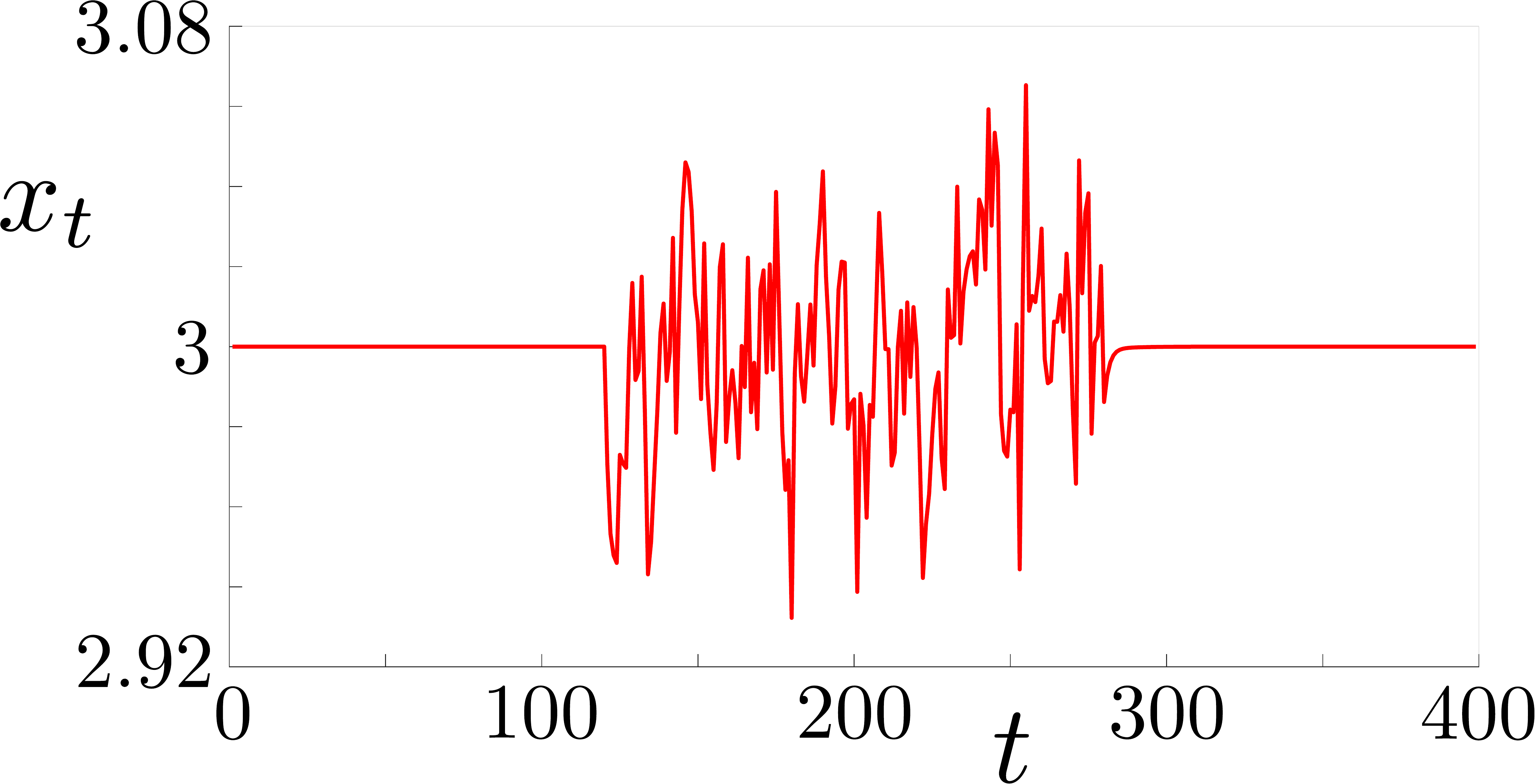}
		\caption{
			}\label{fig:3SC1G1a}
	\end{subfigure}
	\quad
	\begin{subfigure}{0.48\textwidth}
		\includegraphics[width=\textwidth]{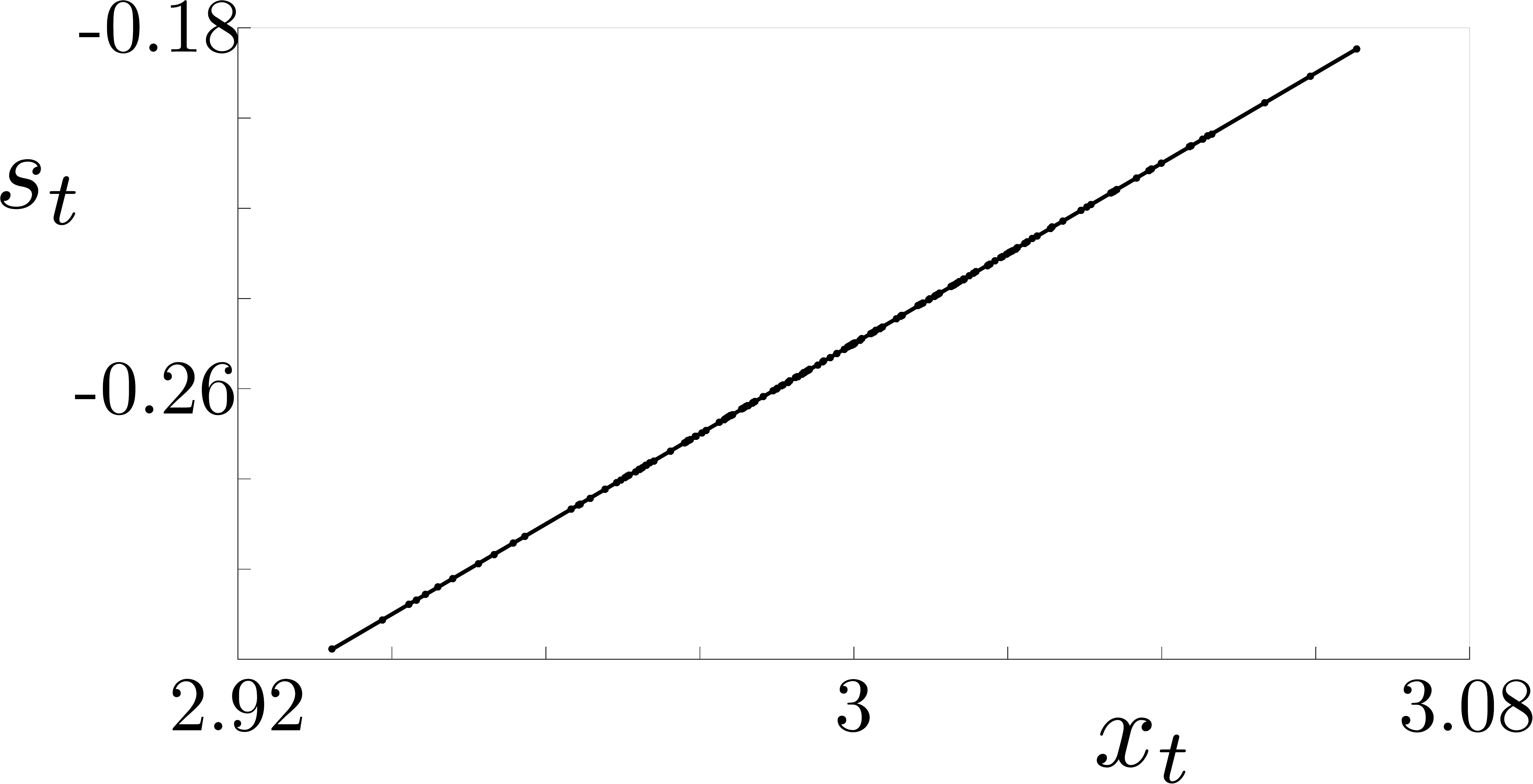}
		\caption{
			}\label{fig:3SC1G1b}
	\end{subfigure}
	\caption{
		Trajectory of the system with 3 agents near an equilibrium state when none of the agents achieves an extreme perception gap (cf.~Figure \ref{fig:C1G1B}(a, d)). Here $c_1>1$. (a) Time trace of inflation. (b) Inflation versus expectation of inflation by any of the agents.}\label{fig:3SC1G1A1}
\end{figure}

\begin{figure}[h!]
	\centering
	\begin{subfigure}{0.3\textwidth}
		\includegraphics[width=\textwidth]{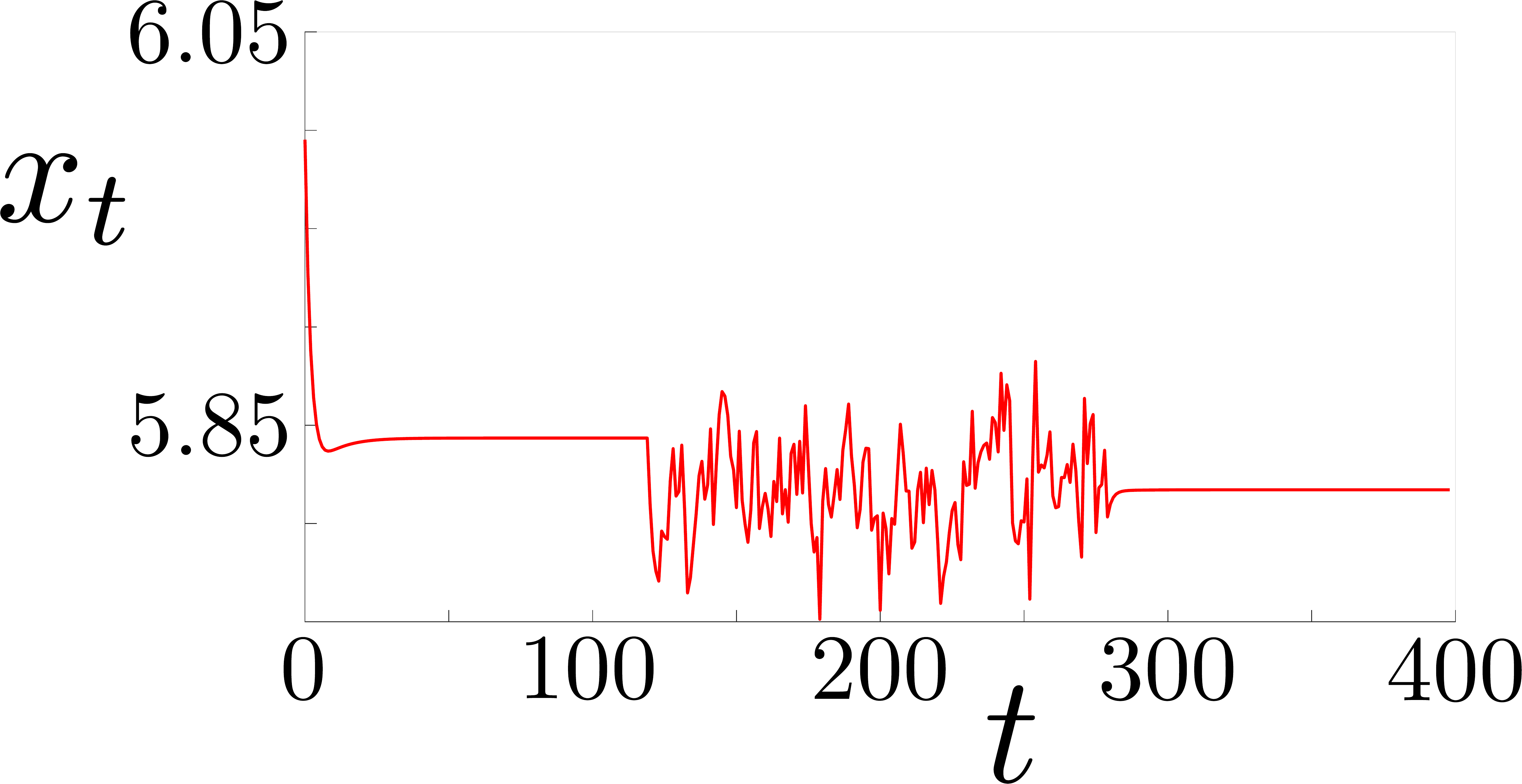}
		\caption{
			}\label{fig:3SC1G1c}
	\end{subfigure}
	\quad
	\begin{subfigure}{0.31\textwidth}
		\includegraphics[width=\textwidth]{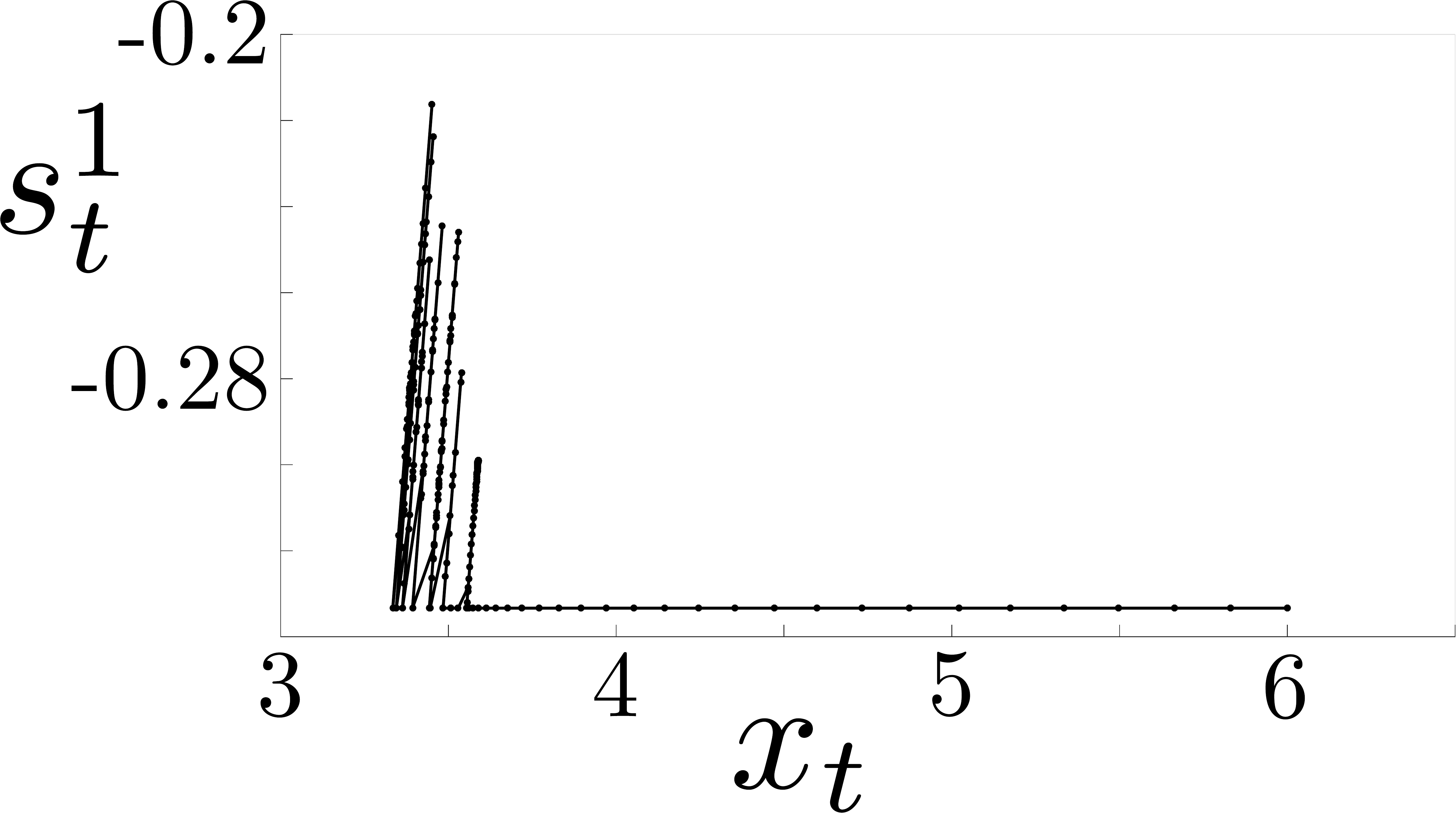}
		\caption{
			}\label{fig:3SC1G1d}
	\end{subfigure}  
	\quad
	\begin{subfigure}{0.31\textwidth}
		\includegraphics[width=\textwidth]{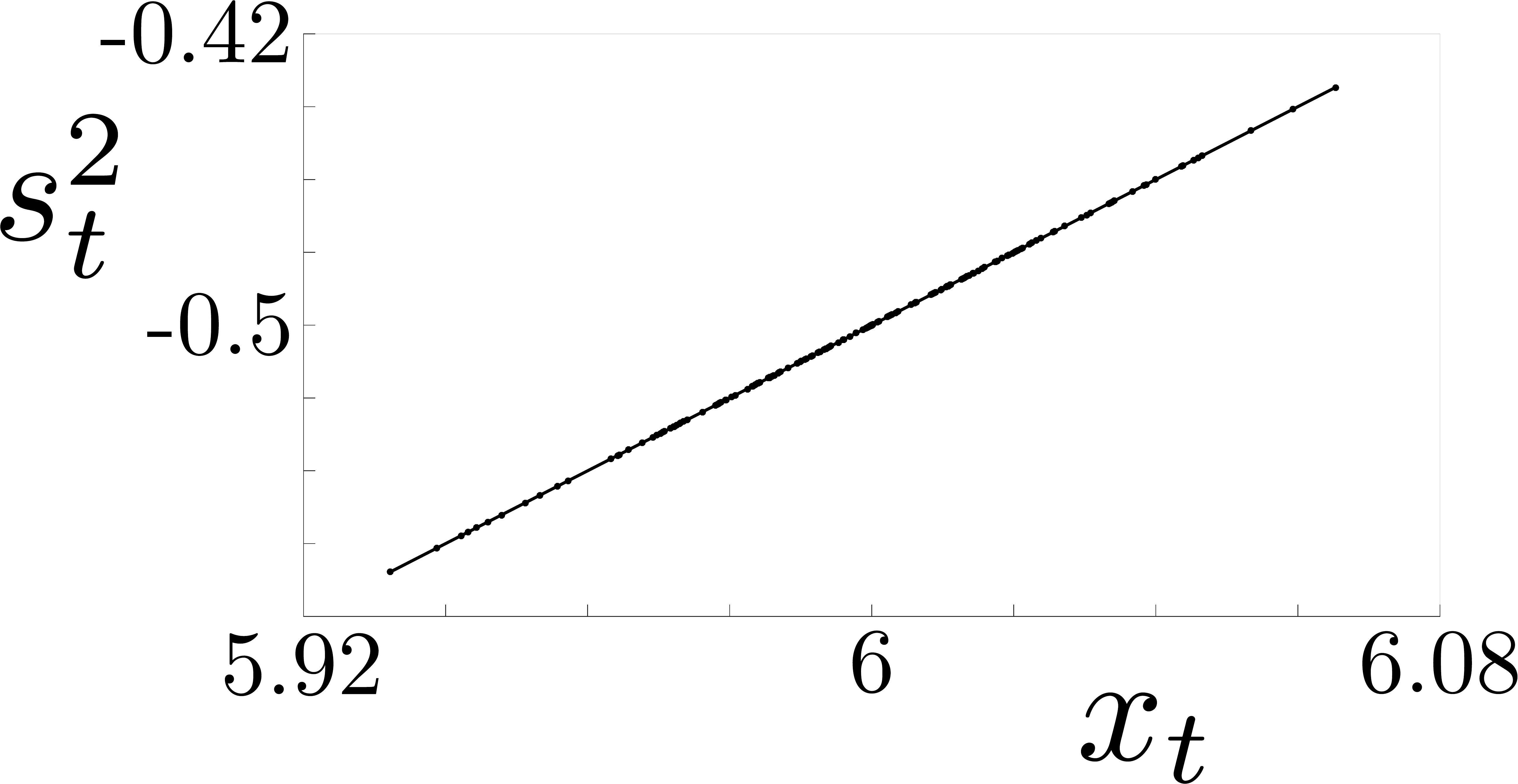}
		\caption{
			}\label{fig:3SC1G1e}
	\end{subfigure}
	\caption{
		Trajectory of the system with 3 agents 
		when 
	the most sensitive agent reaches  an extreme perception gap
        but the two less sensitive agents do not (cf.~Figure \ref{fig:C1G1B}(b, e)). The parameter $c_1$ satisfies $c_1>1$. (a) Time trace of inflation. A change of the equilibrium state occurs. (b) Inflation versus expectation of inflation by the most sensitive agent.
		(c) Inflation versus expectation of inflation by each of the two less sensitive agents.}\label{fig:3SC1G1A2}
\end{figure}

\begin{figure}[h!]
	\centering
	\begin{subfigure}{0.3\textwidth}
		\includegraphics[width=\textwidth]{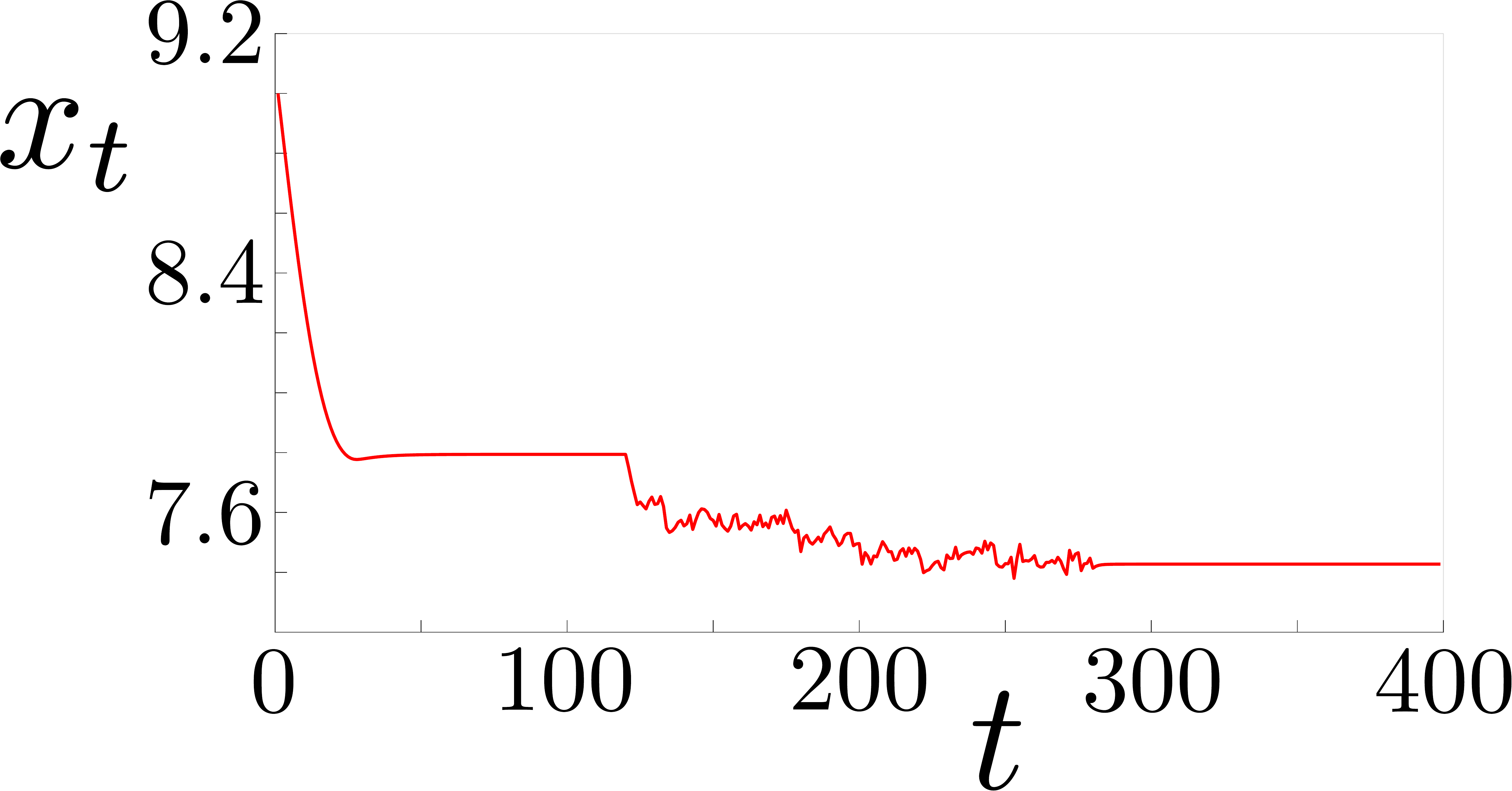}
		\caption{
			}\label{fig:3SC1G1f}
	\end{subfigure}
	\quad
		\begin{subfigure}{0.31\textwidth}
		\includegraphics[width=\textwidth]{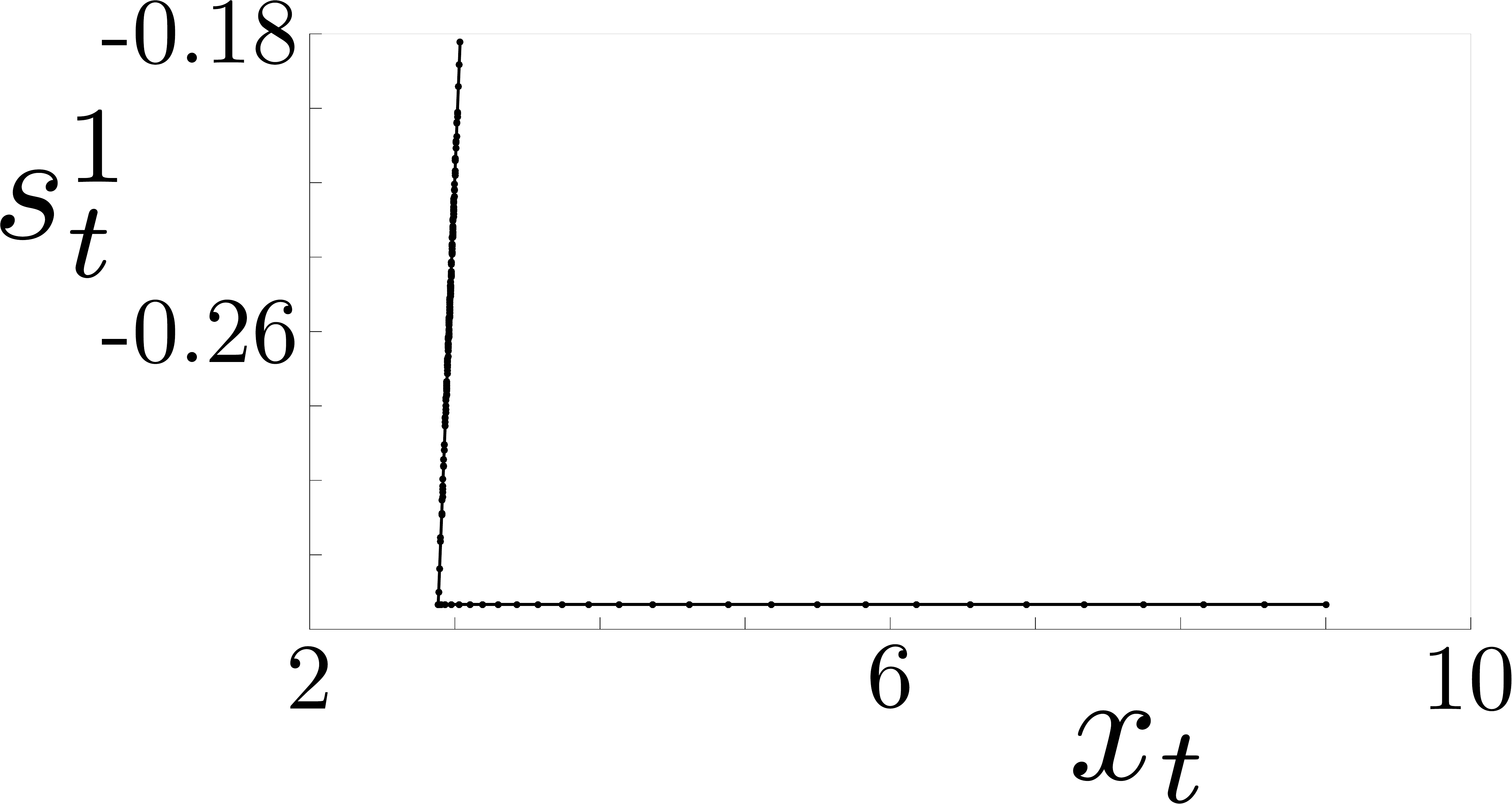}
		\caption{
			}\label{fig:3SC1G1g}
	\end{subfigure}
	\quad
		\begin{subfigure}{0.31\textwidth}
		\includegraphics[width=\textwidth]{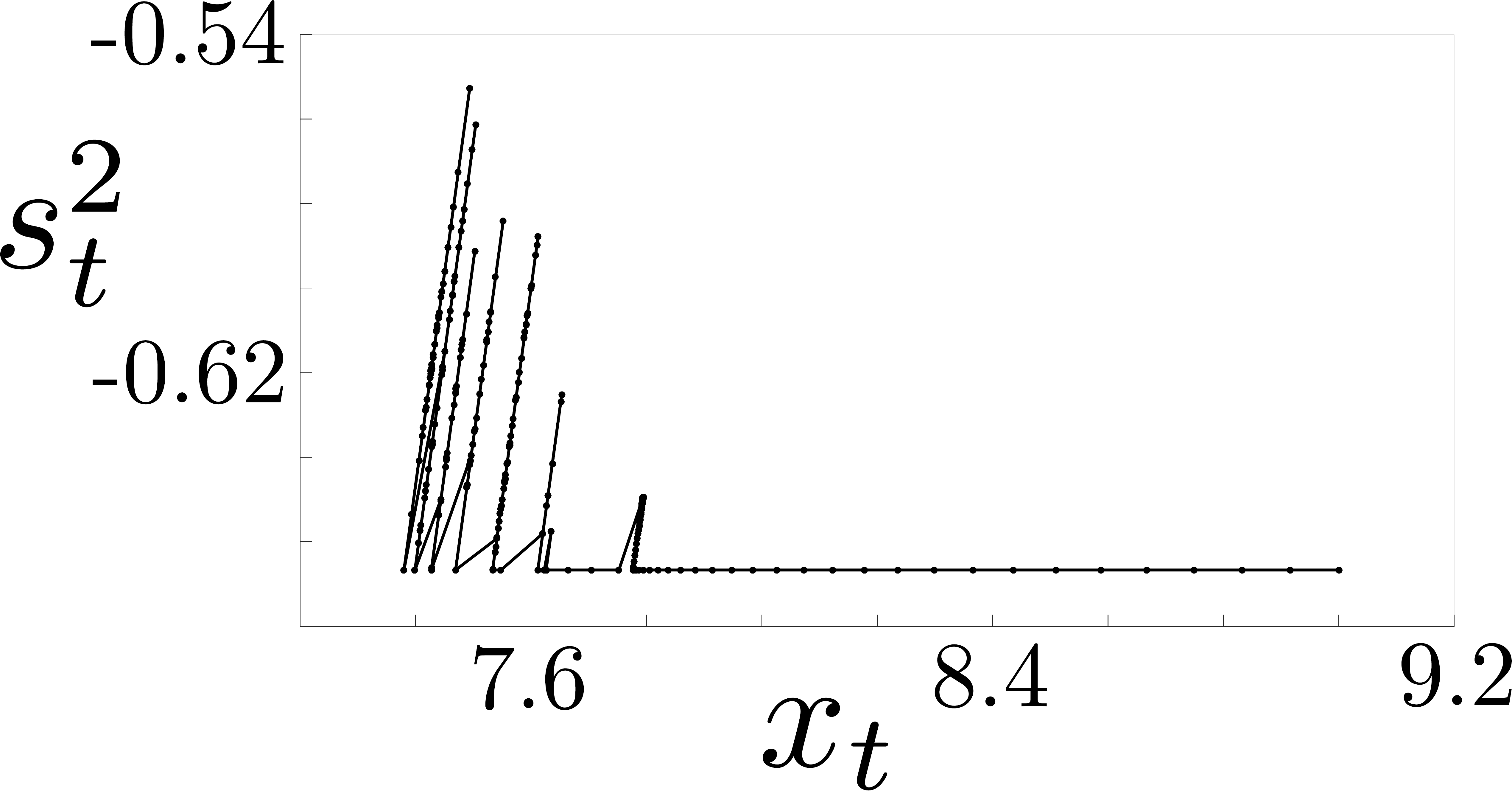}
		\caption{
			}\label{fig:3SC1G1h}
	\end{subfigure}
	\caption{
		Trajectory of the system with 3 agents with the most sensitive agent and the moderately sensitive agent having  an extreme perception gap at the initial (equilibrium) point (cf.~Fig.~\ref{fig:C1G1B}(c, d)). The parameter $c_1$ satisfies $c_1>1$. (a) Time trace of inflation. (b) Inflation versus expectation of inflation for the moderately sensitive agent.
		(c) Inflation versus expectation of inflation for the most sensitive agent.
		The least sensitive agent shows the behavior as in Fig.~\ref{fig:3SC1G1A2}(c).}\label{fig:3SC1G1A3}
\end{figure}

\begin{figure}[h!]
	\centering
	\begin{subfigure}{0.48\textwidth}
		\includegraphics[width=\textwidth]{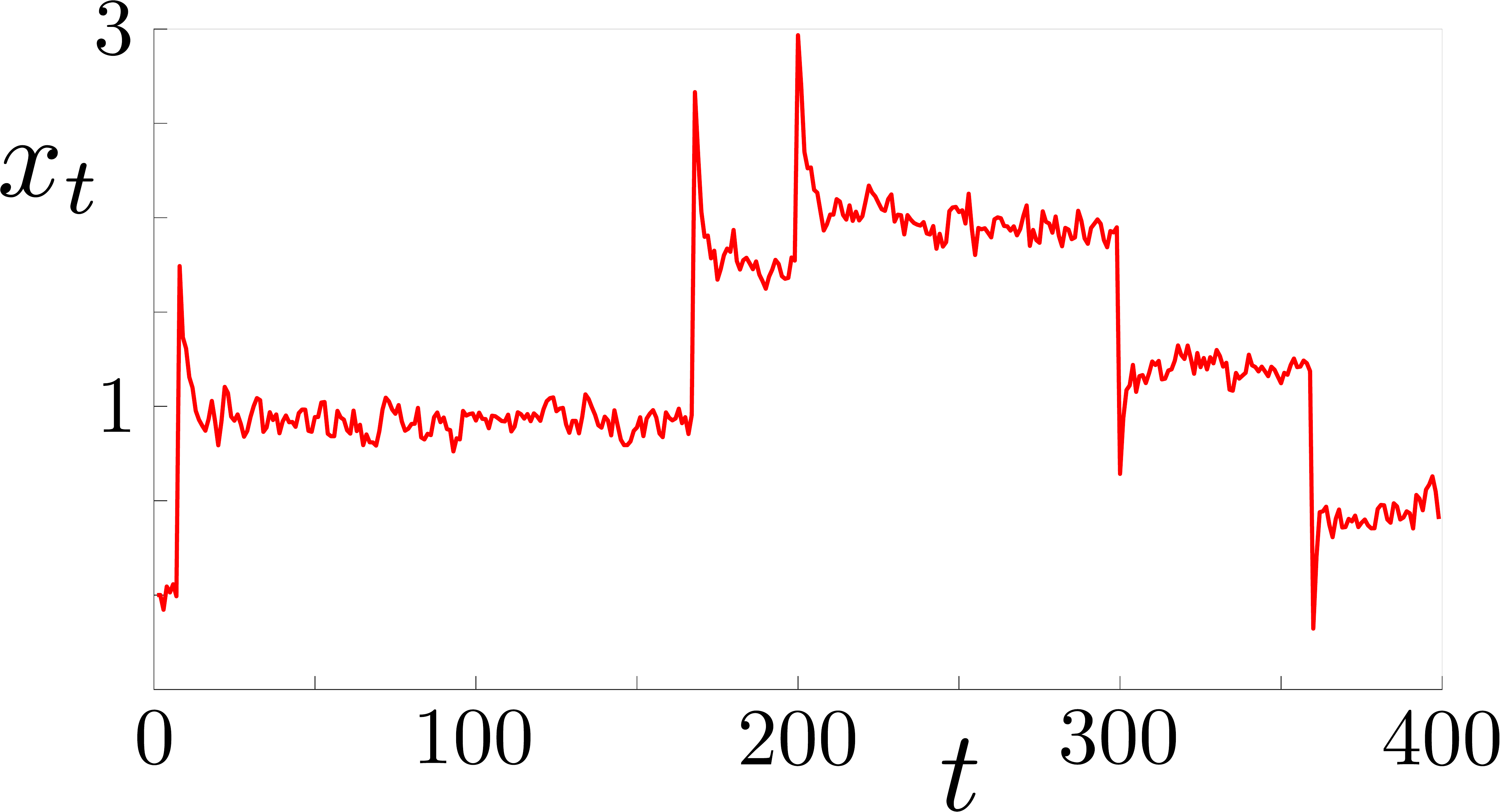}
		\caption{}
		\label{fig:3S2DRS1sx}
	\end{subfigure}
	\quad	
	\begin{subfigure}{0.48\textwidth}
		\includegraphics[width=\textwidth]{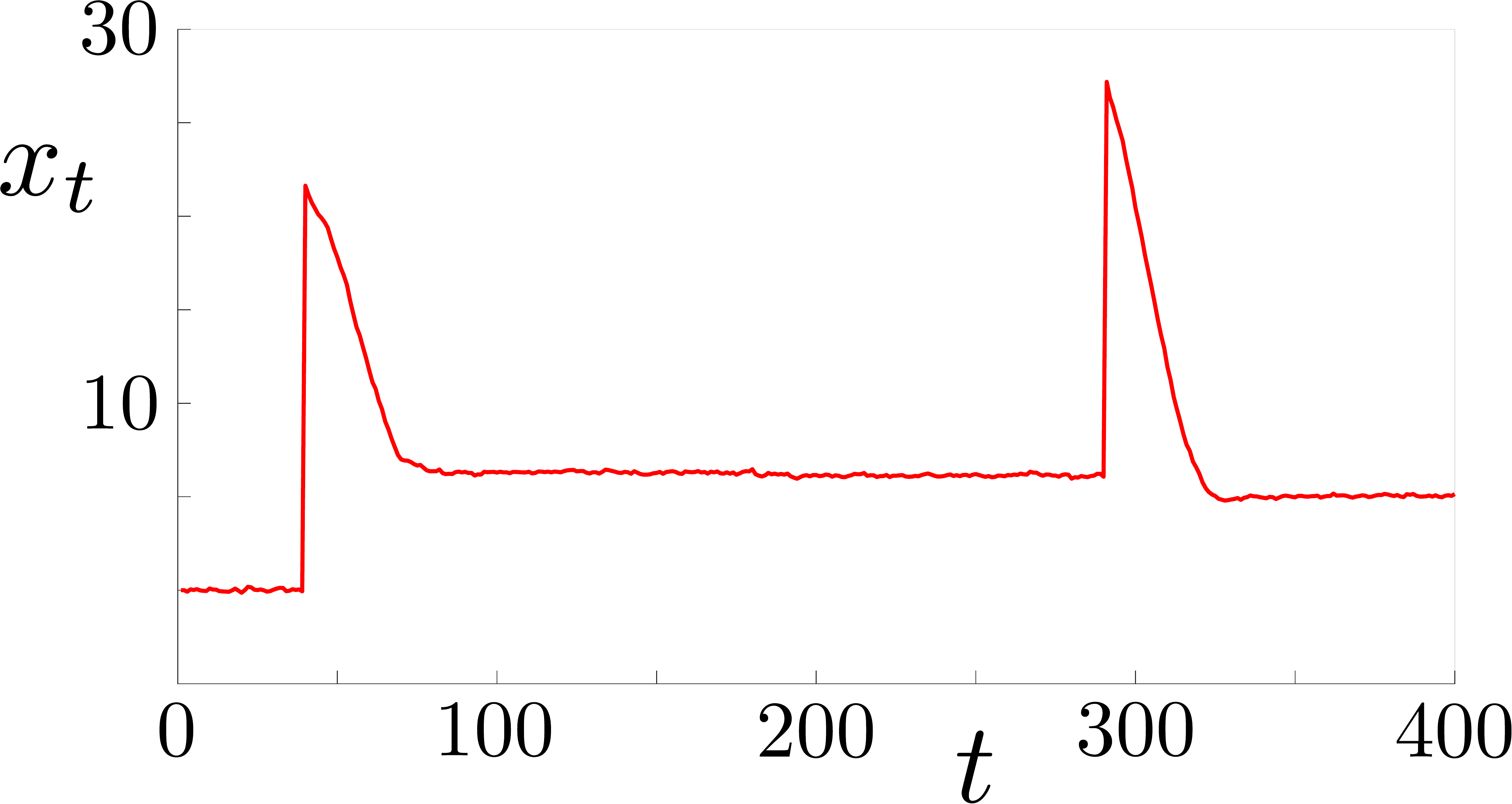}
		\caption{}
		\label{fig:3S2DRS1lx}
	\end{subfigure}
	\caption{Changes of the equilibrium state in the model with 3 agents due to shocks (cf.~Figures \ref{fig:CompareShocks}, \ref{fig:2DRS1}).
		(a) Small shocks. (b) Relatively large shocks.
	}\label{fig:3S2DRS2}
\end{figure}

\begin{figure}[h!]
	\centering
	\includegraphics[scale=0.15]{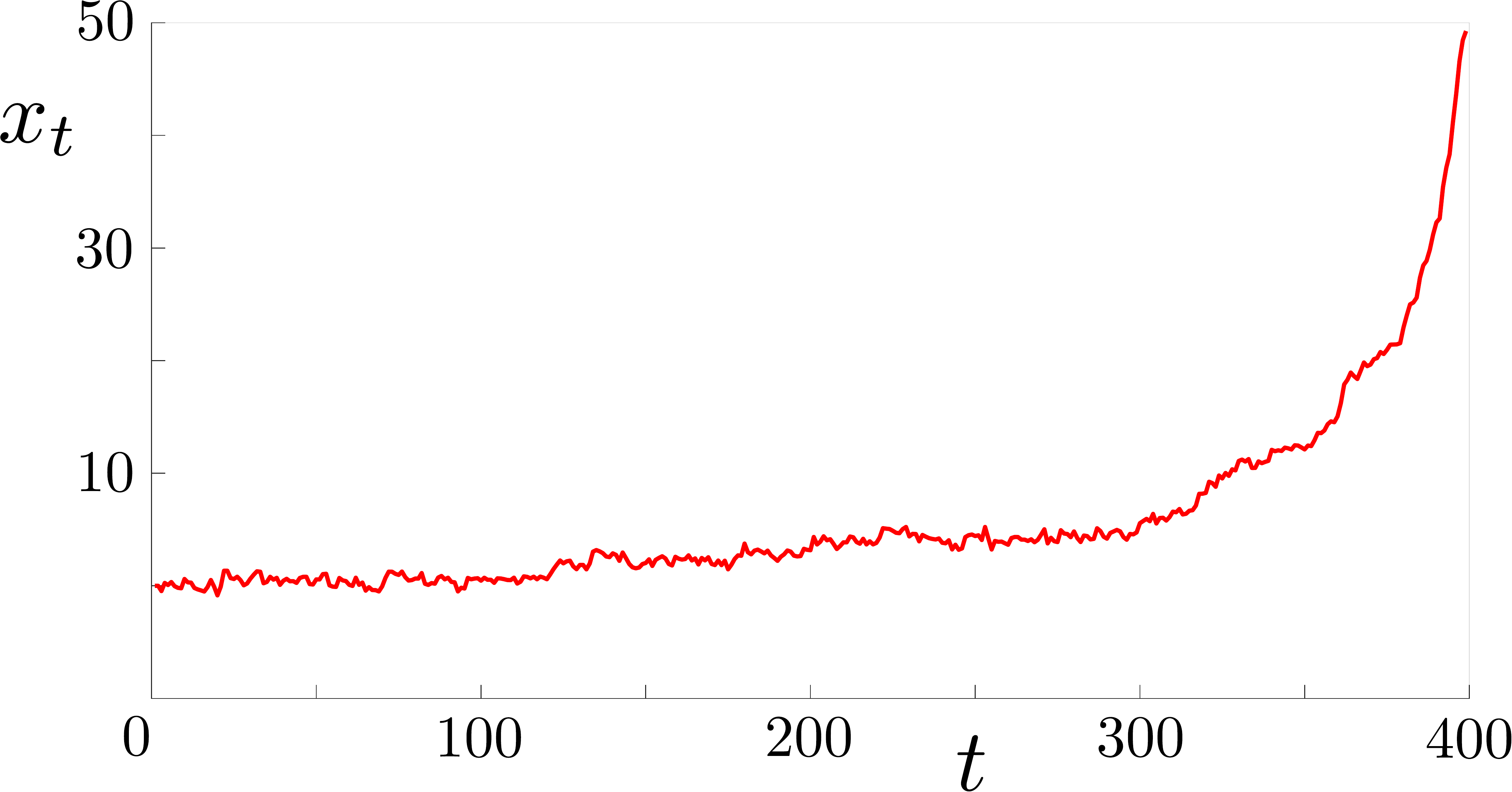}
	\caption{The run-away inflation scenario in the model with 3 agents in the case $c_1<1$ (cf.~Fig.~\ref{fig:C1L1B}).} \label{fig:3SRunAX1X}
\end{figure}

\begin{figure}[h!]
	\centering
	\begin{subfigure}{0.48\textwidth}
		\includegraphics[width=\textwidth]{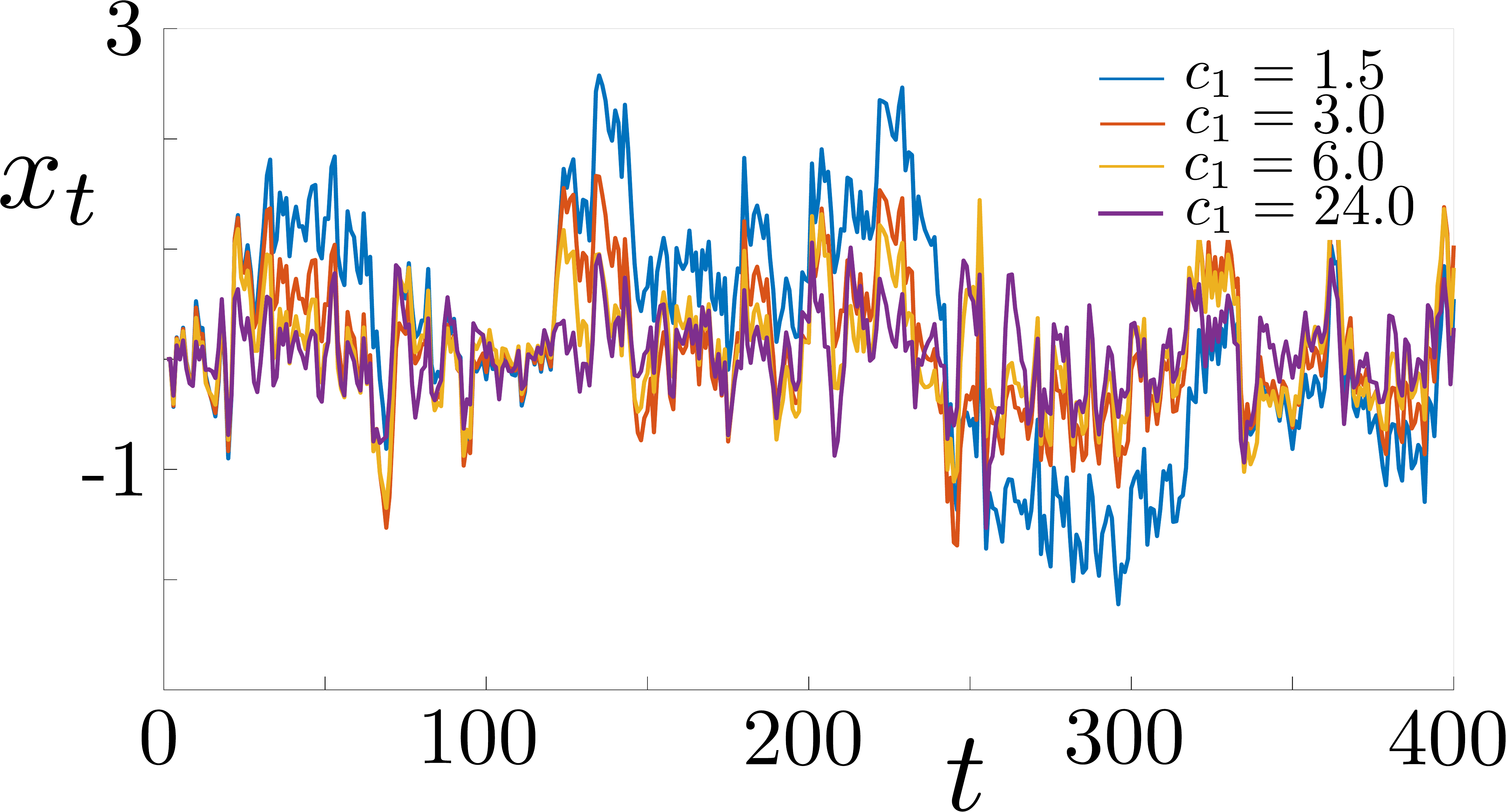}
		\caption{ }
		\label{fig:3SC1Inc1a}
	\end{subfigure}
	\quad
	\begin{subfigure}{0.48\textwidth}
		\includegraphics[width=\textwidth]{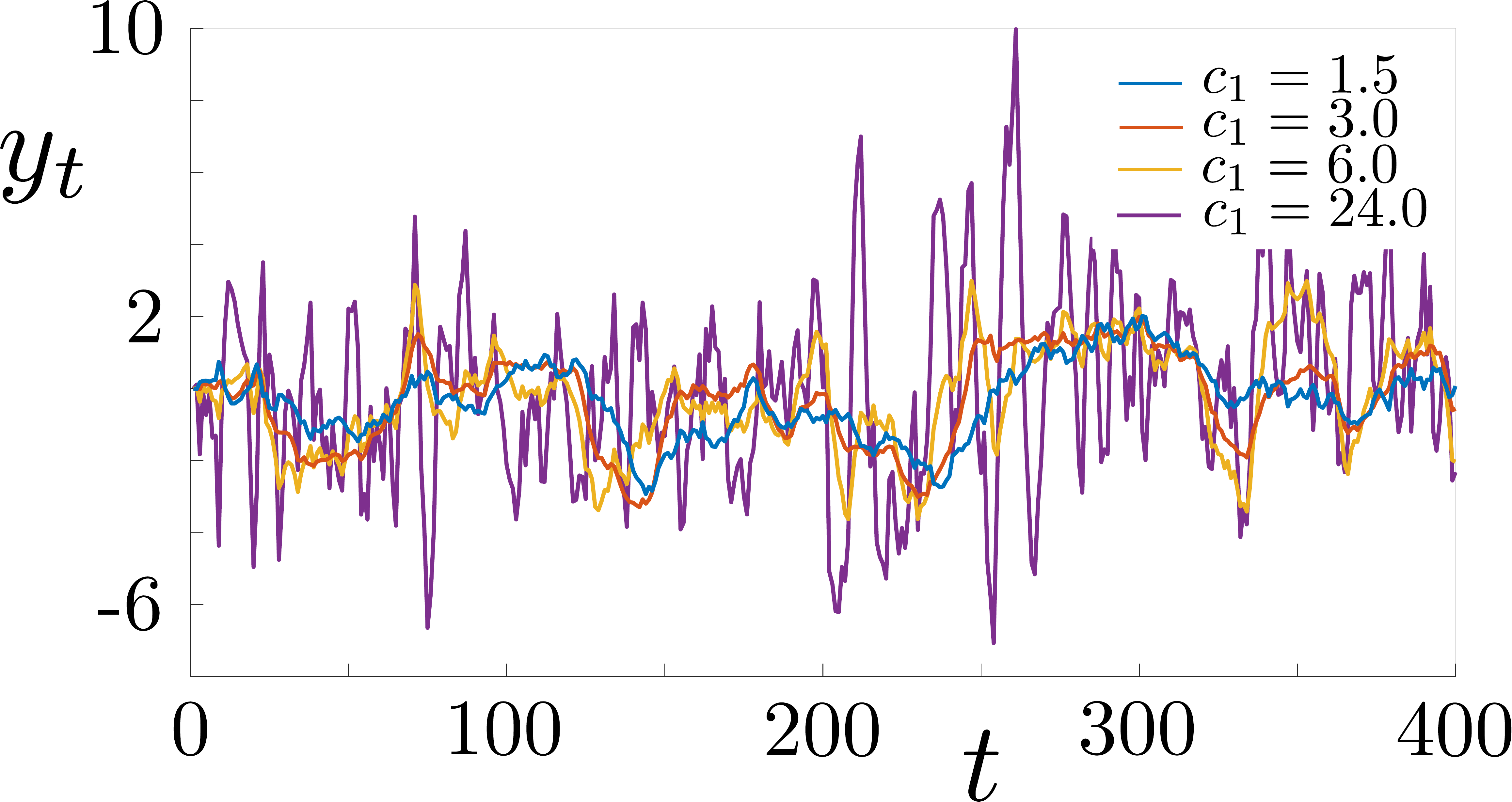}
		\caption{}
		\label{fig:3SC1Inc1b}
	\end{subfigure}        
	\caption{
	Trade-off between the inflation and output gap volatility in the model with 3 agents as the inflation targeting parameter $c_1$ in the Taylor rule is varied (cf.~Fig.~\ref{fig:C1Inc3}). (a) Trajectories of $x_t$. (b) Trajectories of $y_t$.   }\label{fig:3SC1Inc1}
\end{figure}

\begin{figure}[h!]
	\centering
	\begin{subfigure}{0.48\textwidth}
		\includegraphics[width=\textwidth]{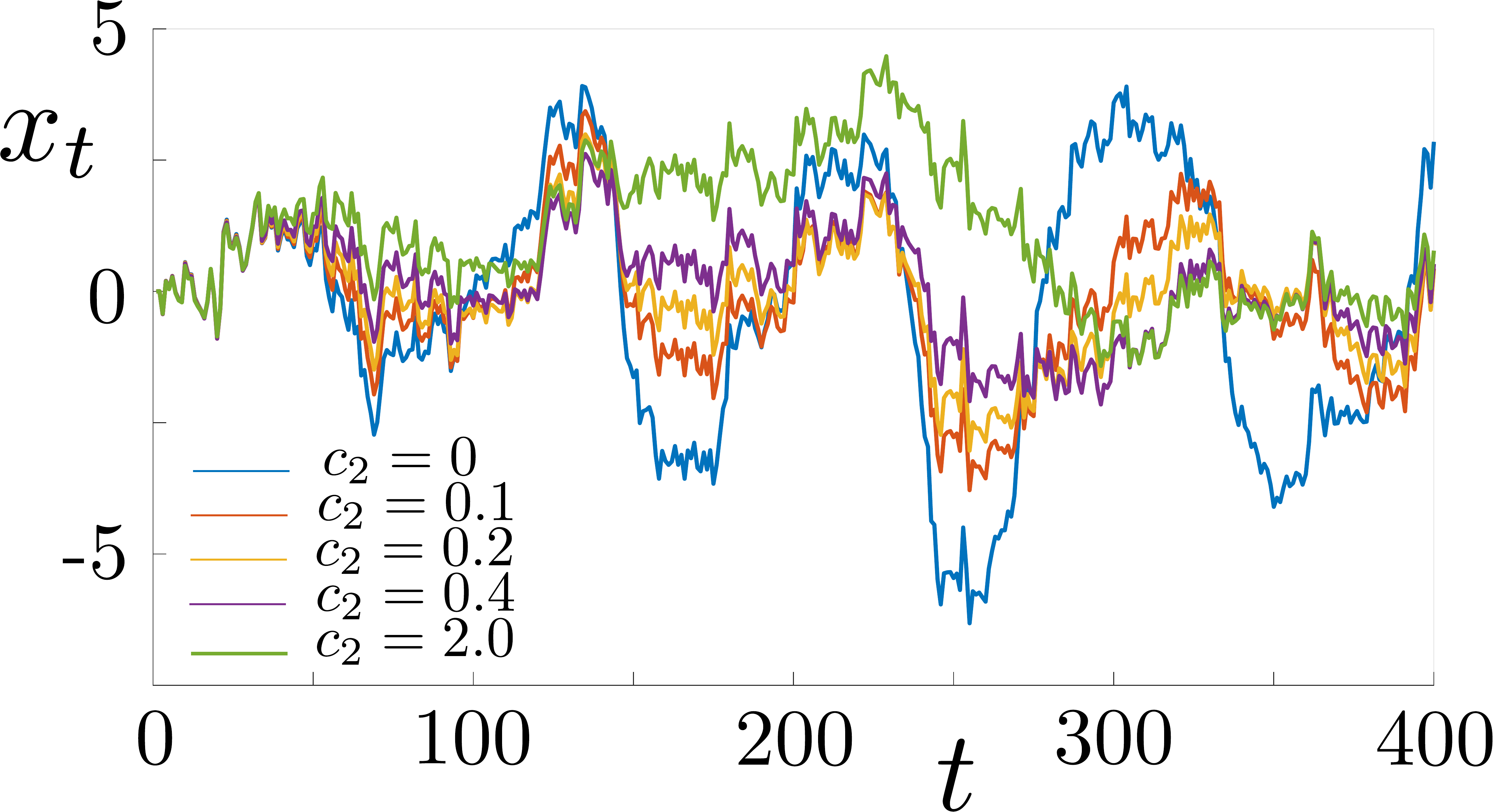}
		\caption{ }
		\label{fig:3SC2Inc1a}
	\end{subfigure}
	\quad
	\begin{subfigure}{0.48\textwidth}
		\includegraphics[width=\textwidth]{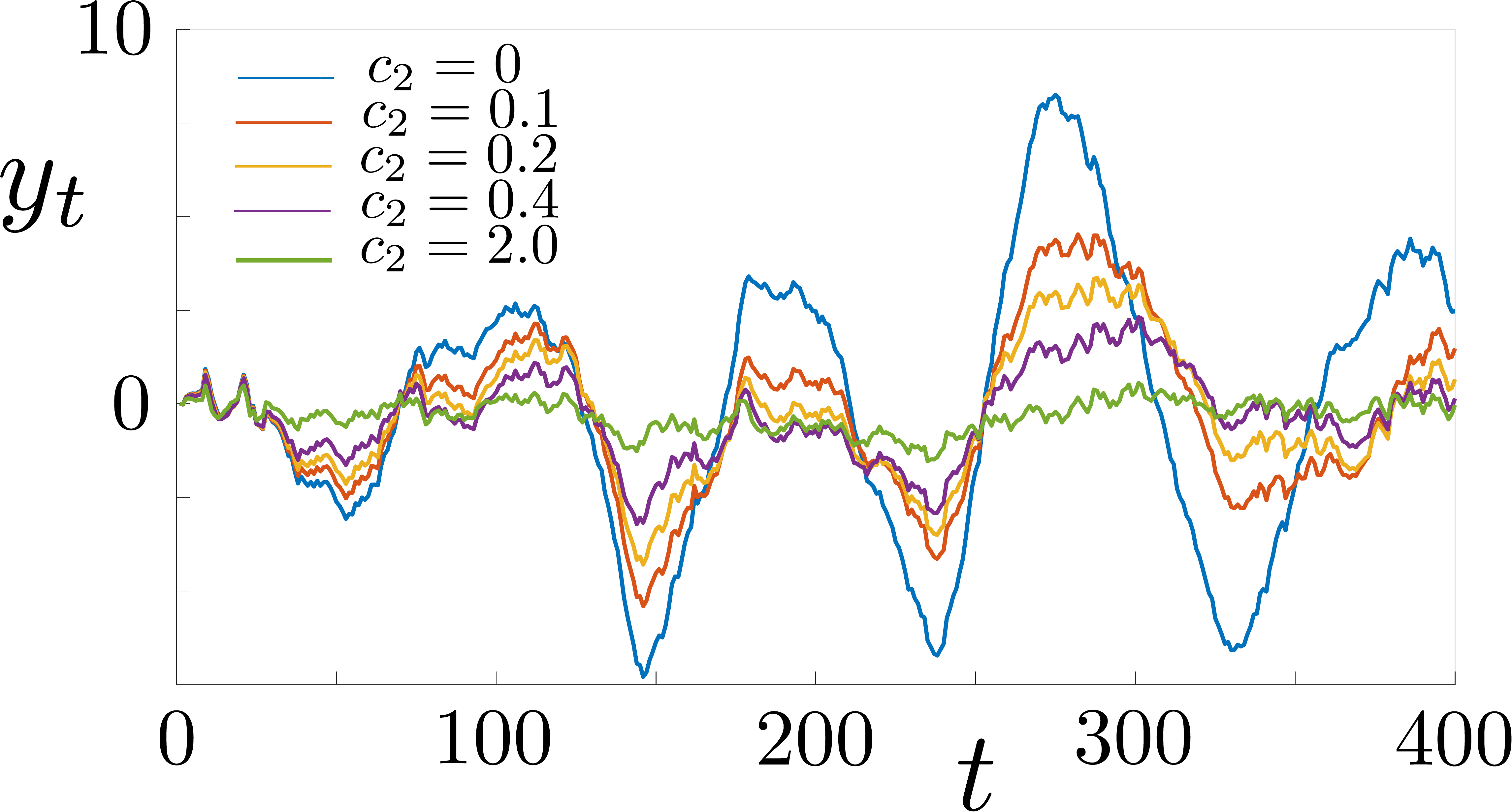}
		\caption{}
		\label{fig:3SC2Inc1b}
	\end{subfigure}           
	\caption{
		Trade-off between the inflation rate and  output gap volatility in the model with 3 agents as the output gap targeting parameter $c_2$ in the Taylor rule is varied (cf.~Fig.~\ref{fig:C2Inc3A}).   (a) Trajectories of $x_t$. (b) Trajectories of $y_t$.  }\label{fig:3SC2Inc1}
\end{figure}

\begin{figure}[h!]
	\centering
	\begin{subfigure}{0.48\textwidth}
		\includegraphics[width=\textwidth]{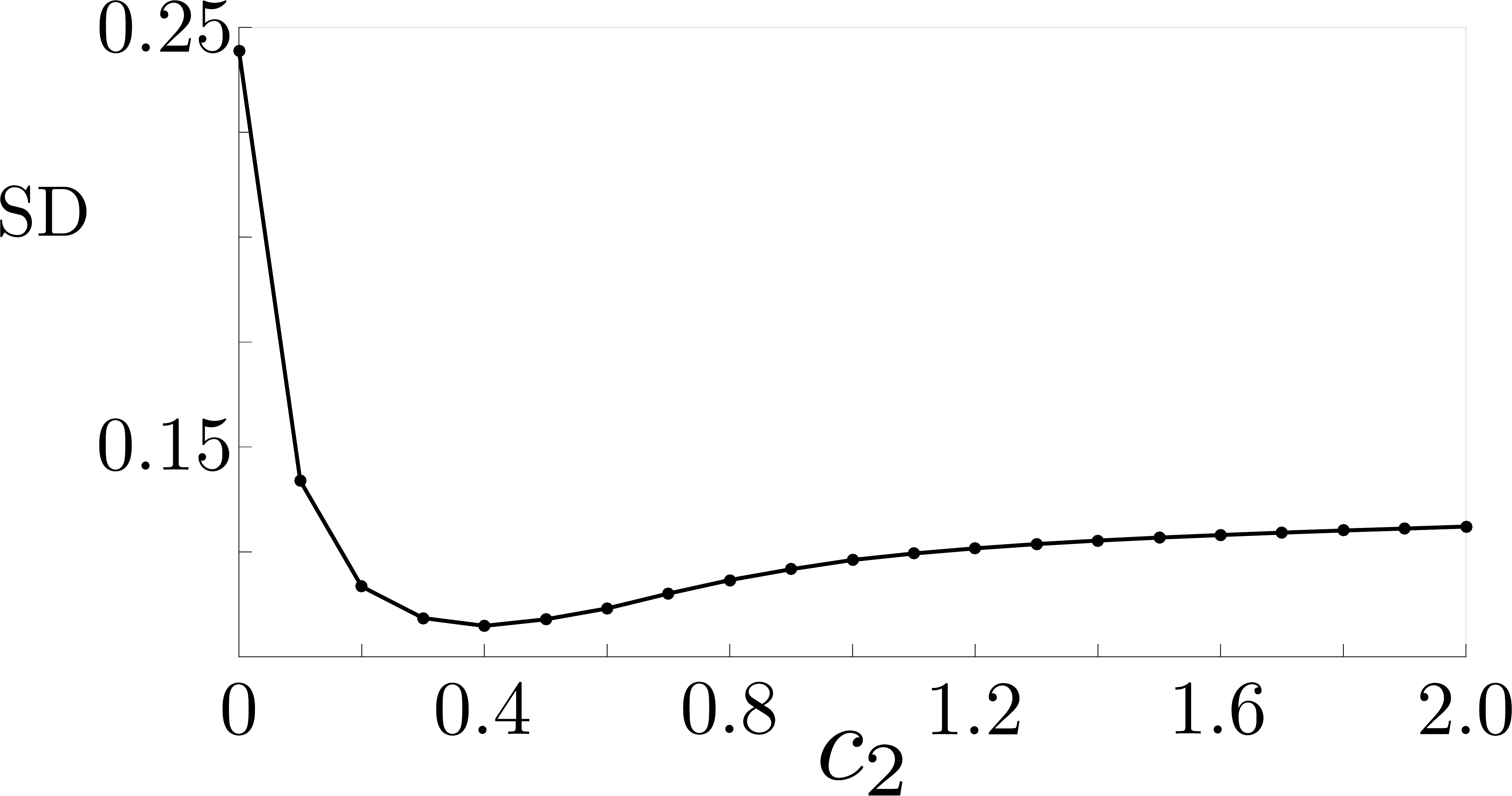}
		\caption{ }
		\label{fig:3SC2Inc1c}
	\end{subfigure}
	\quad
	\begin{subfigure}{0.48\textwidth}
		\includegraphics[width=\textwidth]{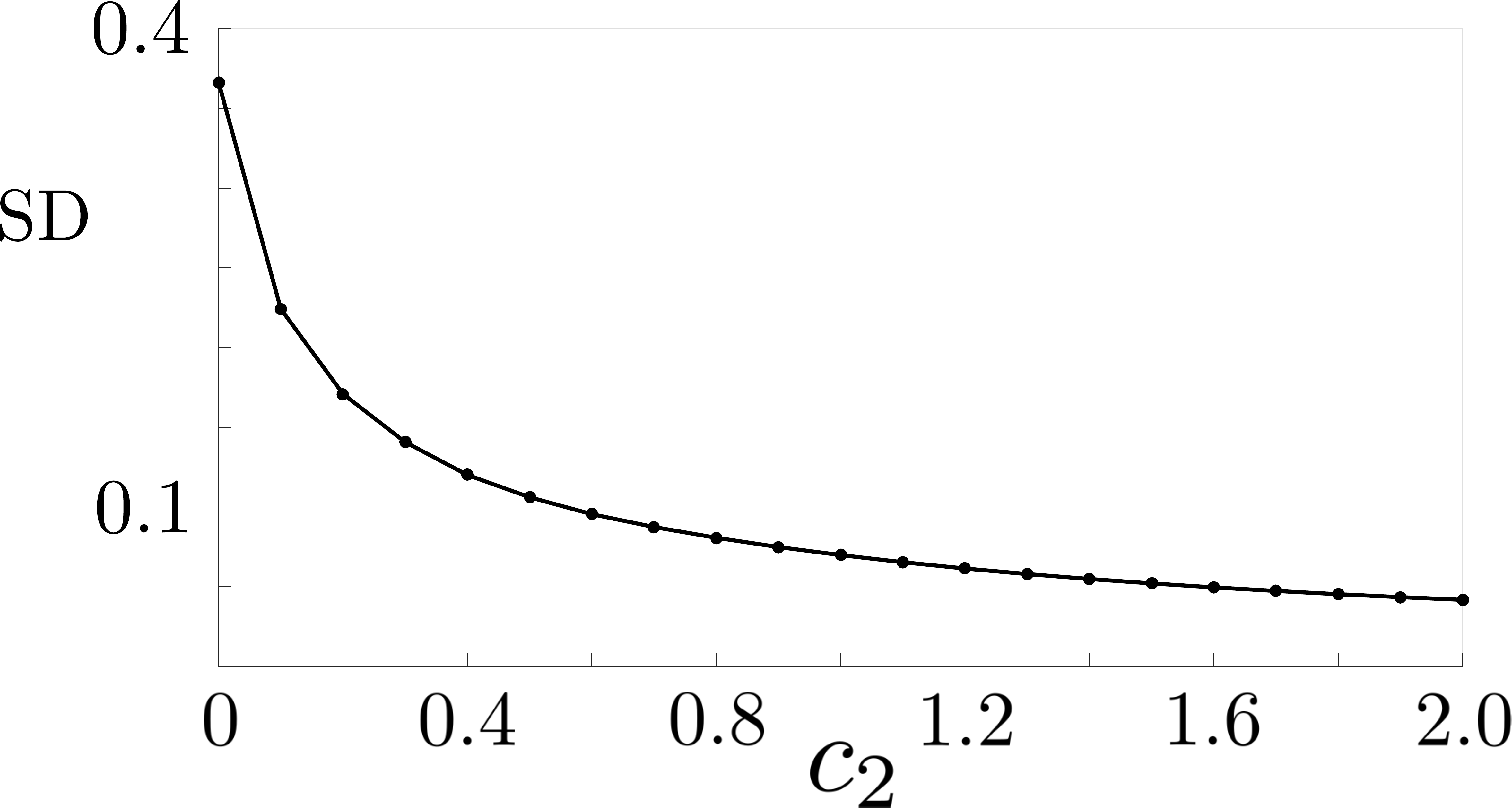}
		\caption{}
		\label{fig:3SC2Inc1d}
	\end{subfigure}           
	\caption{Measure of the effect of $c_2$ on volatility of (a) inflation rate, $x_t$ and (b) output gap, $y_t$   with  standard deviation (SD) 
		(cf.~Fig.~\ref{fig:C2STDInc3A}).}
	\label{fig:3SC2Inc2}
\end{figure}

Overall, numerical results obtained for model \eqref{eqn:M1}, \eqref{eqn:M1'}, \eqref{eqn:MStop}
with three agents are qualitatively similar to the results described above for the model with one  agent, see Figs.~ \ref{fig:3SC1G1A1} -- \ref{fig:3SC2Inc2}, which are counterparts of Figs.~\ref{fig:C1G1B} -- \ref{fig:C2STDInc3A}, respectively.


\subsection{A sticky Central Bank model}\label{interest}

The Central Bank policy can presumably exhibit stickiness too.
To explore this scenario in this Section we shall replace the Taylor rule
\eqref{eqn:M1'} with the relation
\begin{equation}\label{M1'}
r_t=\mathcal{P}_\sigma[c_1x_t+c_2y_t]
\end{equation}
also involving a play operator. But at the same time, for the sake of
simplicity and in order to isolate the effect of stickiness in the
Central Bank response upon the system, we remove the play operator from
equations \eqref{eqn:M1} thus assuming that the aggregate expectation
of inflation equals to the current actual inflation rate, $p_t=x_t$;
this corresponds to setting $\rho=0$ in equations \eqref{eqn:M1}. In
this case,
\begin{equation}\label{M1}
\begin{array}{l}
y_t=y_{t-1}-a(r_t-x_t)+\epsilon_t,\\
x_t=x_{t-1}+\frac{b_2}{1-b_1}y_t+\eta_t.\\
\end{array}
\end{equation}
It would be interesting to consider the model with both sticky inflation expectation and sticky Central Bank response,
however this is beyond the scope of this paper.

System \eqref{M1'}, \eqref{M1} can be written in the form \eqref{eqn:SA1} with 
\[
s_t=\mathcal{S}_\sigma[c_1x_t+c_2y_t],
\]
the matrix $A$ defined by \eqref{AA}, $N=A$, and $d=(a(1-b_1), ab_2)^\top/\Delta$ with $\Delta$ defined by \eqref{Delta}. 
The technique presented in Subsection \ref{inve} can be adapted to convert the 
implicit system \eqref{M1'}, \eqref{M1} into a well-defined explicit system 
provided that 
\begin{equation}
\label{condiinv}
1-b_1-ab_2>0.
\end{equation}
(see Appendix E). Hence, we assume that this condition is satisfied.

Equilibrium states of system \eqref{M1'}, \eqref{M1} with zero noise terms form the line segment
\begin{equation}\label{seg}
(y_*(s_*),x_*(s_*))=\Bigl(0, \frac{s_*}{c_1-1}\Bigr),\qquad s_*\in[-\sigma,\sigma].
\end{equation}
Notice that the output gap value is zero for all the equilibrium states, while the equilibrium inflation rate
ranges over an interval of values. Notably, the local stability analysis (see Appendix E) shows that 
all the equilibrium states with $s_*\in(-\sigma,\sigma)$ are {\em unstable} for any set of parameter values.
That is, stickiness in the Taylor rule leads to destabilization of equilibrium states.

\begin{figure}[ht]
	\begin{center}
		\begin{subfigure}{0.32\textwidth}
			\includegraphics[scale=0.08]{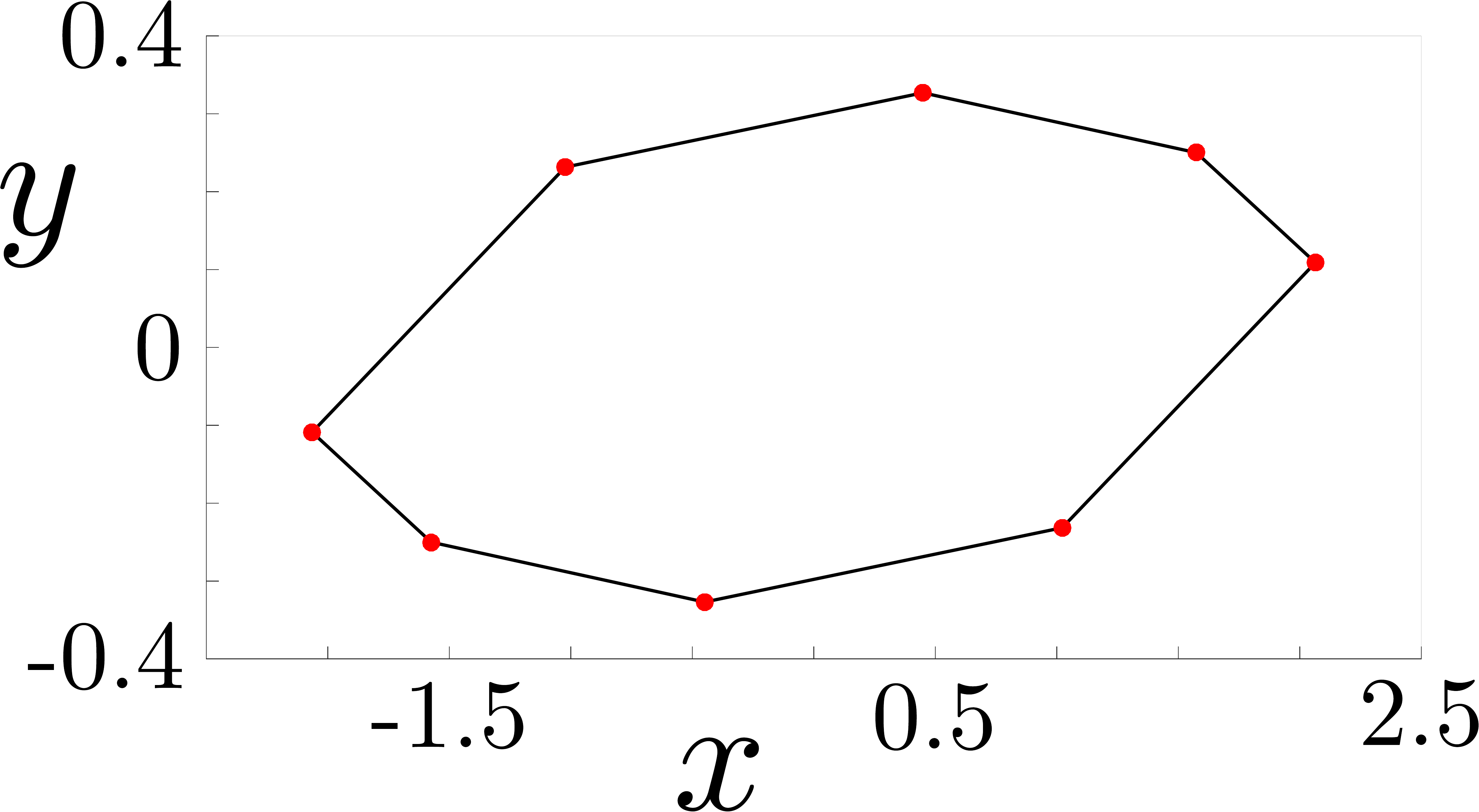}
			\caption{}
		\end{subfigure}\quad
		\begin{subfigure}{0.31\textwidth}
			\includegraphics[scale=0.08]{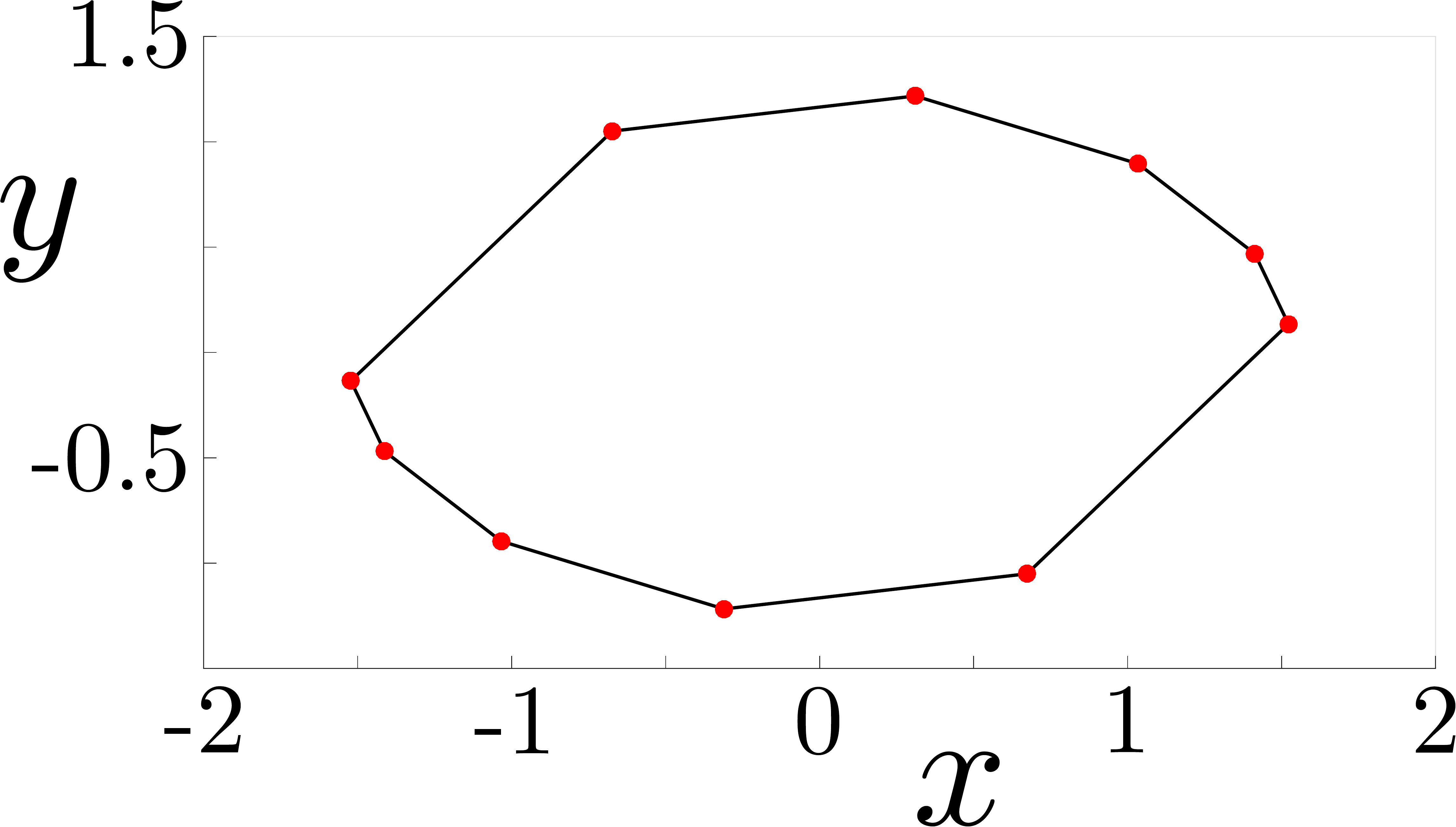}
			\caption{}
		\end{subfigure}\quad
		\begin{subfigure}{0.31\textwidth}
			\includegraphics[scale=0.08]{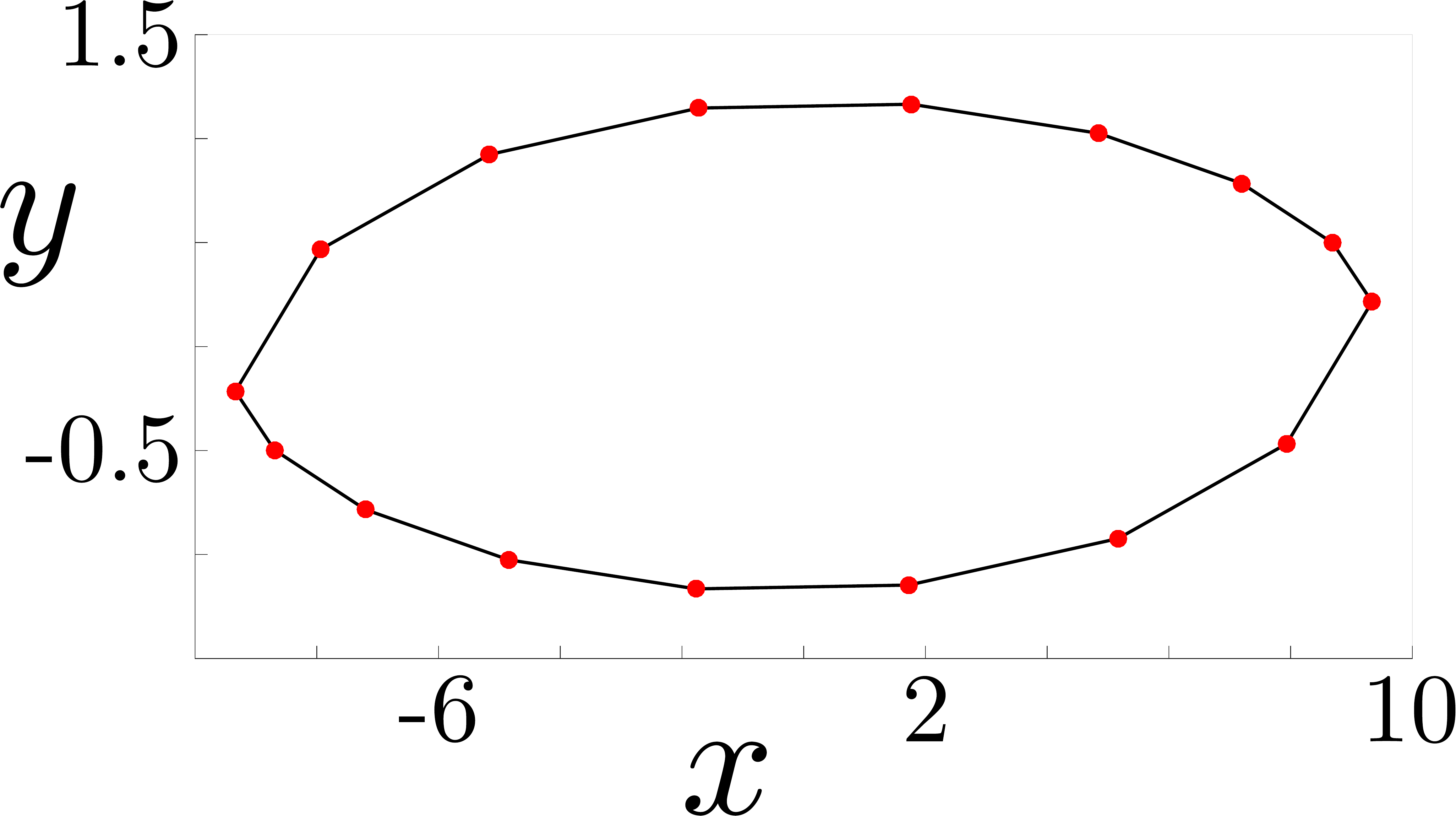}
			\caption{}
		\end{subfigure}
   	\end{center}
   
		
	\begin{center}
		\begin{subfigure}{0.32\textwidth}
			\includegraphics[scale=0.08]{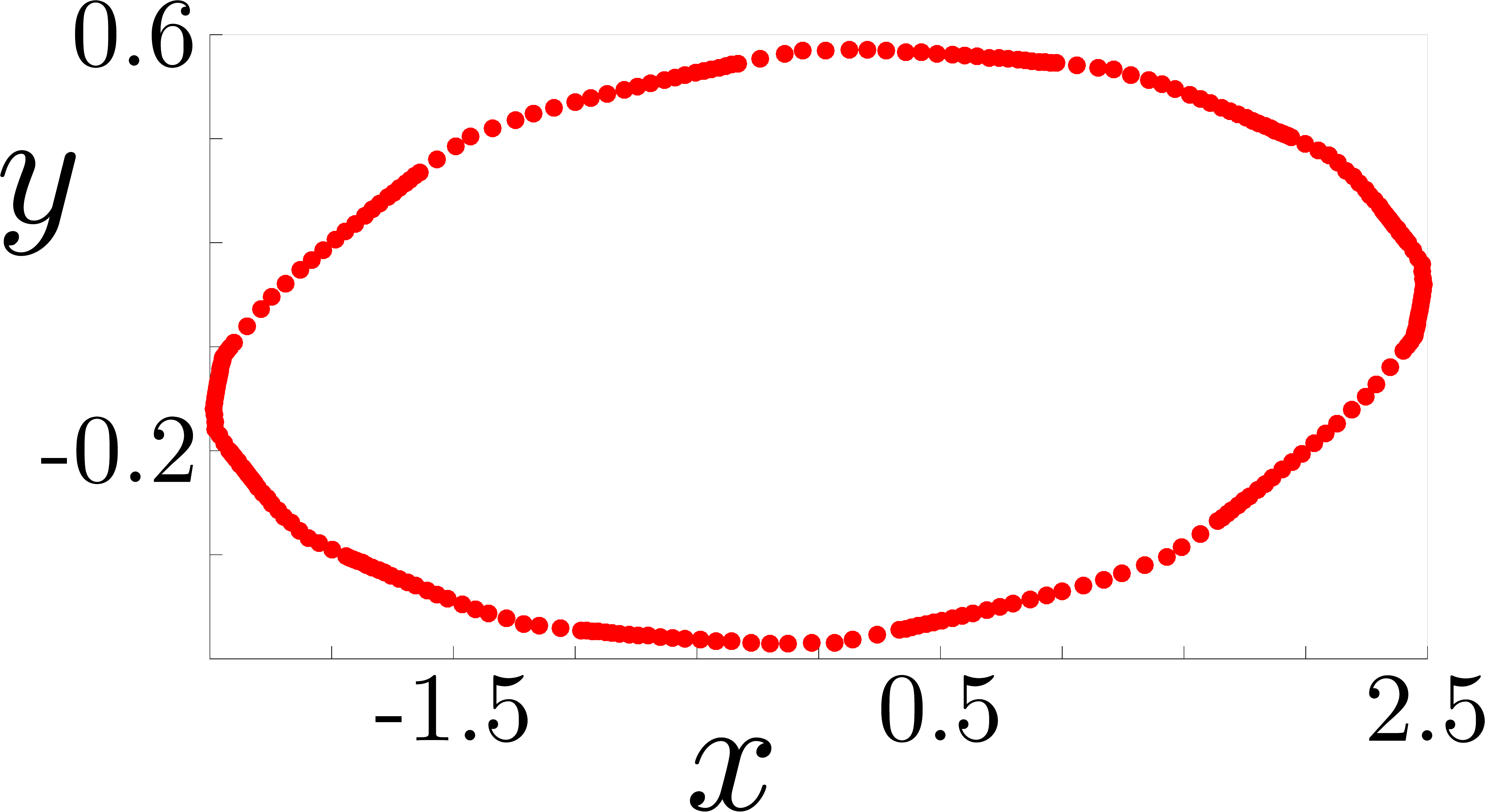}
			\caption{}
		\end{subfigure}\quad
		\begin{subfigure}{0.31\textwidth}
			\includegraphics[scale=0.08]{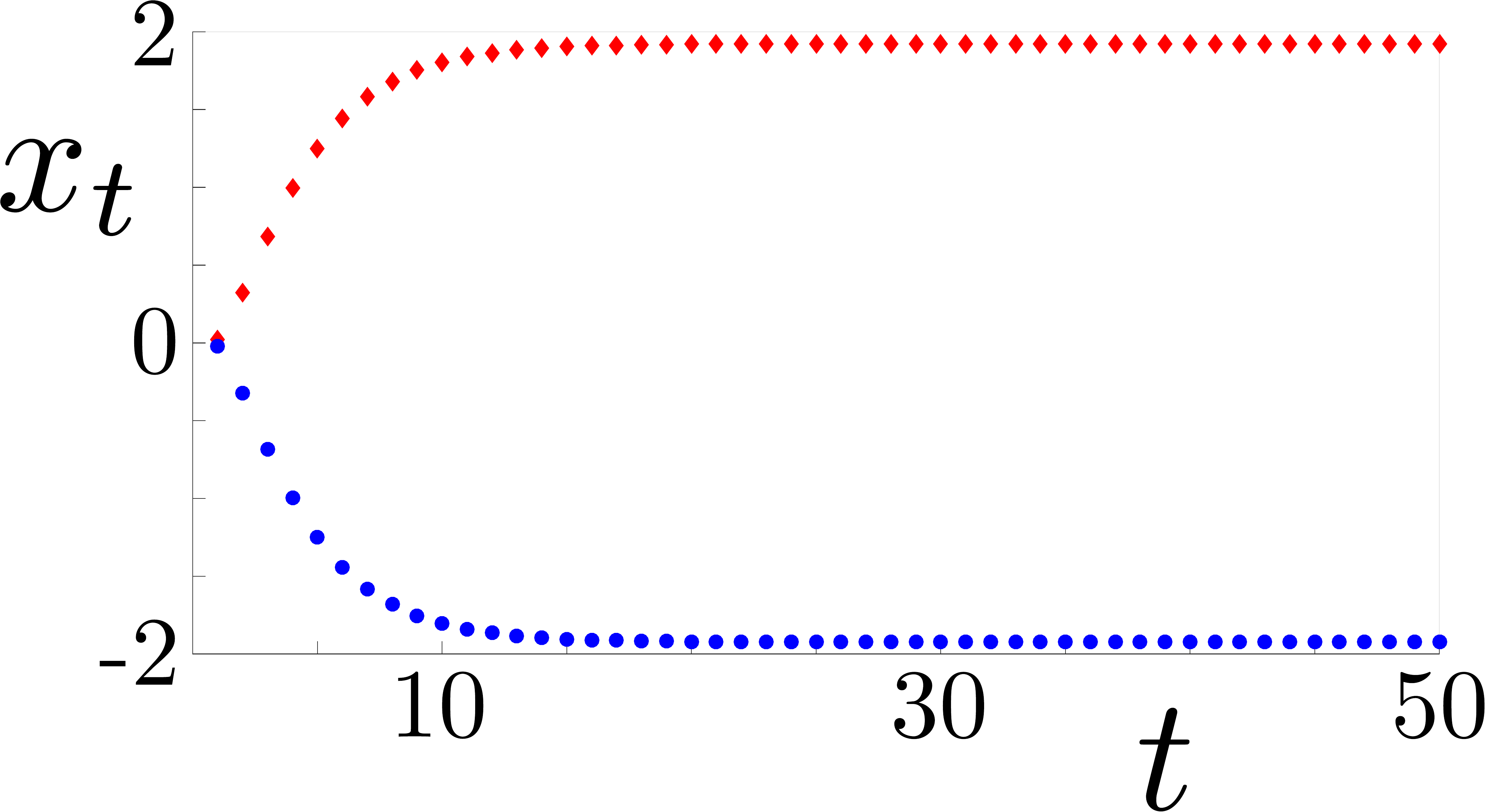}
			\caption{}
		\end{subfigure}\quad
		\begin{subfigure}{0.31\textwidth}
			\includegraphics[scale=0.08]{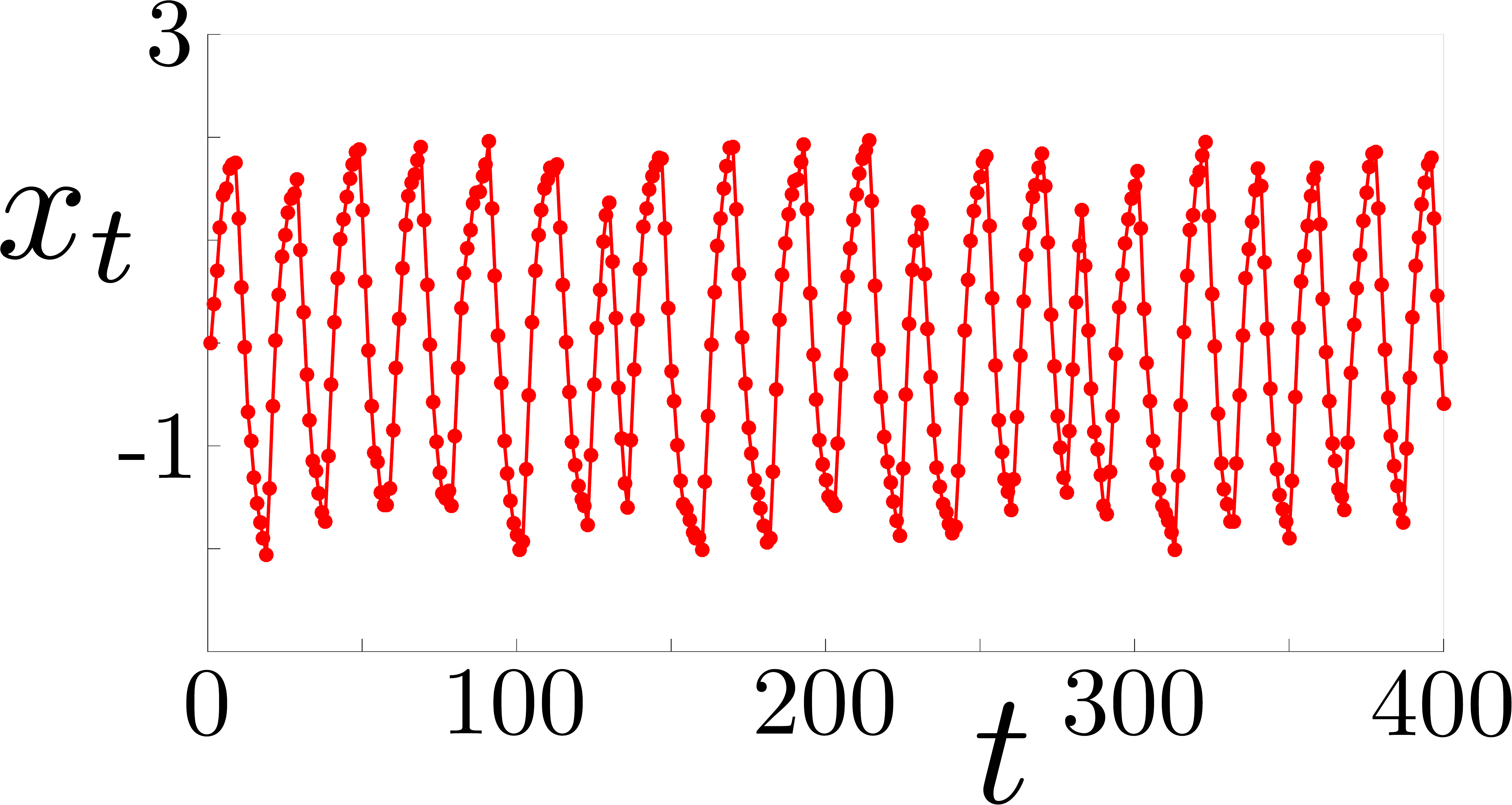}
			\caption{}
		\end{subfigure}
		\vspace{-1mm}
		\caption{An attractor of system \eqref{M1'}, \eqref{M1} for several parameter sets.		
			(a -- c) A periodic orbit (with period 8, 10, 16, respectively)  shown on the $(x,y)$ plane
			for the system without noise. (d) A quasiperiodic orbit.		
			(e) Two equilibrium states corresponding to $s_*=\pm\sigma$ (the time trace
			of $x_t$ shown for 2 trajectories). (f) Time trace of $x_t$ for a trajectory of the system with noise  for the same parameters as in (e).
		} \label{fifi}
	\end{center}
\end{figure}

On the other hand, for large values of $z_t=(y_t,x_t)^\top$, the system can be approximated by equation \eqref{eqn:SA1'''}, which is exponentially stable (as shown in Appendix B).
This ensures that in the system \eqref{M1'}, \eqref{M1},
in the absence of noise, all trajectories converge to a bounded domain $\Omega$
surrounding the segment of equilibrium states and, upon entering this
domain, remain there. However, since the equilibria are all unstable,
more complicated bounded attracting orbits (such as periodic,
quasiperiodic, or even chaotic atractors) must occur.  
Fig.~\ref{fifi} shows a few possibilities for the attractor of system \eqref{M1'}, \eqref{M1} 
obtained for different sets of parameter values. The attractor belongs
to $\Omega$ whose size is controlled by the parameter $\sigma$ of the sticky Taylor rule \eqref{M1'}.
This size can be estimated using the Lyapunov function introduced in Subsection \ref{staglo}. 

Finally, we note that in the presence of noise, a trajectory will most
likely wander unpredictable around $\Omega$ unless kicked outside
temporarily by a
fluctuation.

\section{Conclusions}

In this paper we rigorously analyzed a simple macroeconomic model
with sticky inflation expectations. Perhaps surprisingly, although the
model is nonlinear it can be considered as a hybrid (piecewise)
linear system and analysed using a mostly linear mathematical toolkit.

For such a simple model, defined via a single (and conceptually quite
elementary) change from a standard one, the sticky play operator
introduces surprisingly complicated, subtle-yet-recognizable phenomena into
the dynamics.  Some of the more detailed conclusions of our
simulations may be model-specific but, based upon the mathematics
presented here and additional numerical simulations with more complex
variants of the model, we believe at least the following two
qualitative features to be generic and robust.

Firstly, the presence of an entire continuum of equilibria rather
than a unique one (or even finite numbers of them as occurs in many
New-Keynesian models). This causes permanence
and path dependence at a deep level. It should be noted that in more
sophisticated models, with more variables and more play operators, the
set of possible equilibria may be extremely complicated with the possibility of
`cascades' where one play operator starting to drag causes  others
to do so (the analogy with earthquakes made in the Introduction then
becomes even closer).

Secondly, the existence of different modes depending upon whether
particular play operators are currently `stuck' or `dragged' --- in
our case the `inner'
and `outer' modes. If some modes are less stable than others 
(in our main model the outer mode is 
less stable than the
inner one)
then a large enough shock may move the system far enough away from the
set of equilibria that the route back to an equilibrium is both long
and unpredictable. It may even move the system into an unstable regime
--- in this case runaway inflation --- without any change in the
system parameters. 
Both these features are highly significant not just because they
correspond closely to actual economic events but they have
implications for forecasting and policy prescriptions too. 

Our choice of inflation expectations as the candidate for an initial
investigation was influenced by the work of De Grauwe \cite{DeGrauwe2012} on
a different type of boundedly rational expectation formation process
in a simple DSGE model.  However, play operators are also a viable
candidate for modeling other sticky economic variables at both the
micro- and macro-economic levels. To demonstrate this, in our final
model we used one to represent sticky responses by the Central Bank.

The modeling approach presented above can be thought of as a `stress
test' of the usual rationality assumption in the underlying toy model.
Or to put it another way, it is examining the robustness of modeling
assumptions rather than just the stability of the solutions within a
particular model. As such, we believe that the introduction of a new
form of plausible stickiness has intrinsic merit not just as a form of
expectation formation. It provides an additional class of perturbed
models --- ones that are genuinely nonlinear, tractable, and capable
of changing solutions (and potentially policy prescriptions) in a way
that merely changing the parameters of an equilibrium model cannot.

Our second and third models demonstrated that there are various ways
in which this work can be extended, in particular to systems with
multiple agents and multiple play operators. Although it has not been
relevant to this paper play and stop operators, when combined appropriately
\cite{PRL1} can have a remarkably simple aggregated response, even when
connected via a network. This allows for (almost)-analytic solutions
even when cascades and rapid transitions between states are occurring
and will be the subject of future work.

\section*{Appendix}

\subsection*{A.\ \ Derivation of equations \eqref{eqn:SA1}, \eqref{eqn:SA1'}}

Here we show how to obtain equations \eqref{eqn:SA1}, \eqref{eqn:SA1'} from model \eqref{eqn:M1}--\eqref{formula'}. 
To this end, we substitute the equation for $r_t$ into the equation for  $y_t$ and obtain
\begin{equation*}\label{eqn:RM1a}
(1+ac_2)y_t=y_{t-1}-ac_1x_t+ap_t+\epsilon_t.
\end{equation*}
Next, we substitute this equation into the equation for $x_t$ and simplify to obtain
\begin{equation}\label{eqn:RM1}
\gamma x_t-\beta p_t=b_2y_{t-1} +(1-b_1)(1+ac_2)x_{t-1}+ b_2\epsilon_t+(1+ac_2)\eta_t,
\end{equation}
where
\[\gamma=1+ac_2+ab_2c_1, \qquad  \beta=b_1(1+ac_2)+ab_2.\]
Since
$p_t=x_t-s_t$,
equation (\ref{eqn:RM1}) can be rewritten as
\begin{equation}\label{eqn:RM2}
\alpha x_t+s_t=f_t
\end{equation}
with $\alpha$ and $f_t$ defined by \eqref{alpha}, \eqref{f}.
Therefore, $x_t=\alpha^{-1}(f_t-s_t)$, which combined with \eqref{alpha}, \eqref{f} gives
\begin{equation}
x_t= \frac{b_2 }{\alpha\beta}y_{t-1} +\frac{(1-b_1)(1+ac_2)}{\alpha\beta}x_{t-1}-\frac{1}{\alpha }s_t 
+\frac{b_2}{\alpha\beta}\epsilon_t+\frac{ 1+ac_2}{\alpha\beta}\eta_t.\label{eqn:RMPt}
\end{equation}
Subsequently, substituting equation \eqref{eqn:RMPt} into  equation \eqref{eqn:RM1a} gives
\begin{eqnarray}
y_t&=&\frac{ab_2(1-c_1)+\alpha\beta}{\alpha\beta(1+ac_2)}y_{t-1}+\frac{a(1-c_1)(1-b_1)}{\alpha\beta}x_{t-1}\nonumber\\
&&+\frac{a(c_1-1-\alpha)}{\alpha(1+ac_2)}s_t+\frac{\alpha\beta+ab_2(1-c_1)}{\alpha\beta(1+ac_2)}\epsilon_t
+\frac{a(1-c_1)}{\alpha\beta}\eta_t.\label{eqn:RMYt}
\end{eqnarray}
Equations \eqref{eqn:RMPt}, \eqref{eqn:RMYt} can be written as system \eqref{eqn:SA1} with the matrices $A$, $N$ and the vector $d$ defined by formulas \eqref{AA}.


Equation \eqref{eqn:SA1'} can be obtained from relation
\eqref{eqn:RM2} using the inversion formula for the play
operator. This inversion formula is presented for a more general
Prandtl-Ishlinskii (PI) operator, including the play operator as a particular case, in Appendix D. 

\subsection*{B.\ \ Local stability analysis}

The characteristic polynomial of matrix $B$ is 
\[P_B(\lambda)=\lambda^2-\lambda\left(\frac{2+ac_2-b_1(1+ac_2)}{1+a(b_2c_1+c_2)}\right)+\frac{1-b_1}{1+a(b_2c_1+c_2)}.
\]
Applying Jury's stability criterion to the characteristic polynomial gives the following set of inequalities:
\begin{eqnarray}
\hspace{-5mm} P_B(1)&=&1-\frac{2+ac_2-b_1(1+ac_2)}{1+a(b_2c_1+c_2)}+\frac{1-b_1}{1+a(b_2c_1+c_2)}>0,\nonumber\\
\hspace{-5mm} P_B(-1)&=&1+\frac{2+ac_2-b_1(1+ac_2)}{1+a(b_2c_1+c_2)}+\frac{1-b_1}{1+a(b_2c_1+c_2)}>0,\nonumber\\
1&>&\frac{1-b_1}{1+a(b_2c_1+c_2)}.\nonumber
\end{eqnarray}
It is easy to see that all the three inequalities above are satisfied for 
any set of parameters
$
a, b_2,c_1,c_2>0
$
and $0< b_1<1$,
hence every equilibrium $z_*(s_*)$ with $|s_*|<\rho$ is locally stable.

Now, let us consider the system without stiction.
The characteristic polynomial of matrix $A$ is 
\[P_A(\lambda)=\Delta \lambda^2- {(1-b_1)(2+ac_2)}\lambda+{1-b_1}\]
with $\Delta$ defined by \eqref{Delta}.
Applying Jury's stability criterion, we obtain
\begin{eqnarray}
\hspace{-5mm} P_A(1)&=&1-\frac{(1-b_1)(2+ac_2)}{\Delta}+\frac{1-b_1}{\Delta}>0,\nonumber\\
\hspace{-5mm} P_A(-1)&=&1+\frac{(1-b_1)(2+ac_2)}{\Delta}+\frac{1-b_1}{\Delta}>0,\nonumber\\
1&>&\frac{1-b_1}\Delta.\nonumber
\end{eqnarray}
Taking into account the constraints $
a, b_2,c_1,c_2>0
$
and $0< b_1<1$,  these conditions result in the relationship
\[c_1>1.\]

Note that the system $z_{t}=Az_{t-1}$ is the linearization of sticky system \eqref{eqn:SA1} 
at infinity, hence it describes the return of the sticky system
towards near equilibrium dynamics after a large perturbation. Thus, the stability condition $c_1>1$ for $A$ agrees with the global stability criterion obtained in Section \ref{staglo}.



\subsection*{C.\ \ The effect of parameters on stability properties}

Here we provide some numerical analysis concerning the effect of the parameters
on stability properties of the equilibrium states.
Stronger stability generally implies lower volatility and
more infrequent transitions between different equilibrium states. We quantify
local stability using the maximum absolute value, $|\lambda_{i,e}|$,
of eigenvalues of the linearized system at an equilibrium point.  
The
subscripts $e$ and $i$ refer to the system without stickiness
($\rho=0$) and with stickiness ($\rho=1$), respectively.

\begin{figure}[h!]
\centering
\includegraphics[scale=0.2]{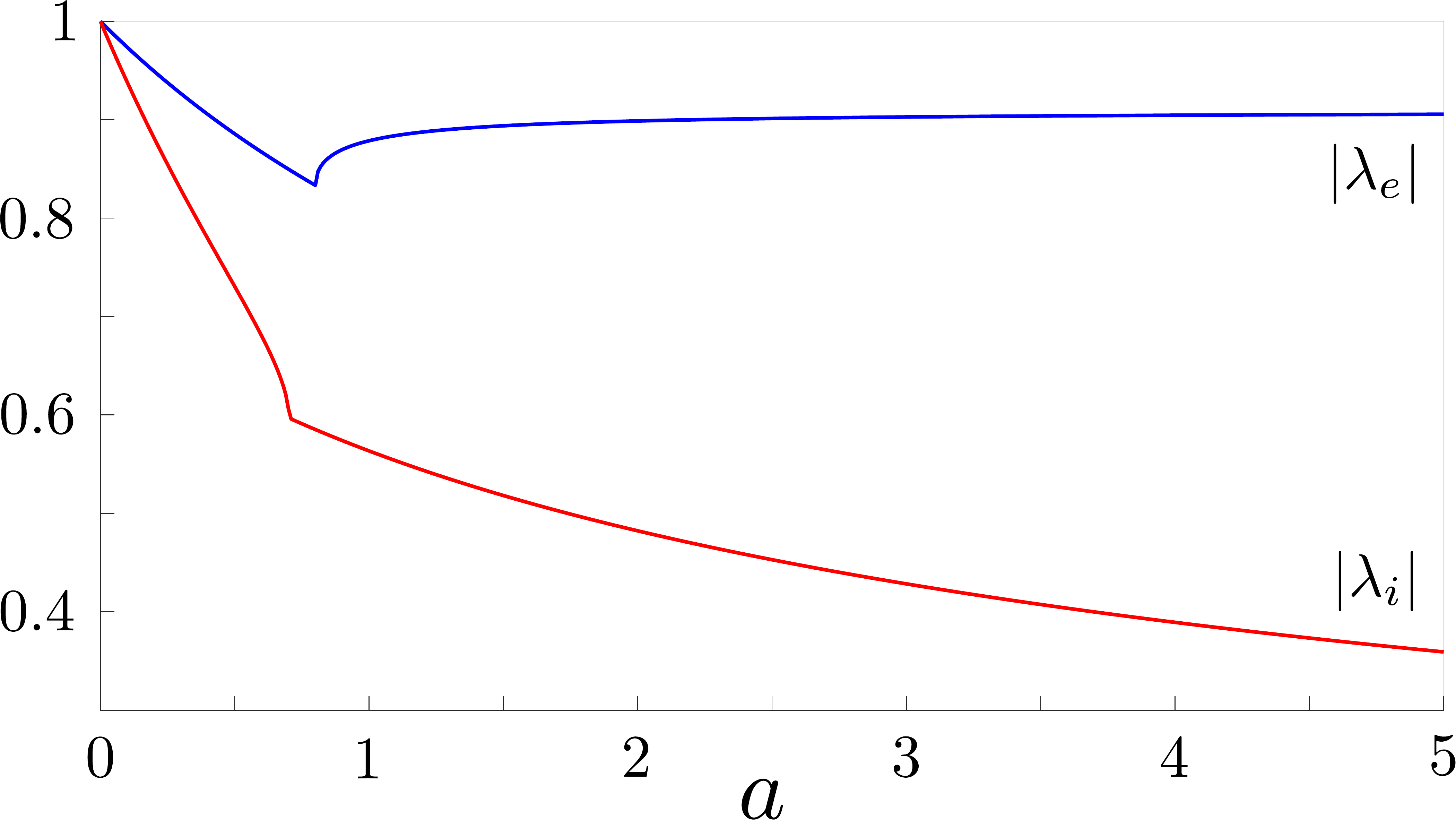}
\caption{Variation of $|\lambda_i|$  and $|\lambda_e|$ with $a$.
	Other parameters are taken from Table \ref{tab:DS1}.}\label{fig:SALA}
\end{figure}

The model contains five other parameters, $a$, $b_1$, $b_2$, $c_1$ and $c_2$.
Fig.~\ref{fig:SALA} shows the dependence of $|\lambda_{i,e}|$ on the
parameter $a$ and 
implies that 
the system with stickiness is more stable than the system without stickiness.
Other parameter values are taken from Table \ref{tab:DS1}.
Interestingly, the system with stickiness becomes more stable for increasing $a$, while
this dependence for the non-sticky system is non-monotone since $|\lambda_e|$ has a minimum at $a\approx 0.8$.

\begin{figure}[h!]
\centering
\begin{subfigure}{0.45\textwidth}
\includegraphics[width=\textwidth]{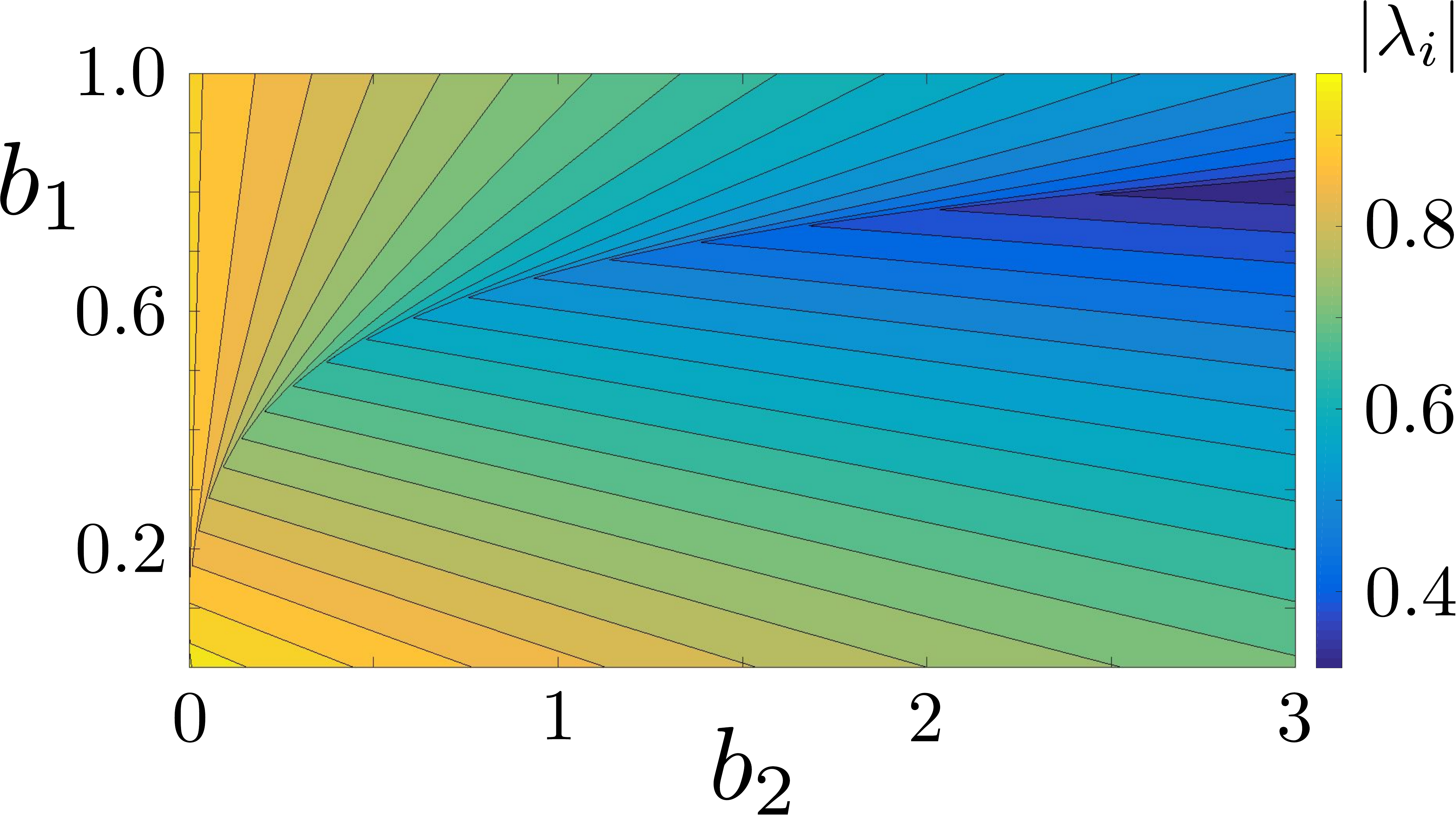}
\caption{}
\label{fig:SALIntB}
\end{subfigure} 
\quad
\begin{subfigure}{0.45\textwidth}
\includegraphics[width=\textwidth]{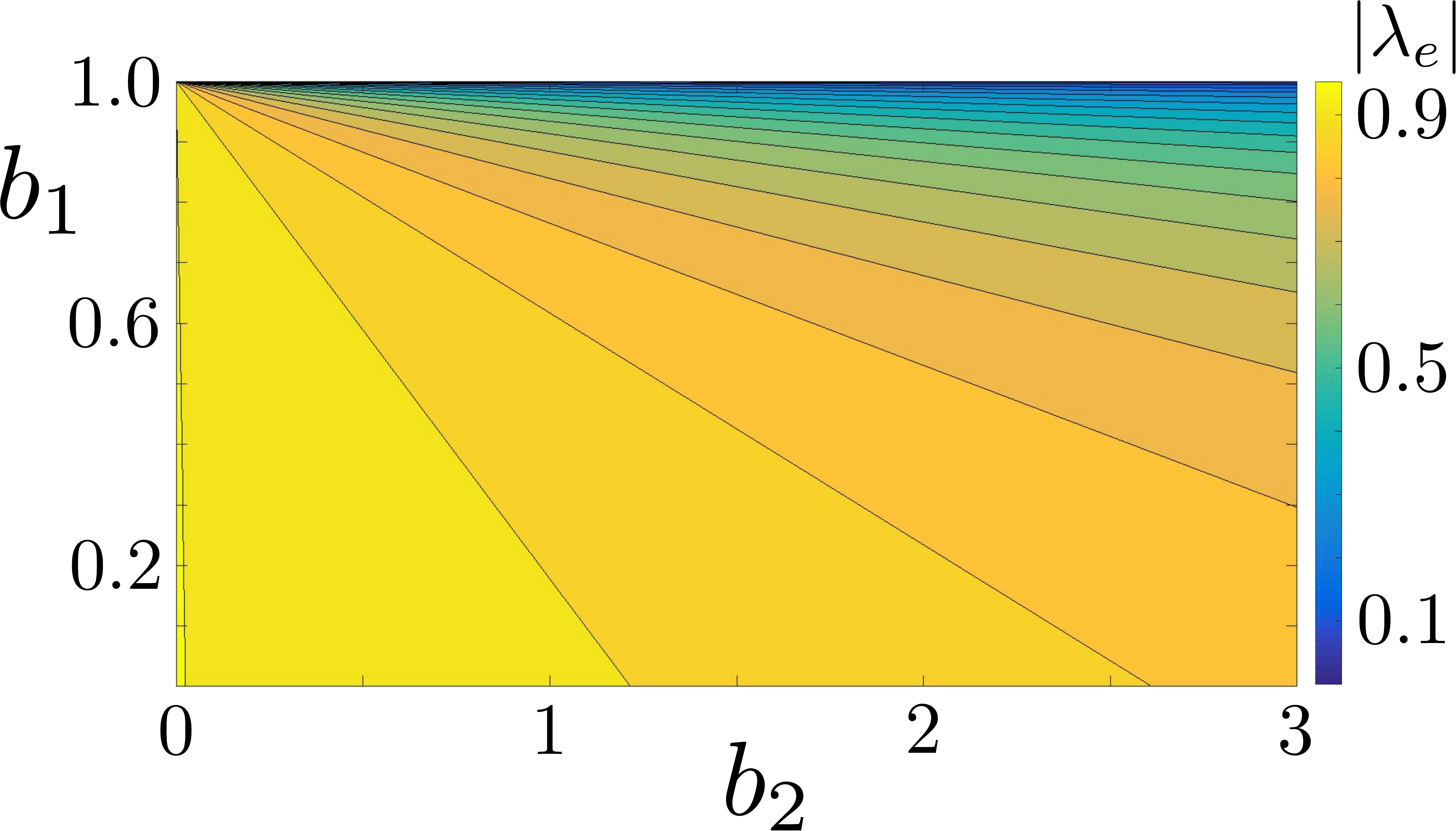}
\caption{}
\label{fig:SALExtB}
\end{subfigure}
\caption{Dependence of (a) $|\lambda_i|$ and (b) $|\lambda_e|$ on $b_1$ and $b_2$.
	Other parameters are taken from Table \ref{tab:DS1}. 
	}\label{fig:SALB}
\end{figure}

The range of output gap equilibrium values is proportional to the
ratio of parameters $b_1$ and $b_2$ according to
\eqref{eqn:SA1E}. Fig.~\ref{fig:SALB} presents the dependence of
$|\lambda_{i,e}|$ on these parameters. The sticky system is more
stable than its non-sticky counterpart for $b_1<0.9$, but becomes less
stable than the non-sticky system as $b_1$ approaches 1 (in the latter
case, the future inflation rate is defined predominantly by
expectations). The dependence of $|\lambda_{i,e}|$ on $b_2$ and the
dependence of $|\lambda_e|$ on $b_1$ is monotone (stronger stability
for larger $b_{1,2}$), while the dependence of $|\lambda_i|$ on $b_2$
is non-monotone. The strongest stability is achieved by the sticky
system for some intermediate value of $b_1$ between 0 and 1.

\begin{figure}[h!]
	\centering
	\begin{subfigure}{0.45\textwidth}
		\includegraphics[width=\textwidth]{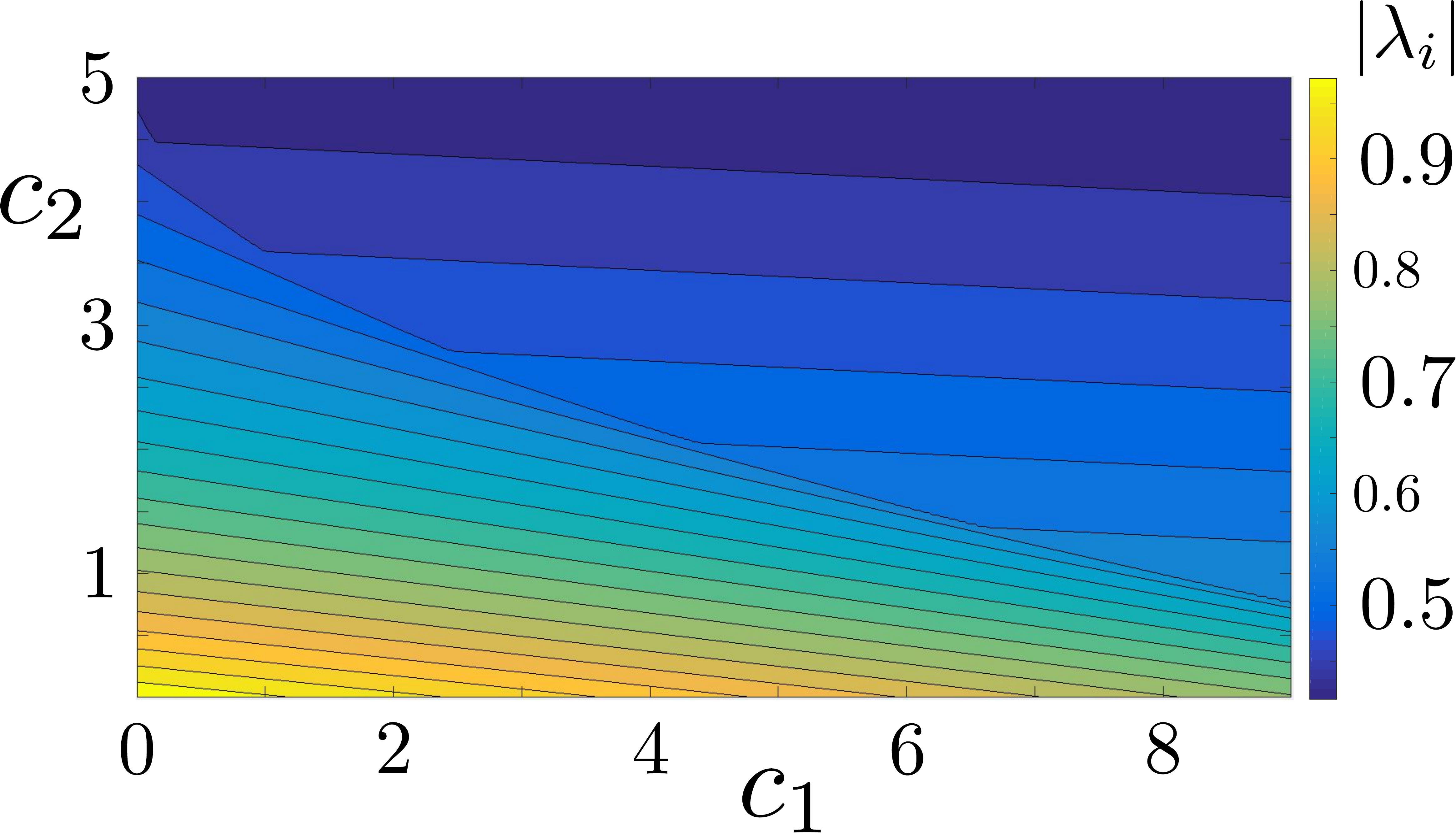}
		\caption{} 
		\label{fig:SALIntC}
	\end{subfigure} 
	\quad
	\begin{subfigure}{0.45\textwidth}
		\includegraphics[width=\textwidth]{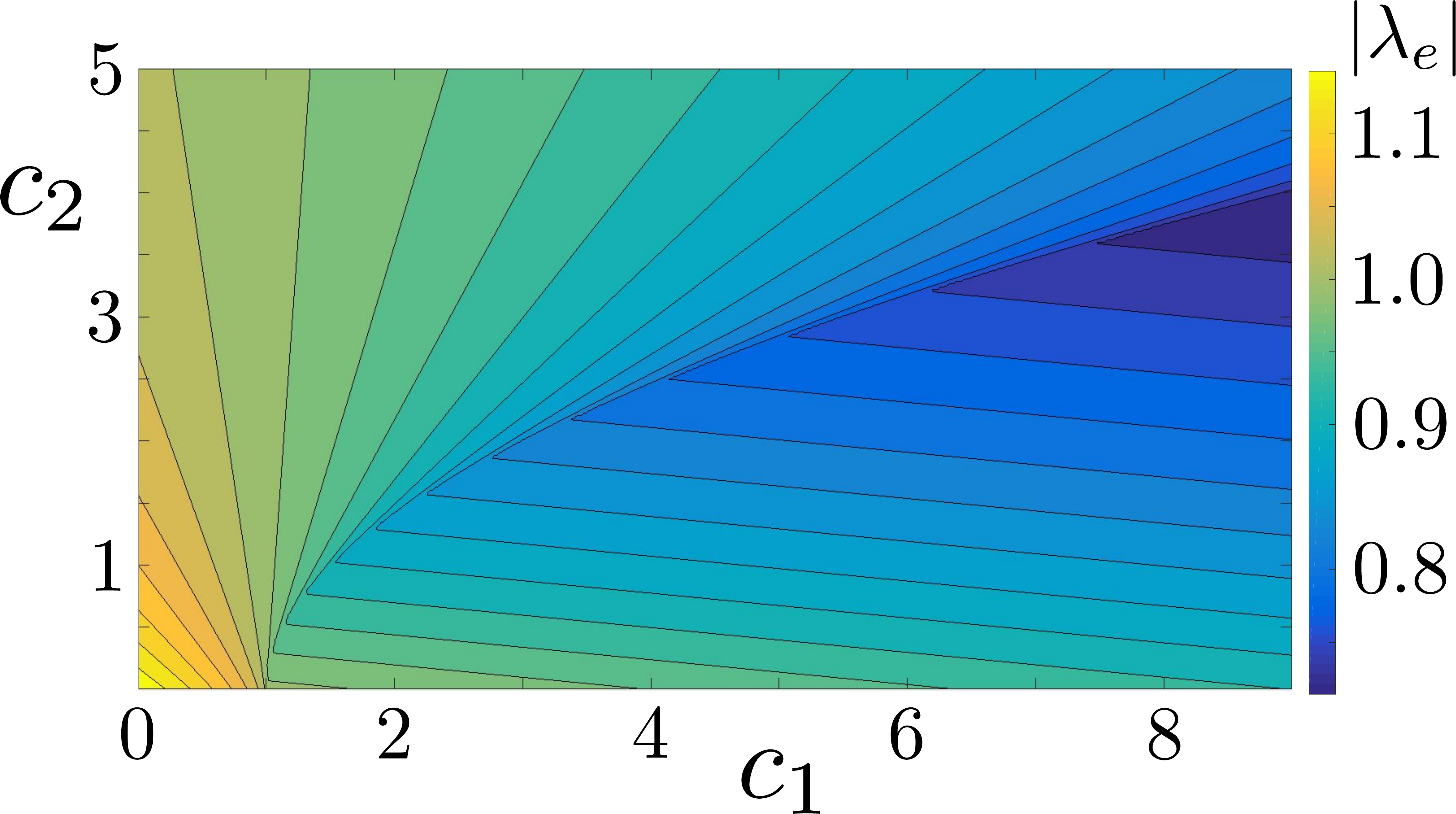}
		\caption{}
		\label{fig:SALExtC}
	\end{subfigure}
	\caption{Dependence of (a) $|\lambda_i|$ and (b) $|\lambda_e|$ on $c_1$ and $c_2$.
		Other parameters are taken from Table \ref{tab:DS1}. 
		 }\label{fig:SALC}
\end{figure}

\begin{figure}[h!]
	\begin{subfigure}{0.45\textwidth}
		\includegraphics[scale=0.2]{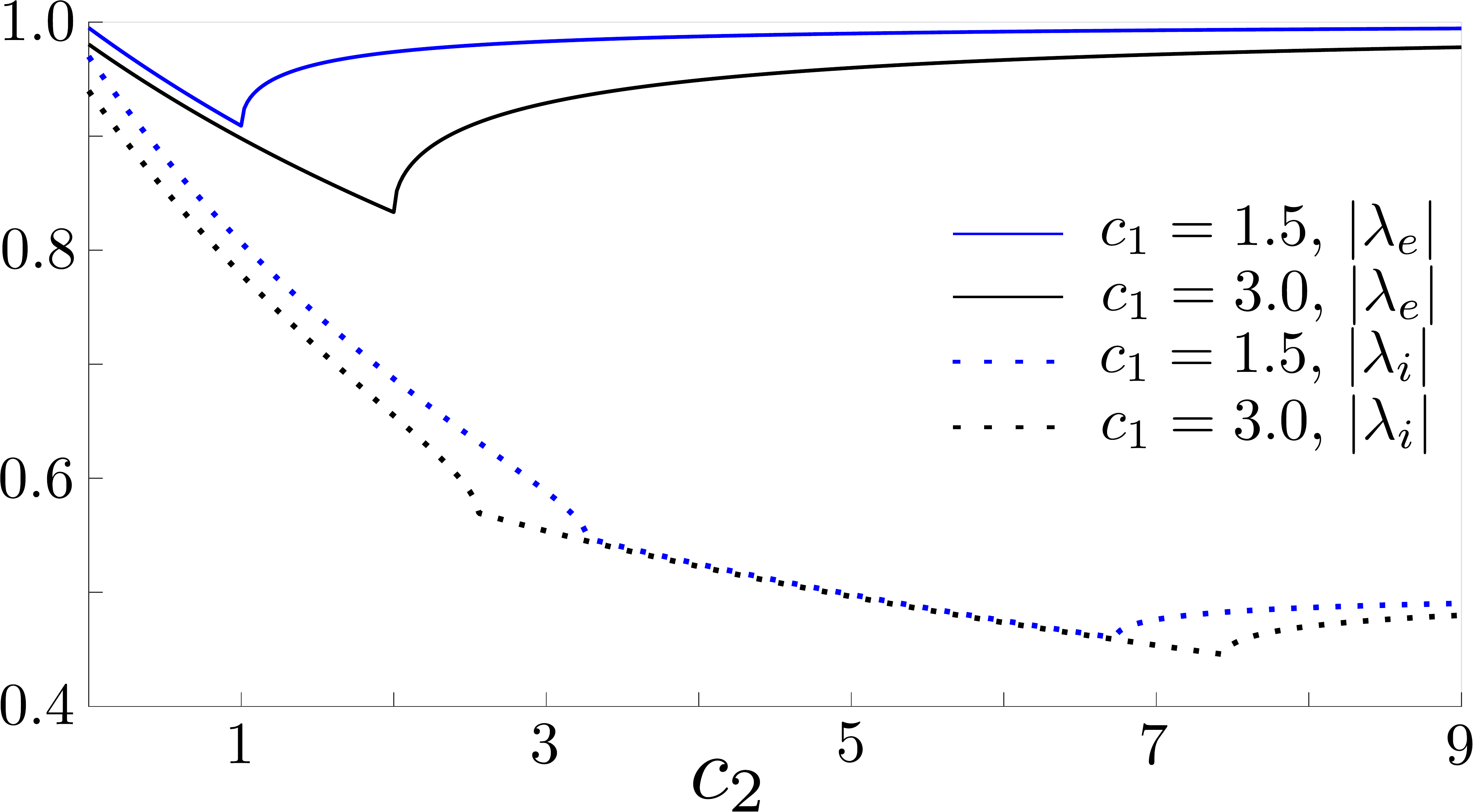}
		\caption{}
		\label{fig:SALIntC1}
	\end{subfigure}
	%
	%
%
%
%


	\begin{subfigure}{0.45\textwidth}
		\includegraphics[scale=0.2]{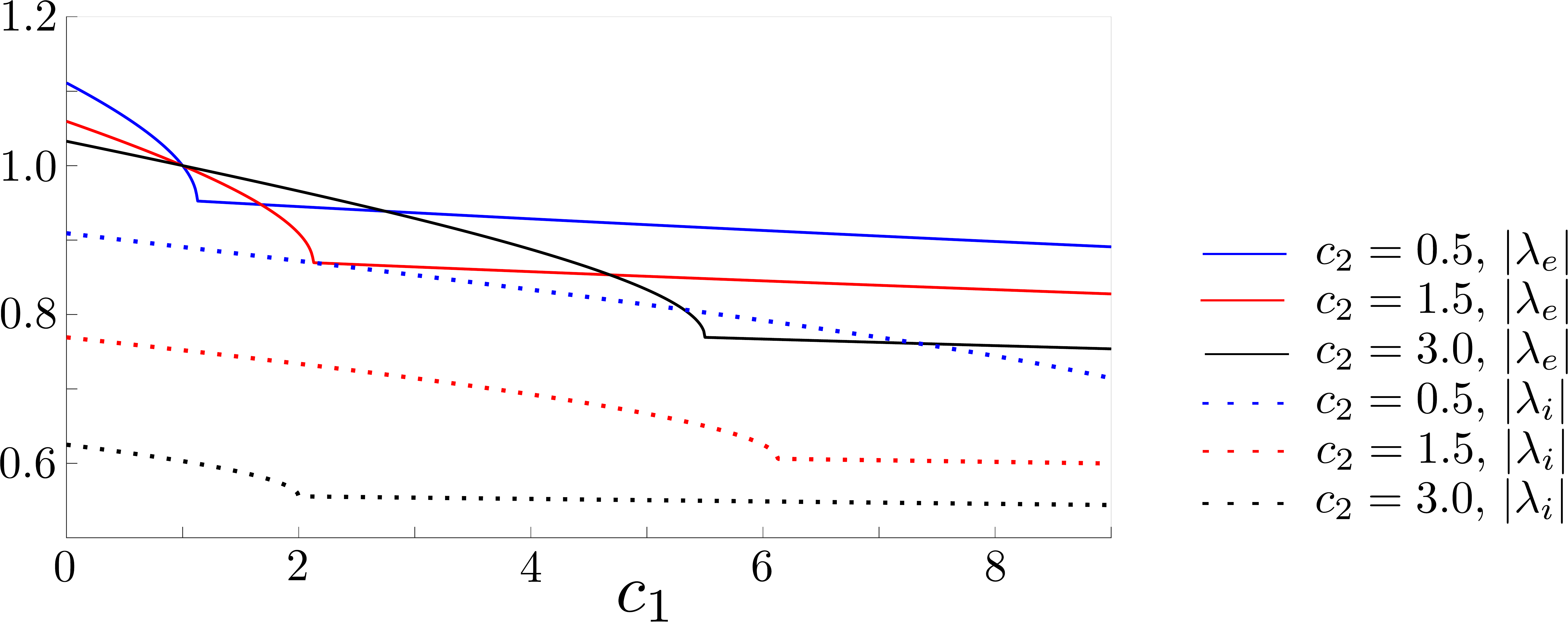}
		\caption{}
		\label{fig:SALIntC2}
	\end{subfigure}


%
%
%
	\caption{Crossections of the plots shown in Fig.~\ref{fig:SALC} (a) for various $c_2$ values and  (b) for various $c_1$ values. 
				}\label{fig:SALC1C2}
\end{figure}

Parameters $c_1$ and $c_2$ control the range of inflation rate equilibrium values
according to \eqref{eqn:SA1E}. This range
contracts when $c_1$ increases (for $c_1>1$) and expands when $c_2$ increases.
Fig.~\ref{fig:SALC} shows that the sticky system is generally more stable than the non-sticky one.
Both systems become more stable with increasing $c_1$ (stronger inflation targeting in Taylor's rule),
see Figs.~\ref{fig:SALC}(a, b) and \ref{fig:SALC1C2}(a, b). The dependence of  $|\lambda_{i}|$ on $c_2$ demonstrates
some slight non-monotonicity for large $c_2$ values, see Figure \ref{fig:SALC1C2}(b).
The non-monotonicity of $|\lambda_{i}|$  with $c_2$ is much more pronounced 
with the minimum achieved for a certain value of $c_2$ depending on $c_1$, see Figs.~\ref{fig:SALC}(b)
and \ref{fig:SALC1C2}(b). This minimum corresponds to the strongest stability and, in this sense, 
optimizes the Central Bank policy. In Fig.~\ref{fig:SALC}(b), the strongest stability 
is achieved on the `parabolic' line.
 

\subsection*{D.\ \ Inversion of the PI operator}

In this section, we consider the inversion of the PI operator, which
is necessary to transform the implicit system \eqref{eqn:M1}, \eqref{eqn:M1'}
coupled with relation \eqref{eqn:MStop} into the explicit form
\eqref{eqn:SA111}. Here we use the term `PI operator' for an
input-output relationship of the form
\begin{equation}\label{PIpi}
f_t=\alpha x_t + \sum_{i=1}^n\mu_i\mathscr{S}_{\rho_i}[x_t],
\end{equation}
where the weights $\mu_i$ are allowed to have any sign, $\alpha\ge 0$, and $\rho_1<\rho_2<\cdots<\rho_n$.
Such an operator is completely defined by the so-called {\em Primary Response} (PR) {\em function} $\phi(x)$,
which describes the output in response to a monotonically increasing input.
\begin{figure}[ht]
\centering
  \includegraphics[scale=0.3]{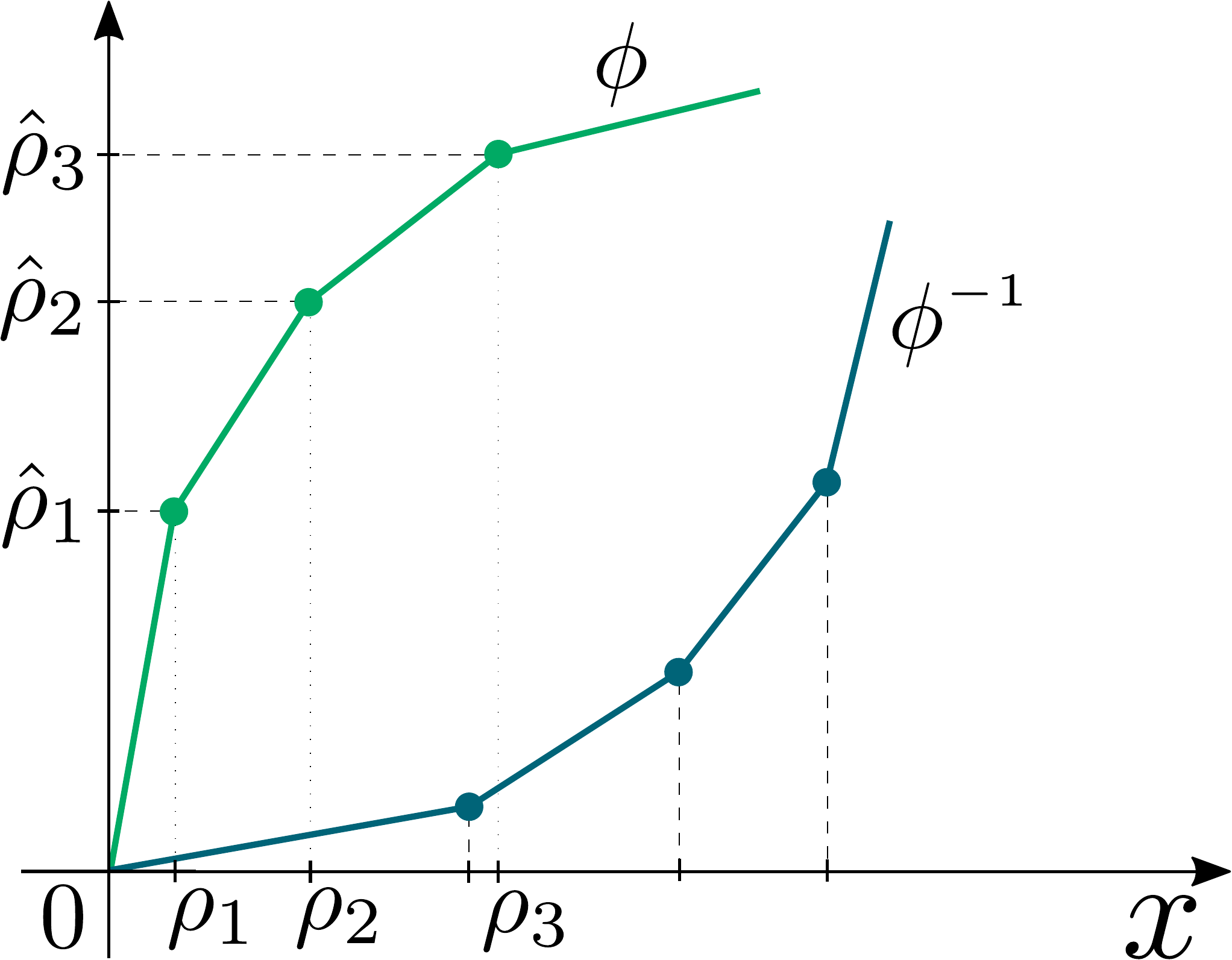}
\caption{PR  function $\phi$ of PI operator \eqref{PIpi}
	and PR function $\phi^{-1}$ of its inverse PI operator \eqref{PIpi'}.}\label{inversion}
\end{figure}
Here, this is a piecewise linear continuous function satisfying $\phi(0)=0$ with the slopes defined by
\[
\phi'(x)=\left\{ 
\begin{array}{ll}
\alpha+\mu_n+\cdots+\mu_2+\mu_1, & 0<x<\rho_1,\\
\alpha+\mu_n+\cdots+\mu_2, & \rho_1<x<\rho_2,\\
\vdots\\
\alpha+\mu_n,& \rho_{n-1}<x<\rho_{n},\\
\alpha,& x>\rho_n,
\end{array}
\right.
\]
see Fig.~\ref{inversion}. As shown in \cite{krejvci1986hysteresis}, if the slopes of $\phi$ are all positive, then 
the PI operator \eqref{PIpi} is invertible, and the inverse relationship is also a PI operator:
\begin{equation}\label{PIpi'}
x_t=\hat\alpha f_t + \sum_{i=1}^n\hat\mu_i\mathscr{S}_{\hat\rho_i}[f_t].
\end{equation}
Further, the PR function of operator \eqref{PIpi'} is the inverse of the PR function $\phi$ of operator \eqref{PIpi}.
This allows one to express the weights $\hat\alpha, \hat{\mu}_i$ and the thresholds $\hat{\rho}_i$
explicitly in terms of the  weights $\alpha, {\mu}_i$ and the thresholds ${\rho}_i$.
In particular, the equation $\alpha x_t + s_t=f_t$ with $s_t= {\mathcal S}_\rho[x_t]$ (see \eqref{eqn:RM2})
can be inverted as
\[
x_t=\frac1\alpha f_t-\frac1{\alpha(1+\alpha)}{\mathcal S}_{(1+\alpha)\rho}[f_t],
\]
and this implies $s_t=\frac1{1+\alpha}{\mathcal S}_{(1+\alpha)\rho}[f_t]$, which is 
equivalent to \eqref{eqn:SA1'} (cf.~Appendix A).

\subsection*{E.\ \ Sticky Taylor rule}
In order to convert system \eqref{M1'}, \eqref{M1} to the explicit form, we replace the variable $y_t$ with the variable $g_t=c_1x_t+c_2y_t$ and  obtain
\begin{align}
g_t=& \hspace{1mm} (c_1+ac_2)x_t+g_{t-1}-c_1x_{t-1}-ac_2\mathscr{P}_\sigma[g_t]+c_2\epsilon_t,\label{eqnXT}\\
x_t=& \hspace{1mm} \frac{c_2(1-b_1)}{b_2c_1+c_2(1-b_1)}x_{t-1}+\frac{b_2}{b_2c_1+c_2(1-b_1)}g_t+\frac{c_2(1-b_1)}{b_2c_1+c_2(1-b_1)}\eta_t.\label{eqnPT}
\end{align}
Further, substituting \eqref{eqnPT} into \eqref{eqnXT} gives
\begin{equation}\label{eqn:InvertPlaywFt}
\alpha g_t+\kappa\mathscr{P}_{\sigma}[g_t]=f_t
\end{equation}
with
\[
\alpha= \hspace{1mm} \frac{c_2(1-b_1-ab_2)}{b_2c_1+c_2(1-b_1)},\qquad
\kappa= \hspace{1mm} ac_2,
\]
\[
f_t= g_{t-1}-c_1x_{t-1}+\frac{c_2(1-b_1)(c_1+ac_2)}{b_2c_1+c_2(1-b_1)}\left(x_{t-1}+\eta_t\right)+c_2\epsilon_t.
\]
Using that $\alpha>0$ due to \eqref{condiinv}, we can invert
\eqref{eqn:InvertPlaywFt} as in Appendix D to obtain
\[ g_t=\frac{1}{\alpha}\left(f_t-\frac{\kappa}{\alpha+\kappa}\mathscr{P}_{\alpha\sigma}[f_t]\right).\]
This equation together with \eqref{eqnPT} defines the explicit system for \eqref{M1'}, \eqref{M1}.
The linearization $z_t=Bz_{t-1}$ of this system at any equilibrium point with $s_*\in(-\sigma,\sigma)$ 
has the matrix
\[
B=\frac1{1-b_1-ab_2}\begin{pmatrix}
{1-b_1}&{a(1-b_1)}\\
{b_2}&{1-b_1}
\end{pmatrix}.
\]
Since 
\[
\det B =\frac{1-b_1}{1-b_1-ab_2}>1,
\]
all these equilibrium states are unstable.

\section*{Acknowledgments}
DR and EK acknowledge the support of NSF through grant DMS-1413223.\\
PK acknowledges the support of GA\v CR through grant GA15-12227S.

\section*{References}

\bibliographystyle{elsarticle-num}

\bibliography{kkpr_bibJun26}

\end{document}